\documentclass[11pt]{article}
\usepackage[utf8]{inputenc}
\usepackage{fullpage,amsmath,amssymb,bm,url,mathrsfs}
\usepackage{algorithmic}
\usepackage{algorithm}
\usepackage{amsthm}
\usepackage[bookmarks]{hyperref}
\usepackage{makecell}
\usepackage{tablefootnote}
\usepackage{comment}
\usepackage{graphics,graphicx}
\usepackage{subcaption,caption,cleveref}
\hypersetup{
    colorlinks=true,
    linkcolor=blue,
    citecolor=blue,
    filecolor=magenta,      
    urlcolor=cyan,
    bookmarksopen=true
}
\usepackage{multirow}
\usepackage[round,colon,authoryear]{natbib}
\usepackage{xcolor}
\usepackage{dsfont}
\usepackage{bbm}
\usepackage{enumitem}

\def\calF{{\mathcal F}}

\def\calM{{\mathcal M}}
\def\calN{{\mathcal N}}

\def\calQ{{\mathcal Q}}
\def\calR{{\mathcal R}}
\def\calS{{\mathcal S}}

\def\calX{{\mathcal X}}

\def\bcalE{{\boldsymbol{\mathcal E}}}

\def\bcalG{{\boldsymbol{\mathcal G}}}

\def\bcalM{{\boldsymbol{\mathcal M}}}

\def\bcalS{{\boldsymbol{\mathcal S}}}

\def\bcalX{{\boldsymbol{\mathcal X}}}

\def\bfrakM{{\boldsymbol{\mathfrak{M}}}}

\def\EE{{\mathbb E}}

\def\II{{\mathbb I}}

\def\OO{{\mathbb O}}
\def\PP{{\mathbb P}}
\def\QQ{{\mathbb Q}}
\def\RR{{\mathbb R}}

\def\s{{\boldsymbol s}}

\DeclareMathOperator*{\argmin}{arg\,min}

\def\bs{\boldsymbol{s}}

\def\calF{{\cal  F}}

\def\calM{{\cal  M}} \def\scrM{{\mathscr  M}}
\def\calN{{\cal  N}}

\def\calQ{{\cal  Q}} 
\def\calR{{\cal  R}} 
\def\calS{{\cal  S}}

\def\calX{{\cal  X}}

\newcommand{\bfm}[1]{\ensuremath{\mathbf{#1}}}

   \def\bB{\bfm B}

\def\be{\bfm e}   \def\bE{\bfm E}  \def\EE{\mathbb{E}}
    
   \def\bG{\bfm G}  
     
   \def\bI{\bfm I}  \def\II{\mathbb{I}}

\def\bm{\bfm m}   \def\bM{\bfm M}  
     
     \def\OO{\mathbb{O}}
     \def\PP{\mathbb{P}}
     \def\QQ{\mathbb{Q}}
     \def\RR{\mathbb{R}}
\def\bs{\bfm s}   \def\bS{\bfm S}  
     
\def\bu{\bfm u}   \def\bU{\bfm U}  
\def\bv{\bfm v}   \def\bV{\bfm V}  
   \def\bW{\bfm W}  
   \def\bX{\bfm X}  
   \def\bY{\bfm Y}  
   \def\bZ{\bfm Z}

\def\bSigma{\bfm \Sigma}

\def\bXi{\bfm \Xi}

\def\bEsai{\bar{\bE}^\ast_{a,i}}
\def\bEsani{\bar{\bE}^\ast_{a,-i}}
\def\bEsa{\bar{\bE}^\ast_a}
\def\bEss{\bar{\bE}^\ast_2}

\def\bEa{\bar{\bE}_a}
\def\bEs{\bar{\bE}_2}

\def\hat{\widehat}
\def\wt{\widetilde}

\newcommand\ind[1]{\II\left( #1 \right)}
\newcommand\inp[2]{\left\langle #1, #2 \right\rangle}
\newcommand\fro[1]{\left\| #1 \right\|_{\rm F}}
\newcommand\op[1]{\left\|#1\right\|}

\newcommand{\mat}[1]{\begin{bmatrix}#1 \\ \end{bmatrix}}

\def\bcalG{\boldsymbol{\mathcal{G}}}
\def\bcalE{\boldsymbol{\mathcal{E}}}

\newtheorem{theorem}{Theorem}

\newtheorem{definition}{Definition}

\newtheorem{assumption}{Assumption}
\newtheorem{lemma}{Lemma}
\newtheorem{corollary}{Corollary}

\newtheorem{conjecture}{Conjecture}
\theoremstyle{remark}

\newcommand{\E}{\mathbb{E}}
\newcommand{\Prob}{\mathbb{P}}

\newcommand{\cM}{\mathcal{M}}

\newcommand\tTop{\rule{0pt}{4ex}}       
\newcommand\tBot{\rule[-2ex]{0pt}{0pt}} 

\newcommand\supt{^{(t)}}
\newcommand\suptt{^{(t-1)}}

\newcommand\sigmasg{\sigma_{\textsf{sg}}}
\newcommand\submax{\textsf{\tiny max}}
\newcommand\submin{\textsf{\tiny min}}

\begin{document}

\title{Optimal Clustering by Lloyd's Algorithm for Low-Rank Mixture Model}

\author{Zhongyuan Lyu and Dong Xia\footnote{Dong Xia's research was partially supported by Hong Kong RGC Grant GRF 16300121 and GRF 16301622.}\\
{\small Department of Mathematics, Hong Kong University of Science and Technology}}

\date{(\today)}

\maketitle

\begin{abstract}
This paper investigates the computational and statistical limit in clustering matrix-valued observations. We propose a low-rank mixture model (LrMM), adapted from the classical Gaussian mixture model (GMM) to treat matrix-valued observations, which assumes low-rankness for population center matrices. A computationally efficient clustering method is designed by integrating Lloyd's algorithm and low-rank approximation. Once well-initialized, the algorithm converges fast and achieves an exponential-type clustering error rate that is minimax optimal. Meanwhile, we show that a tensor-based spectral method delivers a good initial clustering. Comparable to GMM, the minimax optimal clustering error rate is decided by the {\it separation strength}, i.e, the minimal distance between population center matrices. By exploiting low-rankness, the proposed algorithm is blessed with a weaker requirement on the separation strength.  Unlike GMM, however, the computational difficulty of LrMM is characterized by the {\it signal strength}, i.e, the smallest non-zero singular values of population center matrices.  
Evidences are provided showing that no polynomial-time algorithm is consistent if the signal strength is not strong enough, even though the separation strength is strong.  Intriguing differences between estimation and clustering under LrMM are discussed. The merits of low-rank Lloyd's algorithm are confirmed by comprehensive simulation experiments. Finally, our method outperforms others in the literature on real-world datasets. 
\end{abstract}

\section{Introduction}\label{sec:intro}
Nowadays, clustering {\it matrix-valued} observations becomes a ubiquitous task in diverse fields. For instance, each highly variable region (HVR) in the var genes of human malaria parasite \citep{larremore2013network,jing2021community}  is representable by an adjacency matrix and a key scientific question is to identify structurally-similar HVRs by, say, clustering the associated adjacency matrices. The international trade flow of a commodity across different countries can be viewed as a weighted adjacency matrix \citep{lyu2021latent,cai2022generalized}. Finding the similarity between the trading patterns of different commodities is of great value in understanding the global economic structure. This can also be achieved by clustering the weighted adjacency matrices. Other notable examples include clustering multi-layer social networks \citep{dong2012clustering,han2015consistent} and multi-view data \citep{kumar2011co,mai2021doubly}, modeling the connectivity of brain networks \citep{arroyo2021inference,sun2019dynamic}, clustering the correlation networks between bacterial species \citep{stanley2016clustering},  and EEG data analysis \citep{gao2021regularized}, etc. 

Since matrix-valued observations can always be vectorized, a naive approach is to ignore the matrix structure so that numerous classical clustering algorithms, e.g. K-means or spectral clustering, are readily applicable. However, matrix observations are usually blessed with hidden low-dimensional structures, among which low-rankness is perhaps the most common and explored. Network models such as {\it stochastic block model} \citep{holland1983stochastic, jing2021community}, {\it random dot product graph} \citep{athreya2017statistical} and {\it latent space model} \citep{hoff2002latent} often assume a low-rank expectation of adjacency matrix. Low-rank structures have also been successfully explored in brain image clustering \citep{sun2019dynamic}, EEG data analysis \citep{gao2021regularized}, and international trade flow data \citep{lyu2021latent}, to name but a few. Table~\ref{tb:data-summary} presents a summary of datasets analyzed in our paper, where the matrix ranks $r_k$'s (suggested by the numerical performance of our algorithm) are much smaller than the ambient dimensions $(d_1,d_2)$. Without loss of generality, we assume $d_1\geq d_2$.  For these applications, the naive clustering approach becomes statistically sub-optimal since the planted low-dimensional structure is overlooked. 

\begin{table}
\begin{center}
  \begin{tabular}{ c  c  c c c }
    \hline
    Dataset & $n$ & $(d_1,d_2)$ &  $K$ &Ranks \\ \hline
    BHL \citep{mai2021doubly} & 27 & (1124,4)& 3& $\sim\{1,1,1\}$\\
    EEG \citep{zhang1995event}& 122 & (256,64)& 2& $\sim\{2,1\}$\\
    Malaria gene networks  \citep{larremore2013network} & 9 & (212,212)& 6 & $\le 15$ \\ 
    UN trade flow networks \citep{lyu2021latent} & 97 & (48,48) & 2 & $\sim\{3,2\}$\\
    \hline
  \end{tabular}
\end{center}
\caption{Summary of datasets. Here, $n$ is the sample size, $(d_1, d_2)$ is the dimension of each matrix observation, and $K$ is number of clusters. The underlying rank $(r_k's)$ of population center matrices from different clusters can be unequal.}
\label{tb:data-summary}
\end{table}

Motivated by the aforementioned applications, throughout this paper, we assume that {\it each matrix-valued observation has a low-rank expectation} and the {\it expectations are equal for observations from the same cluster}. It is the essence of {\it low-rank mixture model} (LrMM), which shall be formally defined in Section~\ref{sec:method}. Several clustering methods exploiting low-rankness have emerged in the literature. \cite{sun2019dynamic} introduces a tensor Gaussian mixture model and recasts the clustering task as estimating the factors in low-rank tensor decomposition. K-means clustering is then applied to the estimated factors. While a sharp estimation error rate is derived under a suitable signal-to-noise ratio (SNR) condition, the accuracy of clustering is not provided. A tensor normal mixture model is proposed by \cite{mai2021doubly}, where the authors designed an enhanced EM algorithm for estimating the distributional parameters. Under appropriate conditions, sharp estimation error rates are established showing that minimax optimal {\it test} clustering error rate is attainable. However, the {\it training} clustering error is missing, and it is even unclear whether the proposed EM algorithm can consistently recover the true cluster memberships. Aimed at analyzing multi-layer networks, \cite{jing2021community} proposed a mixture multi-layer SBM where a spectral clustering method based on tensor decomposition is investigated. Clustering error rate is established under a fairly weak network sparsity condition, although the rate is likely sub-optimal. More recently, \cite{lyu2021latent} extended the mixture framework to latent space model and a sub-optimal clustering error rate was also derived. Note that \cite{jing2021community} and \cite{lyu2021latent} both require a rather restrictive condition in that $n=O(d_1)$ rendering their theories unattractive in many scenarios. Other representative works include \cite{chen2020global}, \cite{cai2021jointly}, \cite{gao2021regularized} and \cite{stanley2016clustering}, but clustering error rates were not studied.

Note that LrMM reduces to the famous {\it Gaussian mixture model} (GMM) in the dimension $d^{\ast}:=d_1d_2$ if each matrix-valued observation has a full-rank expectation, and the noise matrix has i.i.d. standard normal entries. Under GMM, \cite{loffler2021optimality} proved that a spectral method attains, with high probability,  an average mis-clustering error rate $\exp(-\Delta^2/8)$ that is optimal in the minimax sense.  Here $\Delta$ is the minimal Euclidean distance between the expected centers of distinct clusters (i.e., population center matrices), referred to as the {\it separation strength}. This exponential rate was established by \cite{loffler2021optimality} under a separation strength\footnote{For narration simplicity, we set the number of clusters $K=O(1)$ here.} condition $\Delta\gg 1+d^{\ast}n^{-1}$. \cite{gao2022iterative} investigated a more general iterative algorithm that achieves the same exponential rate under a weaker separation strength condition $\Delta\gg 1+(d^{\ast}/n)^{1/2}$. More recently, \cite{zhang2022leave} applied the leave-one-out method and proved the optimality of spectral clustering under a relaxed separation strength condition. Besides deriving the optimal clustering error rate, prior works also made efforts to  establish the phase transitions in exact recovery, i.e., when the clustering error is zero. \cite{ndaoud2018sharp} investigated a power iteration algorithm for a two-component GMM and proved that exact recovery is attained w.h.p.  if $\Delta^2$ is greater than $\big(1+(1+2d^{\ast}n^{-1}\log^{-1} n)^{1/2}\big)\cdot \log n$. In addition, the author showed that exact recovery is impossible if $\Delta^2$ is smaller than the aforesaid threshold. Later, \cite{chen2021cutoff} established a similar phase transition for general $K$-component GMM based on a semidefinite programming (SDP) relaxation. These foregoing works suggest an intriguing gap in the regime $n=O(d^{\ast})$: \cite{ndaoud2018sharp} and \cite{chen2021cutoff} revealed that exact recovery is achievable beyond the separation strength threshold $(2d^{\ast}n^{-1}\log n)^{1/4}$, whereas the exponential-type clustering error rate \citep{gao2022iterative,zhang2022leave} was derived only beyond the threshold $(d^{\ast}/n)^{1/2}$. To our best knowledge, the gap still exists at the moment. \cite{jin2017phase} proposed a two-component symmetric {\it sparse} GMM and investigated the phase transition in consistent clustering. Specifically, they showed that, ignoring log factors, $\Delta\gg 1+s/n$ is necessary for consistent clustering without restricting the computational complexity. Here $s$ is the sparsity of the expected observation. A recent work \citep{loffler2020computationally} designed an SDP-based spectral method and established an exponential-type clustering error rate  when $\Delta$ is greater than $1+s^{1/2}\log^{1/4}(d^{\ast})n^{-1/4}$.  Moreover, they provided evidence supporting the claim that no polynomial-time algorithm can consistent recover the clusters if $\Delta$ is smaller than the aforesaid threshold, i.e., there exists a statistical-to-computational gap for clustering in sparse GMM. Both \cite{jin2017phase} and \cite{loffler2020computationally} implied that the necessary separation strength primarily depends on the intrinsic dimension $s$ rather than the ambient dimension $d^{\ast}$. We remark that there is a vast literature studying the clustering problem for GMM. A representative but incomplete list includes \cite{lu2016statistical,balakrishnan2017statistical,dasgupta2008hardness,fei2018hidden,hajek2016achieving,verzelen2017detection,witten2010framework,abbe2020ellp} and references therein. 

In contrast, the understanding of the limit of clustering for LrMM is still at its infant stage. In this paper, we fill the void in the optimal clustering error rate of LrMM and demonstrate that the rate is achievable by a computationally fast algorithm. Challenges are posed from multiple fronts. First of all,  designing a computationally fast clustering procedure that sufficiently exploits low-rank structure is non-trivial. Unlike (sparse) GMM \citep{chen2021cutoff,loffler2020computationally}, convex relaxation seems not immediately accessible for the clustering of LrMM,  especially when there are more than two clusters. Non-convex approaches based on tensor decomposition and spectral clustering \citep{jing2021community,luo2022tensor,xia2019sup} usually cannot distinguish the sample size dimension (i.e., $n$) and data point dimension (i.e., $d_1, d_2$). Their theoretical results become sub-optimal when the sample size is much larger than $d_1$.  On the technical front, low-rankness makes deriving an exact exponential-type clustering error rate even more difficult. Under GMM \citep{gao2022iterative,loffler2021optimality}, the exponential-type clustering error rate is established by carefully studying the concentration phenomenon of a Gaussian linear form that usually admits an explicit representation. Estimating procedures under LrMM, however, often require multiple iterations of low-rank approximation, say, by singular value decomposition (SVD). Consequently, deriving the concentration property of respective linear forms under LrMM is much more involved than that under GMM. Moreover, prior related works \citep{loffler2020computationally,jin2017phase,zhang2018tensor,lyu2022optimal} provided evidences that imply the existence of a statistical-to-computational gap. It is unclear which model parameter characterizes such a gap and how the gap depends on the sample size and dimensions. For instance, how the low-rankness benefits the separation strength requirement? Interestingly, we discover that the gap is not determined by the separation strength $\Delta$ but rather by the signal strength (to be defined in Section~\ref{sec:method}) of the population center matrices.

Our main contributions are summarized as follows.  First, we propose a computationally fast clustering algorithm for LrMM. At its essence is the combination of Lloyd's algorithm \citep{lloyd1982least,lu2016statistical} and low-rank approximation. Basically, given the updated cluster memberships of each observation, the cluster centers are obtained by the SVD of the sample average within each cluster. The whole algorithm involves only K-means clustering and matrix SVDs. Secondly, we prove that, equipped with a good initial clustering, the low-rank Lloyd's algorithm converges fast and achieves the minimax optimal clustering error rate $\exp(-\Delta^2/8)$ with high probability as long as the separation strength satisfies $\Delta^2\gg 1+d_1r_{\submax}/n$ and the signal strength is strong enough. Here $r_{\submax}$ is the maximum rank among all the population center matrices. This dictates that a weaker separation strength is sufficient for clustering under LrMM if the rank $r_{\submax}=O(1)$. Our key technical tool to develop the exponential-type error rate is a spectral representation formula from \cite{xia2021normal}, which has helped push forward the understanding of statistical inference for low-rank models \citep{xia2021statistical,xia2022inference}. Thirdly, we propose a novel tensor-based spectral method for obtaining an initial clustering. Under similar separation strength and signal strength conditions, this method delivers an initial clustering that is sufficiently good for ensuring the convergence of low-rank Lloyd's algorithm. Lastly, compared with GMM that only requires a separation strength condition \citep{loffler2020computationally,gao2022iterative}, an additional signal strength condition seems necessary under LrMM. We provide evidences, based on the low-degree framework \citep{kunisky2019notes}, showing that if the signal strength condition fails, all polynomial-time algorithms cannot consistently recover the true clusters, even when the separation strength is much stronger than the aforesaid one. It is worth pointing out that, unlike tensor-based approaches \citep{jing2021community,luo2022tensor,xia2019sup},  our theoretical results impose no constraints on the relation between $n$ and $(d_1,d_2)$.

The rest of the paper is organized as follows. Low-rank mixture model is formalized in Section~\ref{sec:method}, and we introduce the low-rank Lloyd's algorithm and a tensor-based method for spectral initialization. The convergence performance of Lloyds' algorithm, minimax optimal exponential-type clustering error rate, and guarantees of a tensor-based spectral initialization are established in Section~\ref{sec:main}. We discuss the computational barriers of LrMM in Section~\ref{sec:comp}. In Section~\ref{sec:Weak-SNR}, we slightly modify the low-rank Lloyd's algorithm and derive the same minimax optimal clustering error rate requiring a slightly weaker signal strength condition. We discuss the difference between estimation and clustering under LrMM in Section~\ref{sec:comparison}. Further discussions are provided in Section~\ref{sec:discussion}. Numerical simulations and real data examples are presented in Section~\ref{sec:numerical}. All proofs and technical lemmas are relegated to the appendix.

\section{Methodology}\label{sec:method}

\subsection{Background and notations} 
For nonnegative $D_1, D_2$ , the notation $D_1\lesssim D_2$ (equivalently, $D_2\gtrsim D_1$) means that there exists an absolute constant $C>0$ such that $D_1\leq C D_2$; $D_1 \asymp D_2$ is equivalent to $D_1\lesssim D_2$ and $D_2\lesssim D_1$, simultaneously. 
Let $\|\cdot\|$ denote the $\ell_2$ norm for vectors and operator norm for matrices, and $\|\cdot\|_{\rm F}$ denotes the matrix Frobenius norm. Denote $\sigma_1(\bM)\geq \cdots\geq \sigma_r(\bM)>0$ the non-increasing singular values of $\bM$ where $r={\rm rank}(\bM)$. We also define $\sigma_{\submin}(\bM):=\sigma_r(\bM)$. A third order tensor is a three-dimensional array. Throughout the paper, a tensor is written in the calligraphic bold font, e.g. $\bcalM\in\RR^{d_1\times d_2\times n}$. We use $\scrM_1(\bcalM)$  to denote the mode-$1$ matricization of $\bcalM$ such that $\scrM_1(\bcalM)\in\RR^{d_1\times (d_2n)}$ and $\scrM_1(\bcalM)(i_1,(i_2-1)n+i_3)=\bcalM(i_1,i_2,i_3), \forall i_1\in[d_1], i_2\in[d_2], i_3\in[n]$. The mode-$2$ and mode-$3$ matricizations are defined in a similar fashion. Then $\big\{{\rm rank}\big(\scrM_k(\bcalM)\big): k=1,2,3\big\}$ are called {\it Tucker rank} or {\it multilinear rank}. The mode-$1$ marginal multiplication between $\bcalM$ and a matrix $\bU^{\top}\in\RR^{r\times d_1}$ results into a tensor of size $r_1\times d_2\times n$, whose elements are 
$$
\big(\bcalM\times_1 \bU^{\top}\big)(j_1,i_2,i_3):=\sum_{i_1=1}^{d_1} \bcalM(i_1,i_2,i_3)\bU(i_1,j_1), \quad \forall j_1\in[r], i_2\in[d_2], i_3\in[n]
$$
Similarly, we can define the mode-$2$ and mode-$3$ marginal multiplication. Given $\bcalS\in\RR^{r_1\times r_2\times r_3}, \bV\in\RR^{d_2\times r_2}, \bW\in\RR^{n\times r_3}$, the multi-linear product $\bcalM:=\bcalS\times_1 \bU\times_2 \bV\times_3 \bW$ outputs a $d_1\times d_2\times n$ tensor defined by,
\begin{equation}\label{eq:Tucker}
\bcalM(i_1,i_2,i_3):=\sum_{j_1=1}^{r_1}\sum_{j_2=1}^{r_2}\sum_{j_3=1}^{r_3} \bcalS(j_1,j_2,j_3)\bU(i_1,j_1)\bV(i_2,j_2)\bW(i_3,j_3)
\end{equation}
More details can be found in \cite{kolda2009tensor}. 
 Denote $\OO_{d,r}$ the set of all $d\times r$ matrices $\bU$ such that $\bU^{\top}\bU=\bI_r$, where $\bI_r$ is the $r\times r$ identity matrix. Eq. (\ref{eq:Tucker}) is known as the {\it Tucker decomposition} if $r_k={\rm rank}\big(\scrM_k(\bcalM)\big)$, $\bU\in\OO_{d_1,r_1}, \bV\in\OO_{d_2,r_2}$, and $\bW\in\OO_{n,r_3}$.

\subsection{Low-rank sub-Gaussian mixture}
Suppose that the $d_1\times d_2$ matrix-valued observations $\bX_1, \cdots, \bX_n$ are i.i.d., and each of them has a latent label $s_i^{\ast}\in[K]$. Here $K$ denotes the number of underlying clusters, and without loss of generality, assume $d_1\geq d_2$. We assume that there exists $K$ deterministic but {\it unknown} matrices $\bM_1,\cdots,\bM_K$ such that, conditioned on $s_i^{\ast}=k$,  $\bX_i$ has i.i.d.  zero-mean sub-Gaussian entries with the mean matrix $\bM_k$.
This implies that $\bX_i|s_i^{\ast}=k$ is equal to $\bM_k+\bE_i$ {\it in distribution} where the noise matrix $\bE_i$ satisfies the following assumption:
\begin{assumption}(Sub-Gaussian noise)\label{assump:sub-Gaussian}
The noise matrix $\bE_i$ has i.i.d. zero-mean entries and unit variance, and for $\forall \bM\in\RR^{d_1\times d_2}$, the following probability holds
$$
\PP(\langle \bM, \bE_i\rangle\geq t)\leq e^{-t^2/(2\sigmasg^2\cdot \|\bM\|_{\rm F}^2)},\quad \forall t> 0,
$$
where $\sigmasg>0$ is the sub-Gaussian constant. 
\end{assumption}
Throughout the paper,  we let $\sigmasg^2=1$ without loss generality (say, by substituting $\bX_i$ with $\bX_i/\sigmasg$).  Moreover, we assume that the latent labels $s_1^{\ast},\cdots, s_n^{\ast}$ are i.i.d. and
\begin{equation}\label{eq:mixture}
\PP(s_i^{\ast}=k)=\pi_k,\quad \forall k\in[K];\quad {\rm where}\quad \sum_{k=1}^K \pi_k=1.
\end{equation}
Here the unknown $\pi_k>0$ stands for the mass of $k$-th cluster. Put it differently, the matrix-valued observations have a marginal distribution
\begin{equation}\label{eq:LrMM}
\bX_1,\cdots,\bX_n\ \stackrel{{\rm i.i.d.}}{\sim}\ \sum_{k=1}^k \pi_k\cdot p_{\bM_k,\sigmasg^2}(\bX )
\end{equation}
where $p_{\bM_k,\sigmasg^2}(\bX )$ is the density function of matrix observation $\bX\in\RR^{d_1\times d_2}$ with independent entries of unit variance,  the sub-Gaussian constant $\sigmasg$ and mean matrix $\bM$.  Let $r_k={\rm rank}(\bM_k)$ and assume $r_k\ll d_2$ for all $k$, i.e., all the population center matrices are low-rank. Model (\ref{eq:LrMM}) is referred to as the {\it low-rank mixture model} (LrMM). For simplicity, we treat the ranks $r_k$'s as known and will briefly discuss how to estimate them in Section~\ref{sec:discussion}. We denote the compact SVD of population center matrices by $\bM_k=\bU_k\bSigma_k\bV_k^{\top}$ with $\bU_k\in \OO_{d_1,r_k}$ and $\bV_k\in\OO_{d_2,r_k}$. The {\it signal strength} of $\bM_k$ is characterized by $\sigma_{\submin}(\bM_k):=\sigma_{r_k}(\bM_k)$. We remark that estimating $K$ is a challenging question even under GMM. Hence, throughout this paper, it is assumed that $K$ is provided beforehand. 

\cite{sun2019dynamic} introduced a tensor Gaussian mixture model without specifically imposing low-rank structures on the center matrices. A similar tensor normal mixture model without low-rank assumptions is proposed by \cite{mai2021doubly}. Our LrMM can be viewed as a generalization of mixture multi-layer SBM proposed by \cite{jing2021community} and as an extension of the symmetric two-component case introduced by \cite{lyu2021latent}. Mixture of low-rank matrix normal models have also appeared in \cite{gao2021regularized} for image analysis.

Since our goal of current paper is to investigate the fundamental limits of clustering matrix-valued observations, hereafter, we view the latent labels $s_i^{\ast}, i\in[n]$ as a {\it fixed} realization sampled from the mixture distribution (\ref{eq:mixture}).  Then the matrix-valued observations can be written in the following form:
\begin{equation}\label{eq:LrMM-1}
\bX_i=\bM_{s_i^\ast}+\bE_i,\quad i\in[n]	
\end{equation}
Denote $\bs^\ast=(s^\ast_1,\cdots,s^\ast_n)$ the collection of true latent labels, known as the {\it cluster membership vector}. The size of each cluster is given by
$
n_k^{\ast}:=\sum_{i=1}^n \ind{s_i^{\ast}=k}, \forall k\in[K].
$
 With mild conditions under LrMM, Chernoff bound \citep{chernoff1952measure} guarantees $n_k^{\ast}\asymp n\pi_k$ with high probability.

Given an estimated cluster membership vector $\hat\bs:=(\hat s_1,\cdots,\hat s_n)\in[K]^n$, its clustering error is measured by the {\it Hamming distance} defined by 
\begin{equation}\label{eq:hamming}
h_{\textsf{c}}(\hat\bs,\bs^{\ast})=\min_{\pi: \text{ permutation of }[K]}\ \sum_{i=1}^n\ind{\hat s_i\ne \pi(s_i^{\ast})}
\end{equation}

For technical convenience,  we also define the the following Frobenious error related to $\bcalM$:
$$
\ell_{\textsf{c}}(\hat\bs,\bs^{\ast})=\min_{\pi: \text{ permutation of }[K]}\ \sum_{i=1}^n\big\| \bM_{\hat s_i}-\bM_{\pi (s_i^\ast)}\big\|_{\rm F}^2.
$$

\subsection{Low-rank Lloyd's algorithm}
Lloyd's algorithm \citep{lloyd1982least} or K-means algorithm is perhaps, conceptually and implementation-wise, the most simple yet effective method for clustering. It is an iterative algorithm, which consists of two main routines at each iteration: 1). provided with an estimated cluster membership vector, the cluster centers are updated by taking the sample average within every estimated cluster; 2). provided with the updated cluster centers, every data point is assigned an updated cluster label according to its distances from the cluster centers. The iterations are terminated once converged. The success of Lloyd's algorithm is highly reliant on a good initial clustering or initial cluster centers. It is proved by \cite{lu2016statistical} and \cite{gao2022iterative} that, if well initialized, Lloyd's algorithm converges fast and achieves minimax optimal clustering error for GMM and community detections under stochastic block model. 

The original Lloyd's algorithm updates the cluster centers by taking the vanilla sample average. This approach is sub-optimal under LrMM because the underlying low-rank structure is overlooked. It is well-known that exploiting the low-rankness can further de-noise the estimates. Towards that end, we propose the low-rank Lloyd's algorithm whose details are enumerated in Algorithm~\ref{alg:Lloyd}. Compared with the original Lloyd's algorithm, the low-rank version only modifies the procedure of updating the cluster centers. At the $(t+1)$-th iteration, given the current cluster labels $\hat\bs^{(t)}$ and for each $k$, we calculate the sample average $\bar \bX_k(\hat\bs^{(t)})$ defined as in Algorithm~\ref{alg:Lloyd}, and then update the cluster center by 
$$
\hat\bM_k^{(t+1)}:=\hat \bU_k^{(t)}\hat\bU_k^{(t)\top}\bar \bX_k(\hat\bs^{(t)}) \hat\bV_k^{(t)}\hat \bV_k^{(t)\top}
$$ 
where $\hat \bU_k^{(t)}$ and $\hat\bV_k^{(t)}$ are the top-$r_k$ left and right singular vectors of $\bar{\bX}_k(\hat\bs^{(t)})$, respectively. The update of cluster labels is unchanged compared with the original Lloyd's algorithm.

\begin{algorithm}[!htbp]
\caption{Low-rank Lloyd's Algorithm (lr-Lloyd)}\label{alg:Lloyd}
\begin{algorithmic}
	\STATE{\textbf{Input}: Observations $\bX_1,\cdots,\bX_n\in\RR^{d_1\times d_2}$, initial estimate $\hat\bs^{(0)}$, ranks $\{r_k\}_{k=1}^K$.}
	\FOR{$t = 1,\ldots,T$}
	\STATE{for each $k=1,\cdots,K$: \hspace{2cm} ({\it update cluster centers})
	\begin{equation}
		\hat\bM_k^{(t)}\leftarrow \textrm{best rank-}r_k \textrm{ approximation of } \bar {\bX}_k(\hat\bs^{(t-1)}):=\frac{\sum_{i=1}^n \ind{\hat s_{i}\suptt=k}\bX_i}{\sum_{i=1}^n \ind{\hat s_{i}\suptt=k}}
	\end{equation}\label{algo:est_M}
	}
	\STATE{for each $i=1,\cdots,n$: \hspace{2cm} ({\it update cluster labels})
	$$
	\hat s_i^{(t)}\leftarrow \underset{k\in [K]}{\arg\min}\ \|\bX_i-\hat\bM_k\supt\|_{\rm F}^2
	$$}
	\ENDFOR
	\STATE{\textbf{Output}: $\hat\bs:=\hat\bs^{(T)}$ }
\end{algorithmic}
\end{algorithm}

Conceptually, our low-rank Lloyd's algorithm is a direct adaptation of Lloyd's algorithm to accommodate low-rankness. However, the low-rank update of cluster centers poses fresh and highly non-trivial challenges in studying the convergence behavior of Algorithm~\ref{alg:Lloyd}. The original Lloyd's algorithm simply takes the sample average and thus admits a clean and explicit representation form for the updated centers, which plays a critical role in technical analysis, as in \cite{gao2022iterative}. In sharp contrast, the required SVD in Algorithm~\ref{alg:Lloyd} involves intricate and non-linear operations on the matrix-valued observations, and there is surely no clean and explicit representation form for $\hat\bM_k^{(t)}$. More advanced tools are in need for our purpose, as shall be explained in Section~\ref{sec:main}.

\subsection{Tensor-based spectral initialization}
The success of Algorithm~\ref{alg:Lloyd} crucially depends on a reliable initial clustering. A naive approach is to vectorize the matrix observations, concatenate them into a new matrix of size $n\times (d_1d_2)$, then borrow the classic spectral clustering method as in \cite{loffler2021optimality} and \cite{zhang2022leave}. Unfortunately, the naive approach turns out to be sub-optimal for ignoring the planted low-dimensional structure in the row space. 

Our proposed initial clustering is based on tensor decomposition. Towards that end, we construct a third-order data tensor $\bcalX\in\RR^{d_1\times d_2\times n}$ by stacking the matrix-valued observations slice by slice, i.e., its $i$-th slice\footnote{We follow Matlab syntax tradition and denote $\calX(:,:,i)$ the sub-tensor by fixing one index.} $\bcalX(:,:,i)=\bX_i$. The noise tensor $\bcalE$ is defined in the same fashion. The {\it signal tensor} $\bcalM$ is constructed such that $\bcalM(:,:,i)=\bM_{s_i^{\ast}}$. The tensor form of LrMM (\ref{eq:LrMM-1}) is 
\begin{equation}\label{eq:LrMM-tensor}
\bcalX=\bcalM+\bcalE
\end{equation}
Interestingly, eq. (\ref{eq:LrMM-tensor}) coincides with the famous tensor SVD or PCA model \citep{zhang2018tensor,xia2019sup,liu2022characterizing}. Let $\mathring{r}:=\sum_{k=1}^K r_k$. Indeed, the signal tensor $\bcalM$ admits the following low-rank decomposition
\begin{equation}\label{eq:M-dec1}
\bcalM=\bcalS\times_1 \bU\times_2  \bV\times_3 \bW
\end{equation}
where the $\mathring{r}\times \mathring{r}\times K$ core tensor $\bcalS$ is constructed as 
$$\bcalS(:,:,k):={\rm diag}({\bf 0}_{r_1},\cdots,{\bf 0}_{r_{k-1}}, \bSigma_{r_k}, {\bf 0}_{r_{k+1}},\cdots, {\bf 0}_{r_K})$$ 
and $\bU=(\bU_1,\cdots,\bU_K)\in\RR^{d_1\times \mathring{r}}$, $\bV=(\bV_1,\cdots,\bV_K)\in\RR^{d_2\times \mathring{r}}$, $\bW=(\be_{s_1^{\ast}},\cdots,\be_{s_n^{\ast}})^{\top}\in\{0,1\}^{n\times K}$. Here $\be_k$ denotes the $k$-th canonical basis vector in Euclidean space whose dimension might vary at different appearances. Clearly, the rows of $\bW$ provide the cluster information and is referred to as the cluster membership matrix. Note that (\ref{eq:M-dec1}) is not necessarily the Tucker decomposition since $\bU, \bV$ might be rank-deficient, in which case the decomposition in the form (\ref{eq:M-dec1}) is not unique and $\bU, \bV$ become unrecoverable.

The singular space of $\bcalM$ is uniquely characterized by its Tucker decomposition. To this end, denote $\bU^{\ast}\in \OO_{d_1,r_{\bU}}$ and $\bV^{\ast}\in\OO_{d_2,r_{\bV}}$ the left singular vectors of $\bU$ and $\bV$, respectively. Here, $r_{\bU}$ and $r_{\bV}$ are the ranks of $\scrM_1(\bcalM)$ and $\scrM_2(\bcalM)$, respectively. Define $\bW^{\ast}\in\OO_{n,K}$ by normalizing the columns of $\bW$. Re-compute the core tensor $\bcalS^{\ast}:=\bcalM\times _1 \bU^{\ast\top}\times_2 \bV^{\ast\top}\times_3 \bW^{\ast\top}$ that is of size $r_{\bU}\times r_{\bV}\times K$. Finally, we re-parameterize the signal tensor via its Tucker decomposition
\begin{equation}\label{eq:M-dec2}
\bcalM=\bcalS^{\ast}\times_1 \bU^{\ast}\times_2 \bV^{\ast}\times_3 \bW^{\ast}
\end{equation}
Here $\bU^{\ast}, \bV^{\ast}, \bW^{\ast}$ are usually called the singular vectors of $\bcalM$. Still, the rows of $\bW^{\ast}$ tell the cluster information in  that $\bW^{\ast}(i,:)=\bW^{\ast}(j,:)$ iff $s_i^{\ast}=s_j^{\ast}$, i.e, $i, j$ belongs to the same cluster. We note that there are interesting special cases concerning the values of $r_{\bU}, r_{\bV}$. For instance, if $r_{\bU}=r_{\bV}=r_1$, it implies that all the population center matrices share the same low-dimensional singular space with $\bM_1$, which simplifies theoretical investigate of our proposed initialization method.  Another special case is $r_{\bU}=r_{\bV}=\mathring{r}$, namely the singular spaces of all population center matrices are separated to a certain degree. Intuitively, the clustering problem becomes easier. See Section~\ref{sec:main-init} for discussions of both cases.

We now present our tensor-based spectral method for initial clustering. Unlike the aforementioned naive spectral method, ours is specifically designed to exploit the low-rank structure of $\bcalM$ in the $1$st and $2$nd dimension. Without loss of generality, we treat $r_{\bU}$ and $r_{\bV}$ as known here and shall discuss ways to estimate them in Section~\ref{sec:discussion}. Our method consists of three crucial steps with details in Algorithm~\ref{alg:initial}. Step 1 aims to estimate the singular vectors $\bU^{\ast}$ and $\bV^{\ast}$. Here, higher order SVD (HOSVD) is obtained by applying SVD to the matricizations $\scrM_1(\bcalM)$ and $\scrM_2(\bcalM)$. See, for instance, \cite{de2000multilinear} and \cite{xia2019sup}. The estimated singular vectors are used for denoising in Step 2 by projecting the noise into a low-dimensional space. Step 3 applies the classical K-means clustering \citep{loffler2021optimality,zhang2022leave} to the denoised observations. Note that solving K-means is generally NP-hard \citep{mahajan2009planar}, but there exist fast algorithms \citep{kumar2004simple} achieving an approximate solution.

\begin{algorithm}[!htbp]
\caption{Tensor-based Spectral Initialization (TS-Init)}\label{alg:initial}
\begin{algorithmic}
	\STATE{\textbf{Input}: Observations $\bX_1,\cdots,\bX_n\in\RR^{d_1\times d_2}$ or a tensor $\bcalX\in\RR^{d_1\times d_2\times n}$ by concatenating the matrix observations slice by slice.}
	\begin{enumerate}
		\item Obtain the estimated factor matrices $\hat \bU$ and $\hat \bV$ by applying HOSVD to the tensor $\bcalX$ in mode-$1$ and mode-$2$ with rank $r_\bU$ and $r_\bV$, respectively.
		\item Project $\bcalX$ onto the column space of $\hat\bU$ and $\hat\bV$ by $$\hat\bcalG:=\bcalX\times_1\hat \bU\hat \bU^\top \times_2\hat\bV\hat \bV^\top \in\RR^{d_1\times d_2\times n}$$
		\item Apply k-means on rows of $\hat\bG:=\scrM_3(\hat\bcalG)\in \RR^{n\times d_1d_2}$ to obtain initializer for $\bs^\ast$, i.e.
	$$
	(\hat\bs^{(0)},\{\hat\bM_k^{(0)}\}_{k=1}^K):=\argmin_{\bs\in[K]^n,\{\bM_k\}_{k=1}^K,\bM_k\in\RR^{d_1\times d_2},\forall k}\sum_{i=1}^n\op{[\hat \bG]_{i\cdot }-vec(\bM_{s_i})}^2
	$$
	\end{enumerate}
	\STATE{\textbf{Output}:  $\hat\bs^{(0)}$}
\end{algorithmic}
\end{algorithm}

Algorithm~\ref{alg:initial} improves the naive spectral clustering whenever $\hat\bU$ and $\hat\bV$ are reliable estimates of their population counterparts. 
This suggests that a certain signal strength condition on $\scrM_1(\bcalS^{\ast})$ and $\scrM_2(\bcalS^{\ast})$ is necessary. We remark that the higher order orthogonal iteration (HOOI, \cite{zhang2018tensor}) algorithm for tensor decomposition is not suitable for our purpose since it requires a lower bound on $\sigma_{\submin}\big(\scrM_3(\bcalS^{\ast})\big)$, which is too restrictive under LrMM. See Section~\ref{sec:main-init} for more explanations.

\section{Minimax Optimal Clustering Error Rate of LrMM}\label{sec:main}
In this section, we establish the convergence performance of low-rank Lloyd's algorithm, validate our tensor-based spectral initialization, and derive the minimax optimal clustering error rate for LrMM (\ref{eq:LrMM}). The hardness of clustering under LrMM is determined primary by two quantities: 
\begin{align*}
\textit{Separation strength }& \Delta:=\min_{a\ne b,a,b\in[K]}\fro{\bM_a-\bM_b}
\end{align*}
The separation strength  is a generalization of the minimum $\ell_2$ distance between different population centers under GMM \citep{lu2016statistical,chen2021cutoff,gao2022iterative}, which characterizes the intrinsic difficult in clustering the observations. In fact, the minimax optimal error rate, i.e, the best achievable clustering accuracy, is exclusively decided by $\Delta$.

\subsection{Iterative convergence of low-rank Lloyd's algorithm}
The performance of Lloyd's algorithm also relies on the minimal cluster size \citep{lu2016statistical}. To this end, define $\alpha:= \min_{k\in[K]}n_k^\ast\cdot (n/K)^{-1}$, where recall that $n_k^\ast:=|\{i\in[n]:s_i^\ast=k\}|$ is the size of $k$-th cluster.  The cluster sizes are said to be {\it balanced} if $\alpha\asymp 1$. The hamming distance $h_{\textsf{c}}(\hat\bs, \bs^{\ast})$ is defined as in eq. (\ref{eq:hamming}). Without loss of generality, we assume $r:=r_1$ is the largest amongst $\{r_k: k\in[K]\}$ and $d:=d_1\geq d_2$. 

Due to technical reasons, we define $\kappa_0:={\max_{k\in[K]}\|\bM_k\|}/\min_{k\in[K]} \sigma_{\submin}(\bM_k)$, which can be viewed as the maximum condition number of all population center matrices. It usually does not appear in the literature of GMM, but is of unique importance under LrMM. This quantity plays a critical role in connecting the accuracy of updated center matrix $\hat\bM_k\supt$ to the current clustering accuracy $h_{\textsf{c}}(\hat\bs\suptt, \bs^{\ast})$.  Since $\hat\bM_k\supt$ stems from the SVD of $\bar\bX_k(\hat\bs\suptt)$, whose accuracy is characterized by the strength of signal $\bM_k$ and size of perturbation $\bar\bX_k(\hat\bs\suptt)-\bM_k$. Besides random noise, the latter term, roughly, consists of $(n_a^{\ast})^{-1}h_{\textsf{c}}(\hat\bs^{(t-1)},\bs^{\ast})\big(\bM_{k'\neq k}-\bM_{k}\big)$, whose operator norm can be controlled by $O\big((n_a^{\ast})^{-1}h_{\textsf{c}}(\hat\bs\suptt, \bs^{\ast}) \kappa_0\sigma_{\submin}(\bM_k)\big)$.  Hence $\kappa_0$ is, {\it perhaps}, the unavoidable price to be paid for taking advantage of low-rankness \citep{chen2021learning}.

The following theorem presents the convergence performance of low-rank Lloyd's algorithm (Algorithm \ref{alg:Lloyd}). Due to the local nature of Lloyd's algorithm, its success highly relies on a good initialization. Theorem~\ref{thm:main} assumes the initial clustering is consistent, i.e., initial clustering error approaches zero asymptotically as $n\to\infty$. Under suitable conditions of separation strength and signal strength, the output of Algorithm~\ref{alg:Lloyd} attains an exponential-type error rate. The constant factor $1/8$ in the exponential rate exactly matches the minimax lower bound in Theorem \ref{thm:minimax-lower-bound}. Notice that our result is non-asymptotic, and all asymptotic conditions in Theorem \ref{thm:main} are to guarantee the sharp constant $1/8$ in eq. \eqref{eq:contrac-exp-bound}.  More precisely, through a careful inspection on our analysis, the implicit term $o(1)$ in the exponential rate in Theorem \ref{thm:main} can be chosen at the order $\Omega\left(\left(Kr(d+\log n)(\alpha n)^{-1}/\Delta^2\right )^{1/2-\epsilon}\right)=o(1)$ for any fixed $\epsilon\in (0,1/2)$.


\begin{theorem}\label{thm:main}
Suppose $d\ge C_0\log K$ for some absolute constant $C_0>0$. Assume that 
\begin{enumerate}
	\item[(i)] initial clustering error:
	\begin{equation}\label{init-cond}
		n^{-1}\cdot \ell_{\textsf{c}}(\hat\bs^{(0)},\bs^\ast)= o\left(\frac{\alpha }{\kappa_0^2K}\Delta^2\right)
	\end{equation}
	\item[(ii)] separation strength:
	\begin{equation}\label{cond:min-sep}
		\frac{\Delta^2}{\alpha^{-1}(\kappa_0^2\vee Kr)Kr\left(\frac{d}{n}+1\right)} \rightarrow \infty
	\end{equation}
\end{enumerate}
Let $\hat\bs^{(t)}$ be the cluster labels at $t$-th iteration generated by Algorithm~\ref{alg:Lloyd}. Then, for all $t\geq 1$, we have 
\begin{equation}\label{eq:contrac-exp-bound}
	n^{-1}\cdot h_{\textsf{c}}(\hat\bs^{(t)},\bs)\le \exp\left(-(1-o(1))\frac{\Delta^2}{8}\right)+\frac{1}{2^t}
\end{equation}
with probability at least $1-\exp(-\Delta)-\exp(-c_0d)$ with some absolute constant $c_0>0$. 
\end{theorem}

By Theorem~\ref{thm:main}, after at most $O\big(\min\{\Delta^2, \log n\}\big)$ iterations,
our low-rank Lloyd's algorithm achieves the minimax optimal clustering error rate $\exp(-\Delta^2/8)$, which is the same optimal rate for classical GMM \citep{lu2016statistical,loffler2020computationally,gao2022iterative,zhang2022leave} and is exclusively decided by the separation strength $\Delta$. It is worth noting  that provided with good initialization, lr-Lloyd solely requires separation strength strong enough to achieve such optimal rate.  

{\it Blessing of low-rankness and comparison with GMM.} If low-rankness is ignored so that LrMM is treated as GMM, the exponential-type error rate is established only in the regime of separation strength $\Delta\gg 1+(d_1d_2/n)^{1/2}$ \citep{gao2022iterative,zhang2022leave}. In contrast, our condition (\ref{cond:min-sep}) only requires $\Delta\gg 1+(d_1/n)^{1/2}$ if $r,K,\kappa_0,\alpha=O(1)$.  

{\it Discussions on separation strength  $\Delta$.} The separation strength condition is typical in the literature of clustering problems \citep{vempala2004spectral,loffler2021optimality}. To see why our condition (\ref{cond:min-sep}) is minimal, without loss of generality, consider the case $\alpha\asymp 1$ and $K=2$. Moreover, assume the singular vectors $\bU_1=\bU_2$ and $\bV_1=\bV_2$, and they are already known. One can multiply each observation by $\bU_1^{\top}$ from left and by $\bV_1$ from right, which reduces LrMM to GMM in the dimension $r^2$. Literature of GMM \citep{gao2022iterative,loffler2021optimality,zhang2022leave} all impose a separation strength condition $\Delta\gg 1$. This certifies the constant $1$ in eq. (\ref{cond:min-sep}). To understand the term $(rd/n)^{1/2}$, consider that the true labels of first $n-1$ observations are revealed to us and our goal is to estimate the label of the $n$-th sample $\bX_{n}$. A natural way is to first estimate the population centers utilizing the given labels $\bs^\ast_1,\cdots,\bs_{n-1}^\ast$, denoted by $\hat\bM_1$ and $\hat\bM_2$, respectively.  The literature of matrix denoising \citep{cai2018rate,xia2021normal,gavish2017optimal} tells that the minimax optimal estimation error is at the order $\|\hat\bM_1-\bM_1\|_{\rm F}\asymp \|\hat\bM_2-\bM_2\|\asymp (rd/n)^{1/2}$. Thus $\Delta\gg (rd/n)^{1/2}$ is necessary for consistently distinguishing the two clusters. The above rationale suggests that our separation strength condition (\ref{cond:min-sep}) might be minimal up to the order of $n$, if only the exponential-type error rate is sought.

We explained a gap concerning the separation strength in existing literature of GMM. Under GMM with dimension $d^{\ast}=d_1d_2$ and $n\leq d^{\ast}$, the exponential-type rate \citep{gao2022iterative,zhang2022leave} is established in the regime $\Delta\gg (d^{\ast}/n)^{1/2}$, whereas exact clustering results \citep{ndaoud2018sharp,chen2021cutoff} are attained in the regime $\Delta\gtrsim (d^{\ast}n^{-1}\log n)^{1/4}$. This leaves a natural question under LrMM: is the separation strength condition (\ref{cond:min-sep}) is relaxable to the scale $n^{-1/4}$? Unfortunately, answering this question is perhaps more challenging than that under GMM. We note that \cite{ndaoud2018sharp} and \cite{chen2021cutoff} achieve the $O(n^{-1/4})$ barrier by focusing entirely on clustering and by circumventing the estimation of population centers. Nonetheless, under LrMM, exploiting the low-rank structure demands estimating the population center matrices. We suspect, together with the aforementioned special examples, that condition (\ref{cond:min-sep}) might not be improvable in terms of the order of $n$. Anyhow,  It's unclear whether one can obtain a sharper characterization of $\Delta$ under LrMM using other methods like SDP. Further investigation in this respect is out of the scope of current paper.

\subsection{Guaranteed initialization}\label{sec:main-init}
Besides the separation strength condition, Theorem~\ref{thm:main} requires a consistent initial clustering. We now demonstrate the validity of tensor-based Algorithm~\ref{alg:initial}. Observe that denoising by spectral projection (Step 2 of Algorithm~\ref{alg:initial}) is only beneficial if $\hat\bU$ and $\hat\bV$ are properly aligned with $\bU^{\ast}$ and $\bV^{\ast}$, respectively. For that purpose, the signal strengths of $\scrM_1(\bcalM)$ and $\scrM_2(\bcalM)$, i.e., $\sigma_{\submin}\big(\scrM_1(\bcalM)\big)$ and $\sigma_{\submin}\big(\scrM_2(\bcalM)\big)$, needs to be sufficiently strong.  For simplicity, we let  $\Lambda_{\submin}:=\min_{j=1,2} \left\{\sigma_{\min}(\scrM_j(\bcalM))\right\}$ denote \textit{tensor signal strength in 1st and 2nd modes}  of $\bcalM$, or simply the \textit{tensor signal strength} of $\bcalM$. Note that this is a slightly different definition from classical tensor literature, where the signal strength is usually defined as $\min_{j=1,2,3} \left\{\sigma_{\min}(\scrM_j(\bcalM))\right\}$. See remark after Theorem \ref{thm:spec-initialization}.


\begin{theorem}\label{thm:spec-initialization} Let $\hat\bs^{(0)}$ be the initial clustering output by Algorithm~\ref{alg:initial}. There exists some absolute constant $c,C_1,C_2,C_3,C_4>0$ such that if  
\begin{align}\label{cond:Lambda}
	\Lambda_{\submin}\ge C_1(rK)^{1/2}d^{1/2}n^{1/4},
\end{align}
and
\begin{align}
	\Delta^2\ge  {C_2\alpha^{-1}K^2}\left(\frac{dKr}{n}+1\right),
\end{align}
we get, with probability at least $1-\exp(-c(n\wedge d))$, that
$$
 n^{-1}\cdot h_{\textsf{c}}(\hat\bs^{(0)},\bs^\ast)\le C_3\frac{{K}}{\Delta^2}\left(\frac{dKr}{n}+1\right),
$$
and 
\begin{align*}
	n^{-1}\cdot \ell_{\textsf{c}}(\hat\bs^{(0)},\bs^\ast)\le C_4{\gamma^2K}\left(\frac{dKr}{n}+1\right),
\end{align*}
where $\gamma:=\max_{a\ne b\in[K]}\fro{\bM_a-\bM_b}/\Delta$.
\end{theorem}

Theorem~\ref{thm:spec-initialization} suggests that Algorithm~\ref{alg:initial} delivers a consistent clustering if the separation strength $\Delta^2\gg K(1+rdK/n)$. In terms of loss function $\ell_c(\cdot )$, we have  an additional dependence on $\gamma $, which relates $\Delta$ to \textit{maximum separation strength}.  A similar condition is also casted in the vector GMM \citep{lu2016statistical}.  As argued in \cite{lu2016statistical,jin2016local},  a distant cluster can cause local search to fail which indicates the possibly  unavoidable dependence on $\gamma$.  Furthermore, Theorem~\ref{thm:spec-initialization} imposes a condition on the tensor signal strength $\Lambda_{\submin}$, which is not needed in Theorem \ref{thm:main}. Such an eigen-gap type condition is prevalent in low-rank models \citep{zhang2018tensor,richard2014statistical,levin2019recovering,xia2021normal,lyu2022optimal} as it determines whether the population centers or their singular spaces are estimable by polynomial-time algorithms,  only in which case the low-rank structure can be beneficial. Remarkably, $\Lambda_{\submin}$ also governs the computational and statistical limit under LrMM as will be explained  in Section \ref{sec:comp}. 


Finally, by combining Theorem \ref{thm:spec-initialization} and Theorem \ref{thm:main},  the successes of Algorithm \ref{alg:Lloyd} and Algorithm \ref{alg:initial} require the signal strength and separation strength conditions
$$
\Lambda_{\submin}\ge C_1(rK)^{1/2}d^{1/2}n^{1/4}
$$
and 
$$
\quad \frac{\Delta^2}{\alpha^{-1}\gamma^2(\kappa_0^2\vee Kr)Kr\left(\frac{dKr}{n}+1\right)}\rightarrow\infty  
$$
To facilitate a clearer understanding of $\Lambda_{\submin}$, we introduce the concept of \textit{individual signal strength} denoted by $\lambda$. This quantity, which is common in low-rank matrix literature, is defined as the minimum value of the smallest singular value among $\bM_k$'s, i.e., 
$$ \lambda:=\min_{k\in[K]} \sigma_{\submin}(\bM_k)$$

{\it Relation between tensor signal strength $\Lambda_{\submin}$ and individual matrix signal strength $\lambda$.} Define the condition number of $\bcalM$ in the mode-$j$ as $\kappa_j:=\op{\scrM_j(\bcalM)}/\sigma_{\min}(\scrM_j(\bcalM))$ for $j=1,2$. 
\begin{lemma}\label{lem:sv-lower-bound} 
For $j\in\{1,2\}$, $\sigma_{\min}(\scrM_j(\bcalM))\ge \kappa_j^{-1}(Kr)^{-1/2}\sqrt{n}\lambda$.
\end{lemma}
By Lemma \ref{lem:sv-lower-bound}, a sufficient condition for \eqref{cond:Lambda} to hold can be casted as $\lambda\ge C_0(\kappa_1\vee\kappa_2)rKd^{1/2}n^{-1/4}$. Recall that $\kappa_0$ tells whether {\it individual} population center matrices are well-conditioned. Here $\kappa_1$ ($\kappa_2$, resp.) measures the goodness of alignment among the column (row, resp.) spaces of {\it all} population center matrices. However, the exact relation between $\kappa_1$ and the column spaces $\{{\rm ColSpan}(\bU_k^{\ast})\}_{k=1}^K$ can be intricate.  The following lemma unfolds two special cases. Recall that $r_{\bU}$ and $r_{\bV}$ are the ranks of $\bU=(\bU_1,\cdots,\bU_K)$ and $\bV=(\bV_1,\cdots,\bV_K)$, respectively, and $\mathring{r}=\sum_{k=1}^K r_k$. Denote $\kappa(\bU)$ and $\kappa(\bV)$ the condition numbers of $\bU$ and $\bV$, respectively. The following indicates the connection between $\kappa_j$ and $\kappa_0$.
\begin{lemma}\label{lem:kappa12}
Let $\bcalM$ admits low-rank decomposition~(\ref{eq:M-dec1}). We have 
\begin{align*}
\scrM_1(\bcalM)\scrM_1^{\top}(\bcalM)=&\bU\cdot {\rm diag}\big(\{n_k^{\ast}\bSigma_k^2\}_{k=1}^K\big)\cdot \bU^{\top}\\
\scrM_2(\bcalM)\scrM_2^{\top}(\bcalM)=&\bV\cdot {\rm diag}\big(\{n_k^{\ast}\bSigma_k^2\}_{k=1}^K\big)\cdot \bV^{\top}
\end{align*}
and $\kappa_1\le\kappa_0\kappa(\bU)\cdot (n^{\ast}_{\submax}/n^{\ast}_{\submin})^{1/2}$ and $\kappa_2\le \kappa_0\kappa(\bV)\cdot (n^{\ast}_{\submax}/n^{\ast}_{\submin})^{1/2}$ where $n^{\ast}_{\submin}:=\min_k n_k^{\ast}$ and $n^{\ast}_{\submax}:=\max_k n_k^{\ast}$. 
If $r_{\bU}=r_{\bV}=r_1$, i.e., all the population center matrices share the same singular space with $\bM_1$, we have
$
\max\{\kappa_1,\kappa_2\}\leq \kappa_0 \cdot (K^2/\alpha )^{1/2};
$
if $r_{\bU}=r_{\bV}=\mathring{r}$ and $\bM_k$ has mutually orthogonal singular space, we have $\max\{\kappa_1,\kappa_2\}\leq \kappa_0 \cdot (K/\alpha )^{1/2}$. 
\end{lemma}
According to Lemma~\ref{lem:kappa12}, the unfolded matrices $\scrM_1(\bcalM)$ and $\scrM_2(\bcalM)$ are well-conditioned if $\bU$ and $\bV$ are well-conditioned. Interestingly, this implies that our tensor-based spectral initialization becomes more efficient when the population center matrices $\bM_k$'s have {\it either perfectly aligned singular spaces or nearly orthogonal singular spaces}.

{\it Discussions on tensor signal strength $\Lambda_{\submin}$.} 
Condition \eqref{cond:Lambda} reflects the computational difficulty under LrMM.  This intrinsic computational condition is likely attributed to the tensor method,  which is solely present in the initialization stage (Algorithm \ref{alg:initial}). Once well initialized,  the requirement for $\Lambda_{\min }$ vanishes in Theorem \ref{thm:main} for lr-Lloyd (Algorithm \ref{alg:Lloyd}). Such conditions are common in tensor problems \cite{zhang2018tensor,auddy2022estimating,richard2014statistical,luo2022tensor}.  A more relevant work \cite{lyu2022optimal} provides evidence showing that no polynomial time can consistently {\it estimate} the population centers even in the symmetric two-component LrMM if $\Lambda_{\submin}=o(d^{1/2}n^{1/4})$. In Section~\ref{sec:comp}, evidences are provided showing that the same phenomenon exists for clustering, that is, if $\Lambda_{\submin}=o(d^{1/2}n^{1/4})$, consistent clustering is impossible by any polynomial time algorithms even when the separation strength $\Delta$ is much stronger than the minimal condition (\ref{cond:min-sep}). 

{\it Comparison with HOOI \citep{zhang2018tensor} and the condition number of $\scrM_3(\bcalM)$. } Algorithm~\ref{alg:initial} looks similar to HOOI \citep{zhang2018tensor}, which uses HOSVD for {\it mode-wise} spectral initialization and applies power iterations to further improve the estimates of singular spaces. Indeed, \eqref{cond:Lambda} is analogous to the signal strength condition for HOOI therein to succeed.  However, the {\it mode-wise} HOSVD and subsequent power iterations both require a lower bound on $\sigma_{\submin}\big(\scrM_k(\bcalM)\big), k=1,2,3$. While our Theorem~\ref{thm:spec-initialization} also requires a lower bound  on $\sigma_{\submin}\big(\scrM_1(\bcalM)\big)$ and $\sigma_{\submin}\big(\scrM_2(\bcalM)\big)$, we emphasize that a similar lower bound on $\sigma_{\submin}\big(\scrM_3(\bcalM)\big)$ is too strong and trivialize the whole problem. To see this, just notice via definition that $\Delta\geq \sigma_{\submin}\big(\scrM_3(\bcalM)\big)/2$.

	{\it Comparison with \cite{han2022exact}.} A tensor block model was proposed by \cite{han2022exact},  which can be regarded as an extension of the stochastic block model.  They developed the high-order Lloyd's algorithm (HLloyd) with spectral initialization.  The two works differ drastically from  several aspects.  From the algorithmic perspective,  HLloyd doesn't require  low-rank approximation at all since it explores block structure rather than low-rank structure.  The membership matrix in \cite{han2022exact} (analogous to $\bU_k$ in this paper) lies in the space $\{0,1\}^{d_k\times r_k}$, which is  more informative owing to its discrete structure.  Clearly,  block model is just a special case of low-rank model and HLloyd is inapplicable to our LrMM.    On the technical front,  HLolyd updates the block means simply by the sample average which admits an explicit and clean representation form.  In sharp contrast,  the analysis for lr-Llyod is much more challenging due to the implicit and complicated form of the updated cluster centers $\hat{\bM}_k\supt$ defined in \eqref{algo:est_M}, which calls for more advanced tools.  
	

\subsection{Minimax lower bound}
Theorem~\ref{thm:main} has shown that the low-rank Lloyd's algorithm achieves the asymptotical clustering error rate $\exp(-\Delta^2/8)$. In this section, a matching minimax lower bound is derived showing that the aforesaid rate is indeed optimal in the minimax sense. A lower bound under GMM has been established by \cite{lu2016statistical}.  We follow the arguments in \cite{gao2018community} to establish the minimax lower bound for LrMM. Observe that the error rate only depends on the separation strength $\Delta$ implying that the dimension $d_1, d_2$ and ranks $r_k$'s  play a less important role here.

Define the following parameter space for the population center matrices and arrangements of latent labels:
\begin{align*}
\Omega_\Delta\equiv\Omega(\Delta, d_1,d_2,n,K,\alpha):=\Big\{(\{&\bM_k\}_{k=1}^K,\bs):~\bM_k\in\mathbb{R}^{d_1\times d_2},\text{rank}(\bM_k)=r_k, \bs\in[K]^n,\\
&\min_{k\in[K]} |\{i\in[n]:s_i=k\}|\ge \alpha n/K,\min_{a\ne b}\fro{\bM_a-\bM_b}\ge \Delta\Big\}	
\end{align*}
For notation simplicity, we omit its dependence on the ranks $r_k$'s. 

\begin{theorem}\label{thm:minimax-lower-bound} Let $\bX_1,\cdots,\bX_n$ satisfy LrMM (\ref{eq:LrMM}) with $(\{\bM_k\}_{k=1}^K, \bs^{\ast})\in \Omega_{\Delta}$. Suppose $\{\bE_i\}_{i=1}^n$ has i.i.d $\calN(0,\sigma^2)$ entries. If ${\Delta^2}/\left(\sigma^2\log(K/\alpha)\right)\rightarrow \infty$ as $n\to\infty$, we have
	$$\inf_{\hat \bs}\sup_{(\{\bM_k\}_{k=1}^n, \bs^{\ast})\in\Omega_\Delta}\ \E  \frac{h_{\textsf{c}}(\hat \bs, \bs^{\ast})}{n}\ge \exp\left(-(1+o(1))\frac{\Delta^2}{8\sigma^2}\right)$$
where $\underset{\hat \bs}\inf$ is taken over all clustering algorithms.
\end{theorem}

Compared to Theorem~\ref{thm:main} and Theorem~\ref{thm:spec-initialization}, the minimax lower bound is established only requiring a separation strength $\Delta^2\gg 1$ assuming $K/\alpha=O(1)$. Theorem~\ref{thm:minimax-lower-bound} holds for any signal strength and the infimum is taking over all possible clustering algorithms without considering their computational feasibility.  Here, an algorithm is said {\it computationally feasible} if it is computable within a polynomial time complexity in terms of $n$ and $d_1, d_2$.

\section{Computational Barriers}\label{sec:comp}
We now turn to the computational hardness of LrMM. For simplicity, we set $\alpha, K, r\asymp 1$ throughout this section. 
Our signal strength condition (\ref{cond:Lambda}) in initialization requires a lower bound $\Lambda_{\submin}\gtrsim d^{1/2}n^{1/4}$. The purpose of this section is to provide evidences on its necessity to guarantee computationally feasible clustering algorithms. Our evidence is built on the {\it low-degree likelihood ratio} framework for hypothesis testing proposed by \cite{kunisky2019notes, hopkins2018statistical}, which has delivered convincing evidences justifying the computational hardness under sparse GMM \citep{loffler2020computationally} and for sparse PCA \citep{ding2019subexponential}. 

Suppose that, given i.i.d. observations $\bX_1,\cdots,\bX_n$, one is interested in the computational and statistical limit in distinguishing two hypothesis $\QQ_n$ and $\PP_n$, i.e, 
\begin{equation}\label{eq:HT}
H_0^{(n)}: \bX_1\sim \QQ_n\quad {\rm versus}\quad H_1^{(n)}: \bX_1\sim \PP_n
\end{equation}
The above two hypotheses are said {\it statistically indistinguishable} if no test can have both type I and type II error probabilities vanishing asymptotically. 
The famous Neyman-Pearson lemma tells us that the likelihood ratio test based on $L_n(\bcalX):=d\PP_n/d\QQ_n(\bX_1,\cdots,\bX_n)$ has a preferable power and is uniformly most powerful under some scenarios. A well recognized fact is that $\QQ_n$ and $\PP_n$ are statistically indistinguishable if the quantity $\|L_n\|^2:=\EE_{\QQ_n}[L_n(\bcalX)^2]$ remains bounded as $n\to\infty$. See \cite{kunisky2019notes} for a simple proof. 

While the asymptotic magnitude of $\|L_n\|^2$ is informative for understanding the statistical limit of testing (\ref{eq:HT}), it does not directly reflect the computational limit of testing (\ref{eq:HT}). Towards that end, the low-degree likelihood ratio framework seeks a polynomial approximation of $L_n(\bcalX)$ and investigates the magnitude of the resultant approximation. More exactly, let $L_n^{\leq D}(\bcalX)$ be the orthogonal projection of $L_n(\bcalX)$ onto the linear space spanned by polynomials $\RR^{d_1\times d_2\times n}\mapsto \RR$ of degrees at most $D$. Similarly, define $\|L_n^{\leq D}\|^2:=\EE_{\QQ_n}[L_n^{\leq D}(\bcalX)^2]$. \cite{kunisky2019notes} conjectures that the asymptotic magnitude of $\|L_n^{\leq D}\|^2$ reflects the computational hardness of testing the hypothesis (\ref{eq:HT}). More formally, their conjecture, slightly adapted for our purpose, can be written as follows. It has been introduced in \cite{lyu2022optimal}. Here, a test $\phi_n(\cdot)$ taking value $1$ means rejecting the null hypothesis and takes value $0$ if the null hypothesis is not rejected. Thus $\EE_{\QQ_n}[\phi_n(\bcalX)]$ and $\EE_{\PP_n}[1-\phi_n(\bcalX)]$ stands for type-I and type-II error, respectively. 

\begin{conjecture}[\cite{lyu2022optimal}]\label{conjecture:comp-hardness}
	If there exists $\epsilon>0$ and $D=D_n\ge(\log nd)^{1+\epsilon}$ for which $\op{L_n^{\le D}}=1+o(1)$ as $n\to\infty$, then there is no polynomial-time test $\phi_n:\RR^{d_1\times d_2\times n}\mapsto\{0,1\}$ such that the sum of type-I error and type-II error probabilities 
	$$\EE_{\QQ_n}[\phi_n(\bcalX)]+\EE_{\PP_n}[1-\phi_n(\bcalX)]\rightarrow 0 \quad \text{as}\quad n\rightarrow \infty$$ 
\end{conjecture}

Based on this conjecture, \cite{kunisky2019notes} reproduces the sharp phase transitions for the spiked Wigner matrix model and the widely-believed statistical-to-computational gap in tensor PCA, and \cite{lyu2022optimal} develops a computational hardness theory for estimating the population low-rank matrices under LrMM. 

Note that a specific hypothesis $\PP_n$ is necessary to apply Conjecture~\ref{conjecture:comp-hardness} and investigate the computational barriers in clustering for LrMM. Towards that end, we consider a symmetric two-component LrMM as in \cite{lyu2022optimal}.  It is a special case of model (\ref{eq:LrMM}) with $K=2$, $r_1=r_2=1$, $\bM_1=n^{-1/2}\Lambda_{\submin} \bu\bv^{\top}$ and $\bM_2=-\bM_1=-n^{-1/2}\Lambda_{\submin} \bu\bv^{\top}$. Here $\bu\in\RR^{d_1}$ and $\bv\in\RR^{d_2}$ have unit norms. In this case, the tensor signal strength is $\Lambda_{\submin}>0$. Moreover, the individual signal strength is $\lambda=n^{-1/2}\Lambda_{\submin}$ and separation strength is $\Delta=2n^{-1/2}\Lambda_{\submin}$, i.e., the two quantities are at the same order.   Then the observations can be re-written as 
\begin{equation}\label{eq:rank-one-model}
\bX_i=s_i^{\ast}(n^{-1/2}\Lambda_{\submin} \bu\bv^{\top})+\bE_i,\quad \forall i=1,\cdots,n,
\end{equation}
where $s_i^{\ast}=1$ if $\bX_i$ is sampled from $\calN(\bM_1, \bI_{d_1}\otimes \bI_{d_2})$ and $s_i^{\ast}=-1$ if $\bX_i$ is sampled from $\calN(\bM_2, \bI_{d_1}\otimes \bI_{d_2})$. Note that the rank-one model (\ref{eq:rank-one-model}) is no more difficult than the general K-component case but it suffices for our purpose. The null hypothesis $\QQ_n$ corresponds to the case $\Lambda_{\submin}=0$, i.e., all observations are pure noise. Clearly, the difficulty level of distinguishing $\QQ_n$ and $\PP_n$ is characterized by signal strength $\Lambda_{\submin}$ in eq. (\ref{eq:rank-one-model}). 
Conjecture~\ref{conjecture:comp-hardness} requires the calculation of $\|L_n^{\leq D}\|^2$, which is extremely difficult for generally fixed singular vectors $\bu, \bv$ and deterministic latent labels $\bs^{\ast}$. A prior distribution simplifies the calculation. Finally, our null and alternative hypothesis are formally defined as follows. 

\begin{definition}[Null and alternative hypothesis]\label{def:comp-prior}\
	\begin{itemize}
			\item Under $\QQ_n$, we observe $n$ matrices $\bX_1,\cdots,\bX_n$ generated i.i.d. from \eqref{eq:rank-one-model} with $\Lambda_{\submin}=0$. Equivalently, it means that each $\bX_i$ has i.i.d. standard normal entries.
		\item Under $\PP_n:=\PP_n^{\Lambda_{\submin}}$, we observe $n$ matrices $\bX_1,\cdots,\bX_n$ generated i.i.d. from \eqref{eq:rank-one-model} with $\Lambda_{\submin}>0$, and moreover, each coordinate of $\bu$ and $\bv$ independently uniformly take values from $\{\pm d_1^{-1/2}\}$ and  $\{\pm d_2^{-1/2}\}$, respectively, and the entries of $\bs^\ast$ are independent Rademacher random variables, i.e., taking $\pm1$ with equal probabilities.
	\end{itemize}
\end{definition}

\begin{theorem}\label{thm:comp-lower-bound}
	Consider $\QQ_n$ and $\PP_n$ in Definition \ref{def:comp-prior}. If $\Lambda_{\submin}=o\left(d^{1/2}n^{1/4}\right )$ as $n\to\infty$, then $\op{L_n^{\le D}}=1+o(1)$.	
\end{theorem}

The proof of Theorem~\ref{thm:comp-lower-bound} can be found in \cite{lyu2022optimal}. If Conjecture~\ref{conjecture:comp-hardness} is true, Theorem~\ref{thm:comp-lower-bound} implies that $\QQ_n$ and $\PP_n^{\Lambda_{\submin}}$ are statistically indistinguishable by polynomial-time algorithms as long as the signal strength $\Lambda_{\submin}=o\left(d^{1/2}n^{1/4}\right )$. We now establish the connection of testing the hypothesis to the clustering problem under two-component symmetric LrMM (\ref{eq:rank-one-model}).

For any fixed $\Lambda_{\submin}>0$, define the parameter space of interest by
\begin{align*}
\wt\Omega_{\Lambda_{\submin}}&\equiv\wt\Omega(\Lambda_{\submin},d_1,d_2,n)\\
&=\Big\{(\bM,\bs):~\bM=n^{-1/2}\Lambda^\prime_{\submin}\bu\bv^\top, \bu\in\RR^{d_1}, \bv\in\RR^{d_2}, \bs\in\{\pm 1\}^n, |{\bf 1}^{\top}\bs|\leq n/2, \Lambda^\prime_{\submin}\ge \Lambda_{\submin}\Big\}	
\end{align*}
By Chernoff bound, with probability at least $1-e^{-c_0n}$ where $c_0>0$ is an absolute constant, the i.i.d. observations $\bX_1,\cdots,\bX_n$ generated by $\PP_n^{\Lambda_{\submin}}$ satisfy the rank-one LrMM (\ref{eq:rank-one-model}) with parameters $(\bM,\bs)\in \tilde{\Omega}_{\Lambda_{\submin}}$. The following theorem tells that if consistent clustering is possible for LrMM, so is for distinguishing the hypothesis in Definition~\ref{def:comp-prior}.

\begin{theorem}\label{thm:reduction}
	Suppose there exists a clustering algorithm $\hat\bs_{\textsf{comp}}:\RR^{d_1\times d_2\times n}\mapsto \{\pm 1\}^n$ for LrMM (\ref{eq:rank-one-model}) with runtime $poly(n,d)$ that is consistent under the sequence of signal strength $\big\{\Lambda_{\submin}^{(n)}\big\}_{n\geq 1}$ in the sense that there exists a sequence $\{(\delta_n,\zeta_n)\}_{n\geq 1}\to 0$ such that for all large $n$,
\begin{equation}\label{eq:reduction-assump}
		\sup_{(\bM,\bs^\ast)\in \wt\Omega_{\Lambda_{\submin}^{(n)}}}\PP\left(n^{-1}\cdot h_{\textsf{c}}(\hat \bs_{\textsf{comp}},\bs^\ast)>\delta_n\right)\le \zeta_n
\end{equation}
If the signal strength satisfies $\Lambda_{\submin}^{(n)}\geq C_0(1+\epsilon^{-2})^{1/2}d^{1/2}$ with some absolute constant $C_0>0$ and $\epsilon\in(0,1)$, 	then there exists a test $\phi_n:\RR^{d_1\times d_2\times n}\mapsto\{0,1\}$ with runtime $poly(n,d)$ that consistently distinguishes $\PP_n^{\Lambda_{\submin}^{(n)}}$ from $\QQ_n$ so that
$$
\EE_{Q_n}[\phi_n(\bcalX)]+\sup_{((1-\epsilon)\bM,\bs^\ast)\in \wt\Omega_{\Lambda_{\submin}^{(n)}}}\EE_{(\bM,\bs^\ast)}[1-\phi_n(\bcalX)]\to 0,\quad {\rm as}\ n,d\to \infty.
$$
\end{theorem}

Essentially, Theorem~\ref{thm:reduction} only needs a signal strength $\Lambda_{\submin}\gg d^{1/2}$ to successfully reduce a polynomial-time clustering algorithm to a polynomial-time hypothesis test. Based on Conjecture \ref{conjecture:comp-hardness}, a combination of Theorem \ref{thm:comp-lower-bound} and Theorem \ref{thm:reduction} implies the following result, whose proof is straightfoward and hence omitted.

\begin{corollary}
Suppose Conjecture \ref{conjecture:comp-hardness} holds for $\QQ_n$ and $\PP_n$ in Definition \ref{def:comp-prior}. If the signal strength $\Lambda_{\submin}^{(n)}=o(d^{1/2}n^{1/4})$,  then for any polynomial-time clustering algorithm $\hat \bs_{\textsf{comp}}$, there exist absolute constants $\delta,\zeta>0$ such that
$$
\sup_{(\bM,\bs^\ast)\in \wt\Omega_{\Lambda_{\submin}^{(n)}}}\PP\left(n^{-1}\cdot h_{\textsf{c}}(\hat \bs_{\textsf{comp}},\bs^\ast)>\delta\right)\ge \zeta
$$
as $n\to\infty$. 
\end{corollary}
It is worth pointing out that even though the signal strength $\Lambda_{\submin}=o(d^{1/2}n^{1/4})$, the separation strength $\Delta=2n^{-1/2}\Lambda_{\submin}$ can still be much larger than $d^{1/2}n^{-1/2}$ that is required by Theorem~\ref{thm:main}. This suggests that if signal strength is not strong, consistent clustering by polynomial-time algorithms is still impossible even though the separation strength is very strong.

\section{Relaxing the Signal Strength Condition}\label{sec:Weak-SNR}
Our main theorem in Section~\ref{sec:main} imposes a strong signal strength condition on {\it all} the population center matrices, i.e., $\Lambda_{\submin}$ is lower bounded by $\Omega(d^{1/2}n^{1/4})$, or equivalently, $\lambda$ is lower bounded by $\Omega(d^{1/2}n^{-1/4})$.  While evidences in Section~\ref{sec:comp} show that this condition might be necessary for the two-component symmetric case if only polynomial-time algorithms are sought, this condition appears flawed in the general asymmetric case. This section aims to relax the signal strength condition in the sense that one population center matrix is allowed to be arbitrarily smaller (in spectral norm) than $d^{1/2}n^{-1/4}$, in which case \eqref{cond:Lambda} might fail. 

To simplify the narrative, we focus on the two-component LrMM, i.e., $K=2$ in model (\ref{eq:LrMM}), whose population center matrices are denoted by $\bM_1$ and $\bM_2$, respectively. However, it is straightforward to extend our discussion to the general case.  For $K=2$, it is more intuitive and convenient to express everything in terms of individual signal strength $\bM_1$ and $\bM_2$ instead of the tensor signal strength $\Lambda_{\min }$, even though they are equivalent\footnote{Alternatively, we can impose condition on $\min_{j=1,2} \{\sigma_{\submin }(\bcalM_j)\}$, where $[\bcalM_j]_{\cdot\cdot i}=\II(s_i^\ast=1)\bM_j$.}. Without loss of generality, we assume that $\|\bM_1\|_{\rm F}$ is large so that reliable estimation is possible,  and that $\|\bM_2\|_{\rm F}$ is small so that reliable estimation is impossible. The following assumption is made to clarify this further.
\begin{assumption}\label{assump:relaxed-lloyd-SNR}
There exists a small constant $c>0$ such that 
$$\sigma_1(\bM_2)\le  c\alpha^{-1/2}\left(\sqrt{\frac{d}{n}}+\kappa_0^{-1}\right),$$
and 
$$\frac{\sigma^2_{r_1}(\bM_1)}{\alpha^{-1}(\kappa_0^2\vee r_1)\left(\frac{d}{n}+1\right)} \rightarrow \infty$$

where $\kappa_0$, with slight abuse of notation,  is the condition number of $\bM_1$.
\end{assumption}

If $\kappa_0,\alpha=O(1)$, Assumption \ref{assump:relaxed-lloyd-SNR} can be recasted as $\sigma_{r_1}(\bM_1)\gg d^{1/2}n^{-1/2}+1$ and $\sigma_1(\bM_2)\leq c(d^{1/2}n^{-1/2}+1)$. Note that Assumption~\ref{assump:relaxed-lloyd-SNR} puts no lower bound on $\sigma_1(\bM_2)$. In the extreme case, $\sigma_1(\bM_2)$ is allowed to be zero and 
consistent estimation of $\bM_2$ is unavailable even if the true labels are revealed.  
Assumption~\ref{assump:relaxed-lloyd-SNR} already implies that $\Delta\gg \left(d^{1/2}n^{-1/2}+1\right)$ if the ranks $r_1, r_2$ are both upper bounded by $O(1)$, matching the separation condition \eqref{cond:min-sep} in Theorem \ref{thm:main}. Intuitively, although clustering shall becomes easier as the constant $c$ in Assumption~\ref{assump:relaxed-lloyd-SNR} decreases, this cannot be verified by Theorem~\ref{thm:main} where the signal strength condition (\ref{cond:Lambda}) fails. 

Under Assumption~\ref{assump:relaxed-lloyd-SNR}, it is generally pointless to compute the center matrix $\hat\bM_2$ by SVD in Lloyd's algorithm since $\bM_2$ cannot be reliably estimated. Moreover, the SVD procedure complicates the subsequent theoretical analysis of Lloyd's algorithm. Instead of estimating $\bM_2$ via SVD, we opt to a trivial estimate by setting $\hat\bM_2\supt={\bf 0}$. The detailed steps are enumerated in Algorithm~\ref{alg:weak-SNR-Lloyd}, whose theoretical performance is guaranteed by Theorem~\ref{thm:weak-SNR}.\footnote{We remark that the low-rankness assumption for $\bM_2$ in Theorem \ref{thm:weak-SNR} is not essential, which can be dropped by instead requiring ${\sqrt{r_1}\sigma_{r_1}(\bM_1)}/\fro{\bM_2}\to\infty$.}
\begin{algorithm}
\caption{Low-rank Lloyd's Algorithm under Relaxed SNR Assumption~\ref{assump:relaxed-lloyd-SNR} (rlr-Lloyd)}\label{alg:weak-SNR-Lloyd}
\begin{algorithmic}
	\STATE{\textbf{Input}: Observations: $\bX_1,\cdots,\bX_n\in\RR^{d_1\times d_2}$ where $\bX_i=\bM_{s_i^{\ast}}+\bE_i$ and $s_i^{\ast}\in\{1,2\}$, initial estimate $\hat\bs^{(0)}$, ranks $r_1,r_2$.}

	\FOR{$t = 1,\ldots,T$}
	\STATE{For each $k=1,2$:
	$$
	\hat\bM_k^{(t)}\leftarrow \textrm{best rank-}r_k \textrm{ approximation of } \bar {\bX}_k(\hat\bs^{(t-1)}):=\frac{\sum_{i=1}^n \ind{\hat s_{i}\suptt=k}\bX_i}{\sum_{i=1}^n \ind{\hat s_{i}\suptt=k}}
	$$ }
	\STATE{Set $\hat\bM_2^{(t)}\leftarrow {\bf 0}$ if $\sigma_{1}(\hat\bM_2^{(t)})<\sigma_1(\hat\bM_1^{(t)})$; or set $\hat\bM_1^{(t)}\leftarrow\hat\bM_2^{(t)}$, $\hat\bM_2^{(t)}\leftarrow {\bf 0}$ if $\sigma_{1}(\hat\bM_2^{(t)})>\sigma_1(\hat\bM_1^{(t)})$.}
	\STATE{Re-label by  setting, for each $i\in[n]$:
	$$
	\hat s_i^{(t)}\leftarrow \underset{k\in [2]}{\arg\min}\ \|\bX_i-\hat\bM_k\supt\|_{\rm F}^2
	$$}
	\ENDFOR
	\STATE{\textbf{Output}: $\hat\bs=\hat\bs^{(T)}$}
\end{algorithmic}
\end{algorithm}

\begin{theorem}\label{thm:weak-SNR}
Suppose Assumption \ref{assump:relaxed-lloyd-SNR} holds and $d\ge C_0\log K$ for some absolute constant $C_0>0$. Assume $\hat\bs^{(0)}$ satisfies
\begin{equation}\label{relaxed-init-cond}
		n^{-1}\cdot \ell_{\textsf{c}}(\hat\bs^{(0)},\bs^\ast)= o\left(\frac{\alpha }{\kappa_0^2}\Delta^2\right)
	\end{equation}
Furthermore, if $\sqrt{\frac{r_1}{r_2}}\cdot \frac{\sigma_{r_1}(\bM_1)}{\sigma_1(\bM_2)}\to\infty$, then we have 
$$
n^{-1}\cdot h_{\textsf{c}}(\hat\bs\suptt, \bs^{\ast})\leq \exp\left(-\big(1-o(1)\big)\frac{\Delta^2}{8}\right)+\frac{1}{2^t}
$$
with probability at least $1-\exp(-\Delta)-\exp\big(-c_0d\big)$ for a small absolute constant $c_0>0$. 
\end{theorem}

To ensure  a consistent $\hat\bs^{(0)}$ satisfying \eqref{relaxed-init-cond}, we use a modified version of the tensor initialization discussed in Section \ref{sec:main-init}. The original spectral initialization can be mis-leading if a rank $r_{\bU}$ larger than $r_1$ is adopted. For our purpose, only the top-$r_1$ singular vectors are taken during spectral initialization, i.e., effort is made only for estimating $\bM_1$ whose left and right singular vectors are denoted by $\bU_1$ and $\bV_1$, respectively. See Algorithm \ref{alg:relaxed-initial} for further algorithmic details and Theorem \ref{thm:weak-SNR-init} for theoretical guarantees. 

\begin{algorithm}
\caption{Tensor-based Spectral Initialization Under Relaxed SNR Assumption (rTS-Init)}\label{alg:relaxed-initial}
\begin{algorithmic}
	\STATE{\textbf{Input}: observations: $\bX_1,\cdots,\bX_n\in\RR^{d_1\times d_2}$ where $\bX_i=\bM_{s_i^{\ast}}+\bE_i$ and $s_i^{\ast}\in\{1,2\}$; or a tensor $\bcalX\in\RR^{d_1\times d_2\times n}$ by concatenating the matrix observations slice by slice, ranks $r_1$.}
	\STATE{Spectral initialization:}
	\begin{itemize}
	\item[1.] Obtain the estimated singular vectors $\hat\bU_1$ and $\hat\bV_1$ by applying HOSVD to the tensor $\bcalX$ in mode-1 and mode-2 matricizations with rank $r_1$.
	\item[2.] Project $\bcalX$ onto the column space of $\hat\bU_1$ and $\hat\bV_1$ by $\hat\bcalG:=\bcalX\times_1 \hat\bU_1\hat \bU_1^{\top}\times_2 \hat\bV_1\hat\bV_1^{\top}$
	\item[3.] Apply K-means on the rows of $\hat\bG:=\scrM_3(\hat\bcalG)\in\RR^{n\times d_1d_2}$ and obtain the initial clustering by 
	$$
	(\hat\bs^{(0)}, \{\hat \bM_1^{(0)}, \hat\bM_2^{(0)}\}):=\underset{\bs\in[2]^n; \bM_1, \bM_2\in\RR^{d_1\times d_2}}{\arg\min}\ \sum_{i=1}^n \big\|[\hat\bG]_{i\cdot}-vec(\bM_{s_i}) \big\|^2
	$$
	\end{itemize}
	\STATE{\textbf{Output}: $\hat\bs^{(0)}$}
\end{algorithmic}
\end{algorithm}

\begin{theorem}\label{thm:weak-SNR-init}
Let $\hat\bs^{(0)}$ be the initial clustering output by Algorithm~\ref{alg:relaxed-initial}. Suppose there exists constant $c,C_0>0$ and  large constant $C>1$ such that $n/\kappa_0^4\ge C$,
$$\sigma_{r_1}(\bM_1)\ge C\alpha^{-1/2}\frac{d^{1/2}}{n^{1/4}},\quad \sigma_1(\bM_2)\le  C^{-1}\kappa_0^{-1}\frac{d^{1/2}}{n^{1/4}},$$
then we get, with probability at least $1-\exp(-cd)$, that
$$
 n^{-1}\cdot h_{\textsf{c}}(\hat\bs^{(0)},\bs^\ast)\le \frac{C_0}{\Delta^2}\left(\frac{dr_1}{n}+1\right)\quad\text{and}\quad  n^{-1}\cdot \ell_{\textsf{c}}(\hat\bs^{(0)},\bs^\ast)\le {C_0}\left(\frac{dr_1}{n}+1\right).
$$
Furthermore, if $ n/\kappa_0^4\rightarrow \infty$ and $\alpha\Delta^2/\kappa_0^2\rightarrow \infty$, with probability at least $1-\exp(-cd)$ we have that
$$
n^{-1}\cdot h_{\textsf{c}}(\hat\bs^{(0)},\bs^\ast)= o\left(\frac{\alpha }{\kappa_0^2}\right)\quad\text{and}\quad n^{-1}\cdot \ell_{\textsf{c}}(\hat\bs^{(0)},\bs^\ast)= o\left(\frac{\alpha }{\kappa_0^2}\Delta^2\right).
$$
\end{theorem}

Theorem \ref{thm:weak-SNR-init} serves as a counterpart of Theorem \ref{thm:spec-initialization}, with the distinction that we express the conditions in terms of $\sigma_{r_1}(\bM_1)$ and $\sigma_1(\bM_2)$. Notably, the threshold $d^{1/2}n^{-1/4}$  illuminates the disparity between statistical and computational aspects in the presence of low-rankness structure as discussed in Section \ref{sec:comp}. We emphasize that the gap arises solely due to the initialization procedure similar to the case in Section \ref{sec:main-init}.  Within our framework, Assumption  \ref{assump:relaxed-lloyd-SNR} together with good initializer $\hat\bs^{(0)}$ suffices to guarantee the statistical optimality of Algorithm \ref{alg:weak-SNR-Lloyd} under relaxed a signal strength condition and minimal requirement on the separation strength $\Delta$.

\section{Clustering versus Estimation}\label{sec:comparison}
\cite{lyu2022optimal} investigated the minimax optimal estimation of latent low-rank matrices under two-component symmetric LrMM, which revealed multiple phase transitions and a statistical-to-computational gap. In this section, together with Theorem~\ref{thm:main} and \ref{thm:spec-initialization}, we discuss the differences between estimation and clustering.

\subsection{Example where clustering is more challenging}
For simplicity, we consider the rank-one symmetric two-component LrMM \eqref{eq:rank-one-model} with $d_1=d_2=d$, where the separation strength $\Delta\asymp n^{-1/2}\Lambda_{\min}$ and individual signal strength $\lambda=n^{-1/2}\Lambda_{\min}$ coincides up to a constant factor. To make comparison, in this section we consider  $\Lambda_{\submin}$ instead of $\lambda$. The minimax rate of estimating $\bM$ (up to a sign flip), established in \cite{lyu2022optimal}, is
\begin{equation}\label{eq:est_rate}
\inf_{\hat \bM}\sup_{(\bM, \bs^{\ast})\in\wt\Omega_{\Lambda_{\submin}}}\E  \min_{\eta=\pm 1}\fro{\hat\bM-\eta\bM}\asymp \min\left\{d^{1/2}{\Lambda^{-1}_{\submin}}+d^{1/2}n^{-1/2},n^{-1/2}\Lambda_{\submin}\right\}
\end{equation}
The above rate is achievable by the computationally NP-hard maximum likelihood estimator with almost no constraint on signal strength and by 
a computationally fast spectral-aggregation estimator under the regime of strong signal strength $\Lambda_{\submin}\gtrsim d^{1/2}n^{1/4}$. For a fair comparison, we focus on this computationally feasible regime. The phase transitions under this regime can be summarized as in Table \ref{table:est-phase}. 

\begin{table}[!htbp]
\centering
\resizebox{0.6\textwidth}{!}{
\begin{tabular}{cc|c}
\hline
Sample size                  & Signal strength & Minimax optimal estimation error \\ \hline
\multirow{2}{*}{$d^2\lesssim n$} & $d^{1/2}n^{1/4}\lesssim \Lambda_{\submin}\lesssim n^{1/2}$& $\frac{\sqrt{d}}{\Lambda_{\submin}}$ \tTop\tBot\\ \cline{2-3} 
                   & $\Lambda_{\submin}\gtrsim n^{1/2}$ &  $\sqrt{\frac{d}{n}}$ \tTop\tBot\\ \hline
$d^2\gg n $                  & $\Lambda_{\submin}\gtrsim d^{1/2}n^{1/4}$ & $\sqrt{\frac{d}{n}}$ \tTop\tBot\\ \hline
\end{tabular}
}
\caption{Phase transition in minimax optimal estimation for two-component symmetric LrMM under the regime of strong signal strength $\Lambda_{\submin}\gtrsim d^{1/2}n^{1/4}$. See (\ref{eq:est_rate}) and \cite{lyu2022optimal} for more details.}
\label{table:est-phase}
\end{table}

Without loss of generality, we assume the dimension $d\to\infty$ as $n\to\infty$. The case $d^2\gg n$ is referred to as the high-dimensional setting, and $d^2\lesssim n$ is called the low-dimensional setting. 
An estimator $\hat\bM$ is said {\it strongly consistent} if the relative estimation error $\|\hat\bM-\bM\|_{\rm F}\|\bM\|_{\rm F}^{-1}$ approaches to zero in expectation as $n\to\infty$. Table~\ref{table:est-phase} tells that strongly consistent estimation $\bM$ is {\it always} achievable as long as the signal strength is greater $d^{1/2}n^{1/4}$. A particularly interesting regime is $d^{1/2}n^{1/4}\lesssim \Lambda_{\submin}\lesssim n^{1/2}$. For instance, when $d^2=o(n)$, $\bM$ can still be consistently estimated even when the signal strength $\Lambda_{\submin}\to 0$ as $n\to\infty$.  

It is certainly not the case for clustering. Besides {\it consistent clustering} (see definition in Theorem~\ref{thm:reduction}), we say a clustering algorithm is {\it weakly efficient} if it can beat a random guess, but the mis-clustering error rate does not vanish as $n\to \infty$. When $d^2=o(n)$, Theorem~\ref{thm:minimax-lower-bound} dictates that even weakly efficient clustering is impossible, i.e., $\exp(-\Lambda_{\submin}^2/(2n))$ is at least $1/2$, if $\Lambda_{\submin}\le c_0n^{1/2}$ for some absolute constant $c_0>0$. However, the spectral aggregation estimator \citep{lyu2022optimal} can still consistently estimate the population center matrix $\bM$ in the aforesaid scenario. Moreover, by Theorem~\ref{thm:main}, consistent clustering even requires $\Lambda_{\submin}/n^{1/2}\to\infty$, which is much more stringent than that required by (strongly) consistent estimation. 

The differences of phase transitions in estimation and clustering are enumerated in Table ~\ref{table:est-clust-phase}. Basically, strongly consistent estimation is always possible as long as $\Lambda_{\submin}\gtrsim d^{1/2}n^{1/4}$. In contrast, weakly efficient clustering is possible only when $\Lambda_{\submin}\gtrsim n^{1/2}+d^{1/2}n^{1/4}$, and consistent clustering is possible only when $\Lambda_{\submin}\gtrsim d^{1/2}n^{1/4}$ and meanwhile $\Lambda_{\submin}\gg n^{1/2}$. Note that the gap between estimation and clustering is present only under the low-dimensional setting $n\gtrsim d^2$. The gap vanishes under the high-dimensional setting $d^2\gg n$, in which case the signal strength condition $\Lambda_{\submin}\gtrsim d^{1/2}n^{1/4}$ already implies $\Lambda_{\submin}\gg n^{1/2}$. 

\begin{table}[!tbp]
\centering
\resizebox{0.8\textwidth}{!}{
\begin{tabular}{c|cccc}
\hline
Sample size             & Signal strength & Consistent estimation & Weakly efficient clustering & Consistent clustering \\ \hline
\multirow{3}{*}{$d^2\lesssim n$} & $d^{1/2}{n^{1/4}}\lesssim \Lambda_{\submin}\lesssim n^{1/2}$& Possible & Impossible & Impossible  \tTop\tBot\\ \cline{2-5} 
                   & $n^{1/2}\lesssim \Lambda_{\submin}\lesssim n^{1/2} $ &  Possible & Possible & Impossible \tTop\tBot\\ \cline{2-5} 
                   & $\Lambda_{\submin}\gg n^{1/2}$ &   Possible & Possible & Possible \tTop\tBot\\ \hline
$d^2\gg n $                  & $\Lambda_{\submin}\gtrsim d^{1/2}{n^{1/4}}$ & Possible & Possible & Possible \tTop\tBot\\ \hline
\end{tabular}
}
\caption{The differences of phase transitions in estimation and clustering for two-component symmetric LrMM under the regime of strong signal strength $\Lambda_{\submin}\gtrsim d^{1/2}n^{1/4}$. Here $d^2\gg n$ is referred to as the high-dimensional setting, and $d^2\lesssim n$ as the low-dimensional setting.}
\label{table:est-clust-phase}
\end{table}

We collect these facts to convince that, at least for the two-component symmetric LrMM (\ref{eq:rank-one-model}), clustering is intrinsically more challenging than estimation. The same phenomenon also arises in GMM. See, e.g., \cite{wu2019randomly}.

\subsection{Example where estimation is more challenging}
While, generally, clustering is recognized as being more challenging than estimation, there are examples where clustering is easier than estimation. Similarly as in Section~\ref{sec:Weak-SNR}, consider the two-component LrMM with population center matrices $\bM_1$ and $\bM_2$ so that 
$$
\sigma_{r_1}(\bM_1)\geq C_1\Big(1+\frac{d^{1/2}}{n^{1/2}}+\frac{d^{1/2}}{n^{1/4}}\Big)\quad {\rm and}\quad \sigma_1(\bM_2)\leq C_1^{-1}\cdot\frac{d^{1/2}}{n^{1/2}}
$$
where $C_1>0$ is a large constant and, for simplicity, we assume $\kappa_0, \alpha, r_1, r_2=O(1)$. Observe that 
$$
\sqrt{\frac{r_1}{r_2}}\cdot \frac{\sigma_{r_1}(\bM_1)}{\sigma_1(\bM_2)}\geq \begin{cases}
C_1^2n^{1/4},& \textrm{ if } n\leq d^2;\\
C_1^2(n/d)^{1/2},& \textrm{ if } n>d^2;
\end{cases}
\quad \to \infty, \quad {\rm as }\ n\to\infty
$$
Moreover, 
$$
\Delta:=\|\bM_1-\bM_2\|_{\rm F}\gtrsim C_1\Big(1+\frac{d^{1/2}}{n^{1/4}}\Big) \to \infty
$$
if the constant $C_1>0$ diverges to infinity. Therefore, by Theorem~\ref{thm:weak-SNR}, if $C_1\to\infty$, our Algorithm~\ref{alg:weak-SNR-Lloyd} consistently cluster all observations. 

However, consistent estimation of the population center matrices is more challenging. Even if all the latent labels are correctly identified, estimation of $\bM_2$ is still impossible because of its weak signal strength. Indeed, the low-rank approximation to 
$$
\bar{\bX}_2(\bs^{\ast}):=\frac{1}{n_2^{\ast}}\sum_{i=1}^n \ind{s_i^{\ast}=2} \bX_i
$$
achieves the error rate (in expectation) $O(d^{1/2}n^{-1/2})$ and the relative error rate (in expectation) diverges to infinity as $C_1\to\infty$. Similarly, the trivial estimate by a zero matrix attains the relative error rate $1$ that never vanishes as $n\to\infty$. Consequently, a strongly consistent estimate of $\bM_2$ becomes impossible.

\section{Discussions}\label{sec:discussion}

\subsection{Estimation of  $r_{\bU}$, $r_{\bV}$, $K$  and $r_k$'s}
Our tensor-based spectral initialization method requires an input of ranks $r_{\bU}$, $r_{\bV}$ and the number of clusters $K$, which are usually unknown in practice. Under the decomposition (\ref{eq:M-dec2}), they constitute the Tucker ranks of tensor $\bcalM$. Several approaches are available to estimate the Tucker ranks for tensor PCA model. One typical approach \citep{jing2021community,cai2022generalized} is to check the scree plots \citep{cattell1966scree} of $\scrM_1(\bcalX)$, $\scrM_2(\bcalX)$ and $\scrM_3(\bcalX)$, respectively. Under a suitable signal strength condition as in Theorem~\ref{thm:spec-initialization}, the scree plots of $\scrM_1(\bcalX)$ and $\scrM_2(\bcalX)$ shall serve a reliable estimate of $r_{\bU}$ and $r_{\bV}$, respectively. However, we note that it is statistically more efficient to estimate $K$ by, instead, taking the scree plot of $\scrM_3(\bcalX\times_1 \hat\bU^{\top}\times_2 \hat \bV^{\top})$, where $\hat\bU$ and $\hat\bV$ are obtained in step 1 of Algorithm~\ref{alg:initial}. This additional spectral projection promotes further noise reduction as in Algorithm~\ref{alg:initial}. After obtaining $r_{\bU}$, $r_{\bV}$ and $K$, an initial clustering $\hat\bs^{(0)}$ can be attained by apply Algorithm~\ref{alg:initial}. Similarly, we then estimate the rank $r_k$ by the scree plot of the sample average of matrix observations whose initial labels are $k$. It provides a valid estimate as long as the initial clustering is sufficiently good. 
The aforementioned approach works nicely in 
real-world data applications. See Section~\ref{sec:numerical} for more details.

\subsection{Matrix observation with categorical entries}
Oftentimes, the matrix observations consist of categorical entries. For instance, the Malaria parasite gene networks (see Section~\ref{sec:num_malaria}) have binary entries (Bernoulli distribution); the 4D-scanning transmission electron microscopy \citep{han2022optimal} produces count-type entries (Poisson distribution). Our algorithms are still applicable and deliver appealing performance on, e.g., Malaria parasite gene networks dataset. Unfortunately, our theory can not directly cover those cases, although the noise are still sub-Gaussian. Without loss of generality, let us consider multi-layer binary networks and assume $\bX_i$ has Bernoulli entries. Then the entries of $\bX_i$ have an equal variance only when they have the same expectation, reducing the network to a trivial Erd\H{o}s-R\' enyi graph. Nevertheless, equal noise variance is crucial to establish Theorem~\ref{thm:spec-initialization}. Moreover, the techniques for proving Theorem~\ref{thm:main} are likely sub-optimal since the sub-Gaussian constant $\sigmasg$ is usually not sharp enough to characterize a Bernoulli random variable. We leave this to future works.

\section{Numerical Experiments and Real Data Applications}\label{sec:numerical}

\subsection{Numerical Experiments}
This section presents the empirical performance of lr-Lloyd's algorithm (Algorithm \ref{alg:Lloyd}) and its relaxed variant under weak SNR (Algorithm \ref{alg:weak-SNR-Lloyd}) referred to as the rlr-Lloyd's algorithm. Specifically, we focus on the algorithmic convergence and final clustering error. 

In the first simulation setting \textbf{S1}, we fix the dimension $d_1=d_2=50$ and sample size $n=200$. The latent labels $s_i^{\ast}$ are generated i.i.d. from the model (\ref{eq:mixture}) with equal mixing probabilities, i.e., $\pi_k=1/K$. All the presented results in \textbf{S1} are based on the average of $30$ independent trials. We test the convergence of Algorithm~\ref{alg:Lloyd} under both Gaussian (\textbf{S1-1}) and Bernoulli (\textbf{S1-2}) noise.

\begin{figure}[!h]
	\centering
	\begin{subfigure}[b]{.49\linewidth}
		\includegraphics[width=3in]{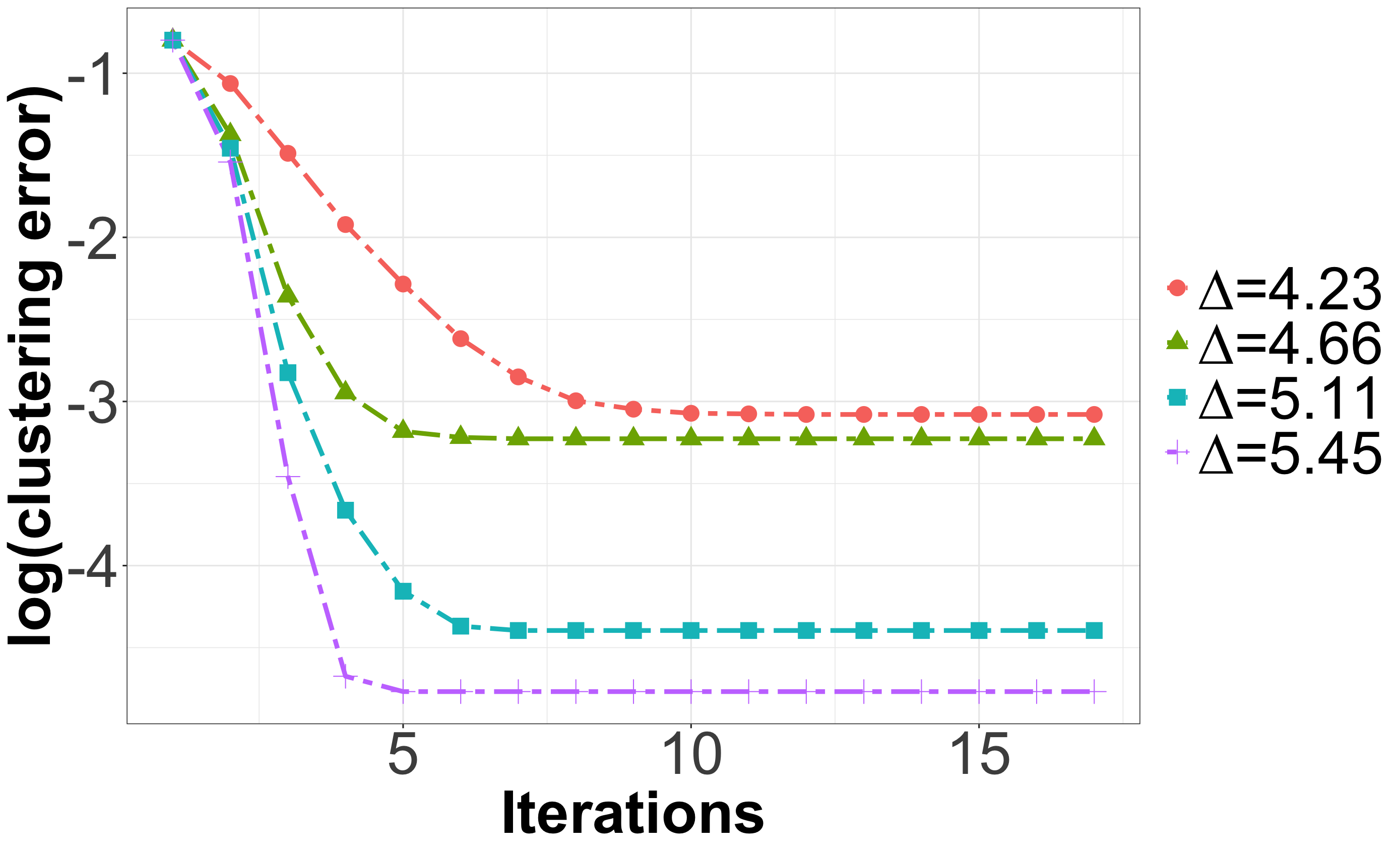}
		\caption{\tiny{Simulation \textbf{S1-1}: Log of clustering error  ($K=2$) with $\Delta$ varying under Gaussian noise.}}
	\end{subfigure}
	\begin{subfigure}[b]{.49\linewidth}
		\includegraphics[width=3in]{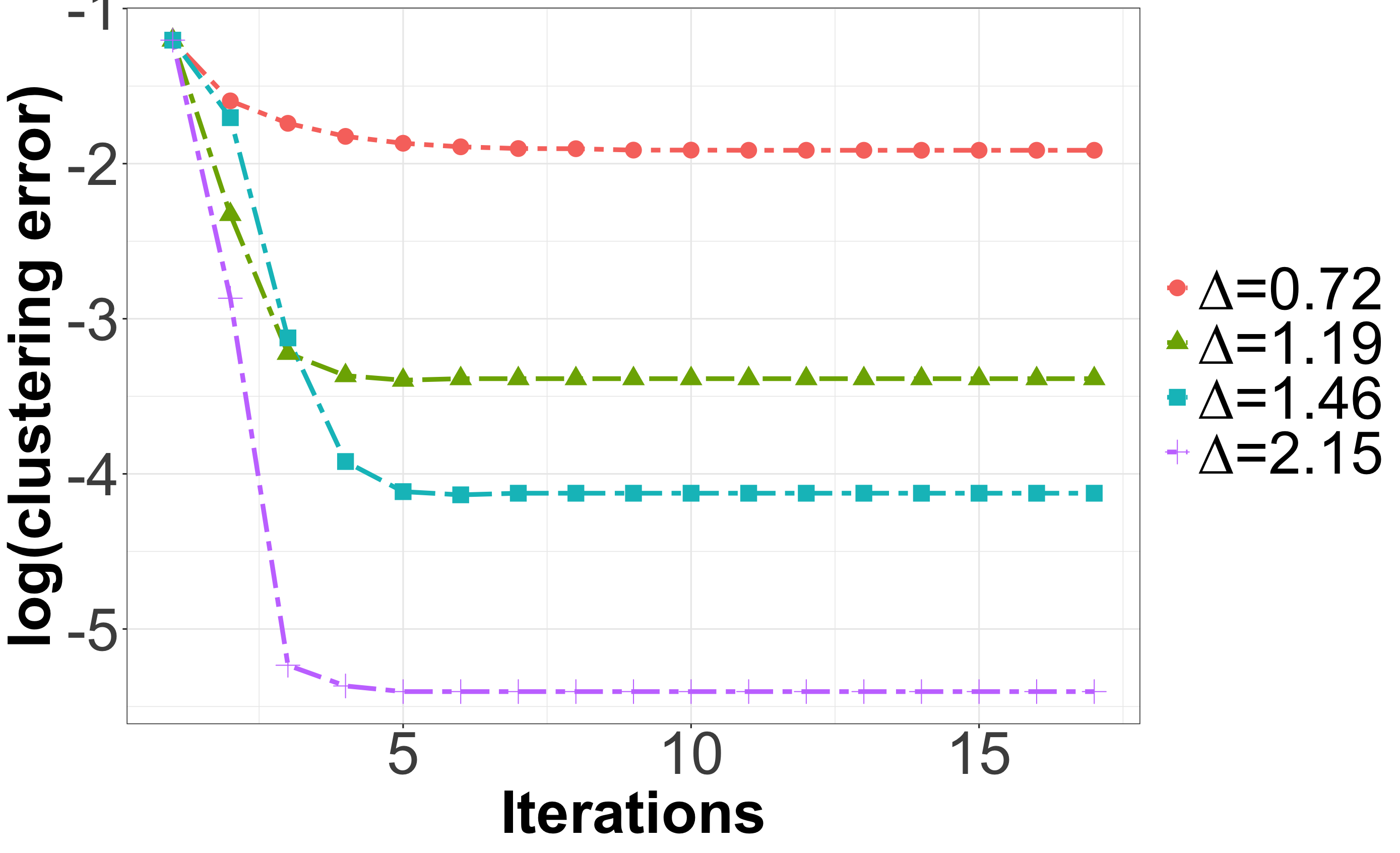}
		\caption{\tiny{Simulation \textbf{S1-2}  Log of clustering error  ($K=3$) with $\Delta$ varying under Bernoulli noise (MMSBM).}}
	\end{subfigure}
	\caption{(Convergence behavior of Algorithm~\ref{alg:Lloyd}) Log of clustering error with $\Delta$ varying under two scenarios: LrMM with Gaussian noise and MMSBM with Bernoulli noise.}
	\label{fig:S1}
\end{figure}

In \textbf{S1-1}, we set $K=2$, $r_1=r_2=2$ and standard Gaussian noise. The population center matrices $\bM_1$ and $\bM_2$ are generated in the following manner. For each $k=1,2$, we independently generate a $d_1\times d_2$ matrix with i.i.d. standard Gaussian entries and extract its top-$2$ left and right singular vectors as $\bU_k$ and $\bV_k$, respectively. The singular values are manually set as $\bSigma_k=\text{diag}\{1.2\lambda,\lambda\}$ for some fix $\lambda>0$. Then the population center matrices are constructed as $\bM_k=\bU_k\bSigma_k\bV_k^\top$. Our experiment tries four levels of signal strength $\lambda\in\{1.9, 2.1, 2.3, 2.5\}$. For each $\lambda$, the population center matrices are generated as above and the separation strength is recorded. The corresponding separation strength are $\Delta\in\{4.22, 4.66, 5.11, 5.45\}$. At each level of signal strength, the observations $\{\bX_i: i=1,\cdots,200\}$ are independently drawn from (\ref{eq:LrMM-1}) with the obtained center matrices $\bM_1$ and $\bM_2$. Here we focus on the convergence behavior of Lloyd's iterations of Algorithm~\ref{alg:Lloyd}, and thus a warm initial clustering $\hat\bs^{(0)}$ is provided before hand. The same initial clustering is used for all simulations and the initial clustering error is $n^{-1}h_{\textsf{c}}(\hat\bs^{(0)}, \bs^{\ast})=0.45$, i.e., slightly better than a random guess. Convergence of Algorithm~\ref{alg:Lloyd} under four levels of signal strength (or, correspondingly, separation strength) is displayed in the left plot of Figure~\ref{fig:S1}. The decreasing of log of clustering error is linear in first few iterations, as expected by our Theorem~\ref{thm:main}. The algorithm converges fast and the final clustering error is reflected by the separation strength $\Delta$. It is worth pointing out that Figure~\ref{fig:S1}(a) also shows that Algorithm~\ref{alg:Lloyd} converges faster when $\Delta$ becomes larger. While this cannot be directly concluded from Theorem~\ref{thm:main}, it can be easily verified by checking the proof.

In \textbf{S1-2}, we test the effectiveness of Algorithm \ref{alg:Lloyd} under non-Gaussian and non-i.i.d. noise. In particular, we consider the mixture multi-layer stochastic block model (MMSBM) introduced in \cite{jing2021community} \footnote{We emphasize that our Theorem~\ref{thm:main} is not directly applicable to MMSBM due to non-i.i.d. noise.}. We set the number of clusters $K=3$. For each $k=1,2,3$, the $k$-th SBM is associated with a connection probability matrix $\bB_k\in[0, 1]^{K\times K}$ and a membership matrix $\bZ_k\in\{0,1\}^{d\times K}$, which are set as $\bB_k:=\bar{p}_k\cdot \bI_K+\bar p_k/2\cdot (\mathbf{1}_K\mathbf{1}_K^\top-\bI_K)$ with $\bar p_k=\bar p\cdot k/K$ and $\bZ_k(i,:)=\be_{s^*_i}$, respectively. Thus each SBM has three cluster of nodes and the population center matrices are $\bM_k=\bZ_k\bB_k\bZ_k^{\top}\in [0,1]^{d\times d}$. Conditioned on the latent label $\bs_i^{\ast}$, the $i$-th observation $\bX_i$ is sampled from SBM($\bZ_{s_i^{\ast}}, \bB_{s_i^{\ast}}$), namely, $\bX_i(j_1,j_2)\sim  \text{Bernoulli}(\bM_{s^*_i}(j_1,j_2))$ and $\bX_i(j_2,j_1)=\bX_i(j_1,j_2)$ for $1\le j_1<j_2\le d$.  Note that $\bX_i$ is symmetric because the network is undirected. We manually set the diagonal entries of $\bX_i$ to zeros so that no self-loop is allowed in the observed network. Clearly, the entry-wise variances of $\bX_i$ are not necessarily equal. Under the above MMSBM, the signal strength and separation strength are characterized by sparsity level $\bar p$. Four sparsity levels $\bar p\in\{0.05, 0.08, 0.10, 0.15\}$ are studied so that the corresponding separation strength are $\Delta\in\{0.75, 1.19, 1.46, 2.15\}$. Similarly, a fixed good initial clustering $\hat\bs^{(0)}$ is used for all simulations and the initial clustering error is $n^{-1}h_{\textsf{c}}(\hat\bs^{(0)}, \bs^{\ast})=0.3$. Convergence behavior of Algorithm~\ref{alg:Lloyd} is displayed in the right plot of Figure~\ref{fig:S1}. Still, Lloyd's iterations converges fast and the final clustering error is decided by the separation strength $\Delta$.

In the second simulation setting \textbf{S2}, we aim to compare the final clustering error of vanilla Lloyd's algorithm and our low-rank Lloyd's algorithm. The dimensions are varied at two cases $d_1=d_2\in\{50, 100\}$, sample size is set as $n\in\{100, 200\}$, number of clusters $K=2$ and ranks $r_1=r_2=3$. The latent labels are generated as in \textbf{S1}. For each $d_1$ and $n$, the simulation is repeated for 100 times and their average clustering error rate is reported. 

In \textbf{S2-1}, the population center matrices $\bM_1$ and $\bM_2$ are constructed such that they share identical singular spaces. More exactly, we extract singular vectors $\bU_1$, $\bV_1$ and singular value matrix $\bSigma_1$ as is done in \textbf{S1-1}. Then the population center matrices are set as $\bM_1=\bU_1\bSigma_1\bV_1^{\top}$ and $\bM_2=\bU_1(\bSigma_1+\text{diag}\{\Delta/3,\Delta/3,\Delta/3\})\bV_1^\top$. Here the signal strength is fixed at $\lambda=10$ and the separation parameter is chosen from $\Delta\in\{1,5,10\}$. The final clustering error and its standard error by four methods are reported in the upper half of  Table~\ref{table:S2}. 
Noted that the initialization of ``vec-Lloyd" in \cite{lu2016statistical} is attained by spectral clustering on $\scrM_3(\bcalX)$. 
We observe that the clustering errors of four methods all decrease as $\Delta$ increases. However, lr-Lloyd initialized by Algorithm~\ref{alg:initial} achieves a much smaller clustering error compared with other methods. This is due to the fact that our proposed tensor-based spectral initialization is capable to capture the low-rank signal whereas both spectral clustering and naive K-means on $\scrM_3(\bcalX)$ ignores the low-rank structure in the other two modes of $\bcalM$. As a result, all the other three methods perform almost the same under current setting. Lastly, the bold-font column in Table~\ref{table:S2} confirms Theorem~\ref{thm:main} in that the clustering error achieved by TS-init initialized lr-Lloyd algorithm is only determined by $\Delta$ regardless of the dimension $d_1, d_2$ or the sample size $n$. 

In \textbf{S2-2}, the singular vectors of $\bM_1$ and $\bM_2$ are generated exactly the same as in \textbf{S1-1}. The singular values of $\bM_1$ and $\bM_2$ are set as $\bSigma_1={\rm diag}(1.2\lambda, 1.1\lambda, \lambda)$ and $\bSigma_2={\rm diag}(0.36, 0.33, 0.30)$, respectively. Then $\sigma_{\submin}(\bM_1)=\lambda$ and $\sigma_1(\bM_2)=0.36$. Here $\lambda$ is varied at $\{1.9, 2.2, 2.5\}$ for the case $d_1=d_2=50$ and $\{2.7, 3.0, 3.3\}$ for the case $d_1=d_2=100$. Consequently, the signal strength of $\bM_2$ is much smaller than $\bM_1$ that corresponds to the weak SNR setting in Section~\ref{sec:Weak-SNR}, and we test the performance of the relaxed lr-Lloyd's algorithm (Algorithm~\ref{alg:weak-SNR-Lloyd}). The results are reported in the lower half of Table~\ref{table:S2}. Clearly, rlr-Lloyd's algorithm outperforms the vanilla Lloyd's algorithm (i.e., the vectorized version). In certain cases, the vanilla Lloyd's algorithm merely beats a random guess whereas the rlr-Lloyd's algorithm almost achieves zero clustering error. We also observe that rlr-Lloyd's algorithm still performs nicely if initialized by K-means on $\scrM_3(\bcalX)$.

\begin{table}[!h]
\centering
\resizebox{\textwidth}{!}{%
\begin{tabular}{c|c|c|cc|cccc}
\hline
Setting & $d_1=d_2$& $n$ &  $\lambda$ & $\Delta$&  \begin{tabular}{@{}c@{}}vec-Lloyd \\  \citep{lu2016statistical}\end{tabular}   & \begin{tabular}{@{}c@{}}lr-Lloyd initialized by \\ TS-Init (Algorithm \ref{alg:initial})\end{tabular}  & \begin{tabular}{@{}c@{}}vec-Lloyd initialized by \\ K-means on $\scrM_3(\bcalX)$ \end{tabular} &\begin{tabular}{@{}c@{}}lr-Lloyd initialized by \\ K-means on $\scrM_3(\bcalX)$\end{tabular} \\ \hline
\multirow{12}{*}{\textbf{S2-1}} & \multirow{6}{*}{50} & \multirow{3}{*}{100} & 10 & 1 & 0.461 (0.032) & \textbf{0.401 (0.058)} & 0.462 (0.030) & 0.459 (0.031)\\
 & &  & 10 & 5 & 0.459 (0.033) & \textbf{0.163 (0.039)}  & 0.456 (0.033) & 0.452 (0.034) \\
 & &  & 10 & 10 & 0.458 (0.034) & \textbf{0.066 (0.025)}  & 0.441 (0.047) & 0.433 (0.054) \\
 \cline{3-9}
 & & \multirow{3}{*}{200} & 10 & 1 & 0.475 (0.019) & \textbf{0.398 (0.056)} & 0.469 (0.025) & 0.466 (0.025)\\
 &  &  & 10 & 5 & 0.473 (0.021) & \textbf{0.152 (0.027)} & 0.462 (0.027) & 0.450 (0.039) \\
 &  &  & 10 & 10 & 0.471 (0.022) & \textbf{0.063 (0.016)} & 0.437 (0.041) & 0.380 (0.082) \\
\cline{2-9}
& \multirow{6}{*}{100} & \multirow{3}{*}{100} & 10 & 1& 0.461 (0.028) & \textbf{0.391 (0.069)} & 0.460 (0.033) & 0.461 (0.033)  \\
&  & & 10 & 5 & 0.461 (0.029) & \textbf{0.157 (0.054)} & 0.455 (0.036) & 0.455 (0.036)\\
&  & & 10 & 10 & 0.460 (0.029) & \textbf{0.063 (0.026)} & 0.458 (0.034) & 0.456 (0.034)\\
\cline{3-9}
&  & \multirow{3}{*}{200} & 10 & 1 & 0.468 (0.023) & \textbf{0.390 (0.064)} & 0.469 (0.023) &  0.467 (0.023) \\
&  &  & 10 & 5 & 0.468 (0.024) & \textbf{0.147 (0.028)} &  0.469 (0.022) &  0.465 (0.026) \\
&  & & 10 & 10 & 0.467 (0.024) & \textbf{0.062 (0.017)} & 0.459 (0.030) & 0.451 (0.037)\\ \hline
Setting & $d_1=d_2$& $n$ &  $\sigma_\submin(\bM_1)$ & $\Delta$& \begin{tabular}{@{}c@{}}vec-Lloyd \\  \citep{lu2016statistical}\end{tabular}    & \begin{tabular}{@{}c@{}}rlr-Lloyd \\ (Algorithm \ref{alg:weak-SNR-Lloyd})\end{tabular}   & \begin{tabular}{@{}c@{}}vec-Lloyd initialized by \\ K-means on $\scrM_3(\bcalX)$ \end{tabular} &\begin{tabular}{@{}c@{}}rlr-Lloyd initialized by \\ K-means on $\scrM_3(\bcalX)$\end{tabular}   \\\hline
\multirow{12}{*}{\textbf{S2-2}} & \multirow{6}{*}{50} & \multirow{3}{*}{100} & 1.9 & 3.68 & 0.434 (0.052) & \textbf{0.314 (0.138)} & 0.418 (0.066) & 0.327 (0.129)\\
 & &  & 2.2 & 4.24 & 0.424 (0.061) & \textbf{0.134 (0.125)} & 0.385 (0.079) & 0.152 (0.138)\\
 & &  & 2.5 &  4.81 & 0.417 (0.068) & \textbf{0.041 (0.051)}  & 0.309 (0.103) & 0.055 (0.091) \\
 \cline{3-9}
 & & \multirow{3}{*}{200} & 1.9 & 3.68 & 0.433 (0.052) & \textbf{0.070 (0.020)} & 0.380 (0.070) &  0.072 (0.046)\\
 &  &  & 2.2 & 4.24 & 0.431 (0.054) & \textbf{0.057 (0.018)} & 0.351 (0.077) &  0.059 (0.048)\\
  &  &  & 2.5 & 4.81 & 0.424 (0.057) & \textbf{0.035 (0.015)} & 0.268 (0.088) & 0.033 (0.014) \\
\cline{2-9}
& \multirow{6}{*}{100} & \multirow{3}{*}{100} & 2.7 & 5.19 & 0.422 (0.056) & \textbf{0.300 (0.169)} & 0.416 (0.057) & 0.301 (0.164)  \\
&  & & 3 & 5.76 & 0.421 (0.059) & \textbf{0.131 (0.164)} & 0.390 (0.077) & 0.176 (0.181)\\
&  & & 3.3 & 6.33 & 0.426 (0.053) & \textbf{0.067 (0.139)} & 0.347 (0.086) & 0.065 (0.130)\\
\cline{3-9}
&  & \multirow{3}{*}{200} & 2.7 & 5.19 & 0.442 (0.040) & \textbf{0.019 (0.010)} & 0.395 (0.071) &  0.022 (0.037) \\
&  &  & 3 & 5.76 & 0.443 (0.041) & \textbf{0.008 (0.006)} &  0.301 (0.089) &  0.008 (0.007) \\
&  &  & 3.3 & 6.33 & 0.440 (0.043) & \textbf{0.003 (0.004)} &  0.190 (0.069) &  0.003 (0.004) \\
\hline
\end{tabular}%
}
\caption{Clustering error of lr-Lloyd (Algorithm~\ref{alg:Lloyd}) and rlr-Lloyd (Algorithm~\ref{alg:weak-SNR-Lloyd}) compared with vanilla Lloyd's algorithm \citep{lu2016statistical} on vectorized data  (vec-Lloyd). The number in brackets represents the standard error over $100$ trials.} 
\label{table:S2}
\end{table}

\subsection{Real Data Applications}
We now demonstrate the merits of our proposed low-rank Lloyd's (lr-Lloyd) algorithm on several real-world datasets and compare with existing methods. 

\subsubsection{BHL dataset}
The BHL (brain, heart and lung ) dataset\footnote{The dataset is publicly available at \url{https://www.ncbi.nlm.nih.gov/sites/GDSbrowser?acc=GDS1083}.}, which had been analyzed in \cite{mai2021doubly}, consists of $d_1=1124$ gene expression profiles of $n=27$  brain, heart, or lung tissues. Each tissue is measured repeatedly for $d_2=4$ times and hence the $i$th sample can be constructed as  $\bX_i\in\RR^{1124\times 4}$ for $i=1,\cdots,27$. Our aim is to correctly identify those $\bX_i$'s belonging to the same type of tissue, i.e., $K=3$. We apply Algorithm \ref{alg:Lloyd} together with an initial clustering $\hat\bs^{(0)}$ obtained by Algorithm \ref{alg:initial}  with $r_\bU=r_\bV=1$. These ranks are chosen based on the scree plots of $\scrM_1(\bcalX)$ and $\scrM_2(\bcalX)$. The final clustering error attained by lr-Lloyd's algorithm is $n^{-1}\cdot h_{\textsf{c}}(\hat\bs,\bs^\ast)=0.03704$. 
As shown in Table \ref{table:BHL}, our lr-Lloyd's algorithm performs the best among all the competitors\footnote{Note that all results except  lr-Lloyd are directly borrowed from \cite{mai2021doubly}, which use $\bX_i$'s after dimension reduction to a size of either $20\times 4$ or $30\times 4$, and we only report the better one here.} that are reported in \cite{mai2021doubly}. 

\begin{table}[]
\centering
{%
\begin{tabular}{ccccccccc}
\hline
 & {lr-Lloyd} & DEEM & K-means & SKM & DTC & TBM & EM & AFPF \\ \hline
Clustering error & \textbf{3.70} & 7.41 & 11.11 & 11.11 & 18.52 & 11.11 & 11.11 & 11.11 \\ \hline
\end{tabular}%
}
\caption{Clustering error on BHL dataset. SKM:  sparse K-means \citep{witten2010framework}; DTC: dynamic tensor clustering  \citep{sun2019dynamic}; TBM:  tensor block model (TBM) \citep{wang2019multiway}; EM: standard EM implemented in \cite{mai2021doubly}; AFPF: adaptive pairwise fusion penalized clustering \citep{guo2010pairwise}.} 
\label{table:BHL}
\end{table}

The improvement can be attributed to two reasons. First, DEEM in \cite{mai2021doubly} is designed based on EM algorithm targeted at Gaussian probability distribution, and hence they need to first perform multiple Kolmogorov-Smirnov tests to drop the columns not following Gaussian distribution, which might lead to potential information loss. In sharp contrast, their procedure is not necessary for our method, as the low-rank Lloyd's algorithm allows for sub-Gaussian noise. Secondly, our algorithm is more suitable for the specific structure of the data. Particularly, the population center matrices are expected to be rank-one as the columns of $\bX_i$ represent repeated measurements for the same sample. However, such planted structure is under-exploited in \cite{mai2021doubly} and others.

\subsubsection{EEG dataset}
The EEG dataset\footnote{The dataset is publicly available at \url{https://archive.ics.uci.edu/ml/datasets/EEG+Database}.} has been extensively studied by various statistical models \citep{li2010dimension, zhou2014regularized,hu2020matrix,huang2022robust}. The goal is to inspect EEG correlations of genetic predisposition to alcoholism. 
The data contains measurements which were sampled at $d_1=256$ Hz for $1$ second, from $d_2=64$ electrodes placed on each scalp of $n=122$ subjects. Each subject, either being \textit{alcoholic} or not, completed 120 trials under different stimuli. More detailed description of the dataset can be found in \cite{zhang1995event}. For our application, we average all the trials for each subject under single stimulus condition (S1) and two matched stimuli condition (S2), respectively, and  construct the data tensor as $\bcalX^{(S_1)}\in\RR^{256\times 64\times 122}$ (or $\bcalX^{(S_2)}\in\RR^{256\times 64\times 122}$) after standardization. Thus each subject is associated with a $256\times 64$ matrix, and we aim to cluster these subjects into $K=2$ groups, corresponding to alcholic group and control group. We apply rlr-Lloyd's algorithm (Algorithm \ref{alg:weak-SNR-Lloyd}) with $r_\bU=r_\bV=3$ and $r_1=2,r_2=1$. Here $r_{\bU}$ and $r_{\bV}$ are selected by the scree plot of $\scrM_1(\bcalX)$ and $\scrM_2(\bcalX)$, and $r_{1}$ and $r_{2}$ are tuned by interpreting the final outcomes. The clustering error of our method and competitors are shown in  Table \ref{table:EEG}. It is worth pointing out that our task of clustering is generally more challenging than classification, which has been investigated on the EEG dataset \citep{li2010dimension, zhou2014regularized,hu2020matrix,huang2022robust}. Those classification approaches often achieve lower {\it classification} error rates. As a faithful comparison, our rlr-Lloyd's algorithm enjoys a superior performance to its competitors in terms of {\it clustering} error rate and time complexity. 

\begin{figure}[h]
	\includegraphics[width=0.45\textwidth]{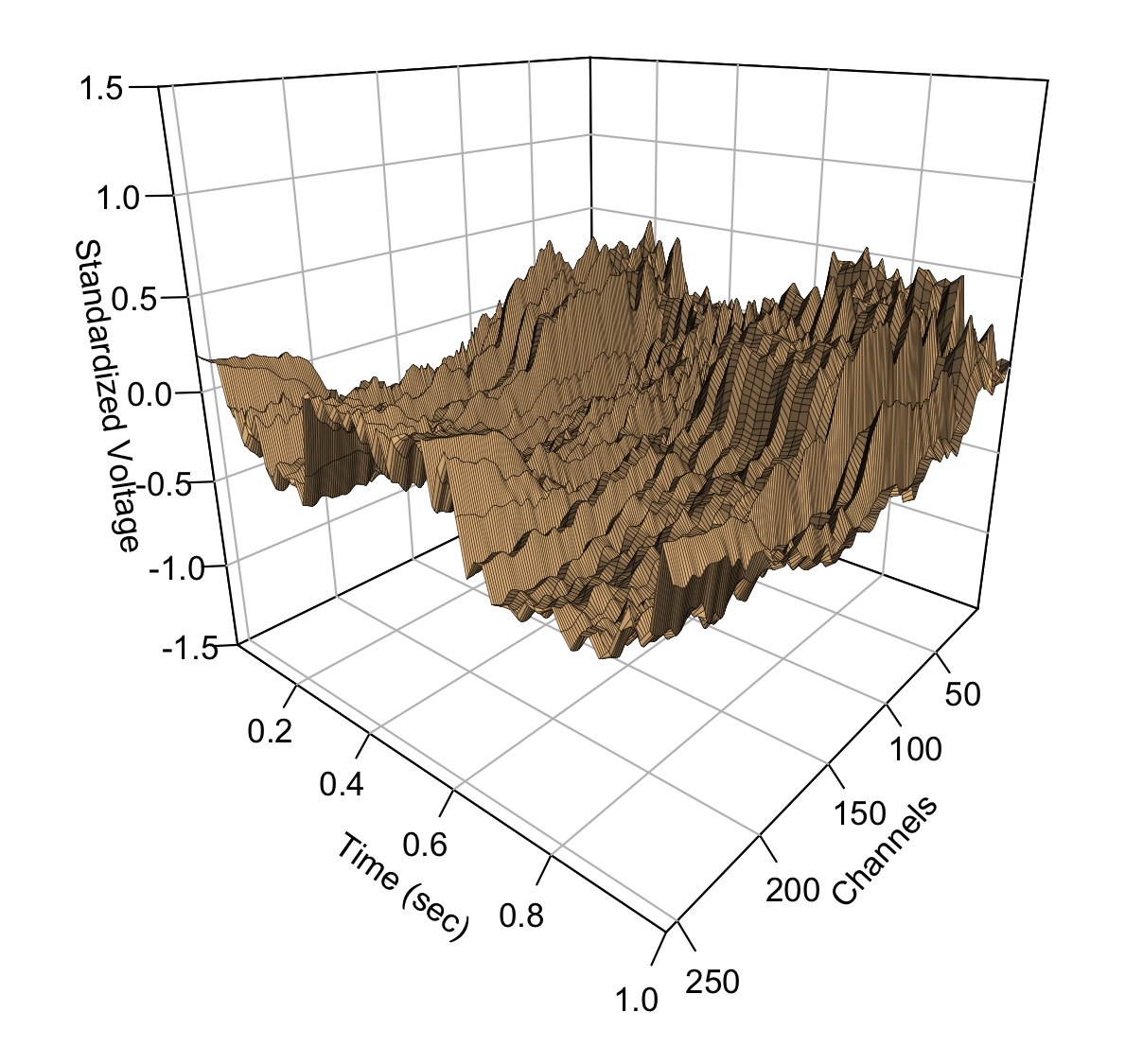}
	\includegraphics[width=0.45\textwidth]{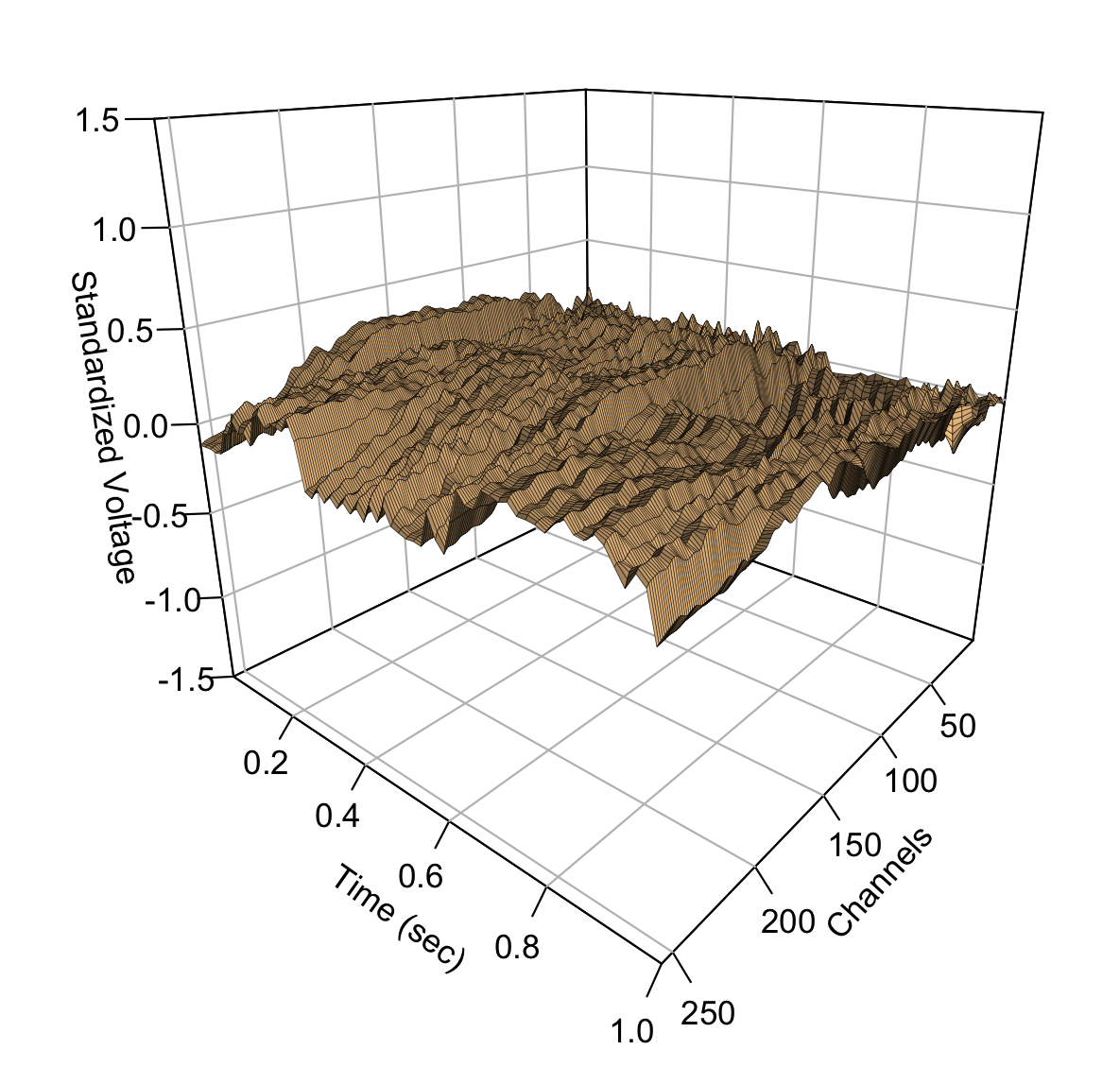}
	\caption{EEG dataset: average of matrix observations for alcoholic group (left) and control group (right) under S2.}
	\label{figure:EEG}
\end{figure}

Surprisingly, we note that the original lr-Lloyd's algorithm (Algorithm \ref{alg:Lloyd}  + Algorithm \ref{alg:initial}) would not deliver a satisfactory result on this dataset. It can be partially explained by Figure \ref{figure:EEG}, which displays the average of all trials under S2 for two groups. It is readily seen that the average matrix of control group is comparatively close to  pure noise, and hence the relaxed version of lr-Lloyd's algorithm can work reasonably well in this scenario.

\begin{table}[h]
\centering
{
\begin{tabular}{cccccc}
\hline
 & {rlr-Lloyd} & {vec-Lloyd} & {SKM}  & DTC &TBM\\ \hline
S1 & \textbf{39.34} & 42.62 & 44.26 &45.08 &43.44\\
S2 & \textbf{28.69} & 35.25 & 36.07& 39.34 & 35.25\\ \hline
\end{tabular}%
}
\caption{Clustering error of EEG dataset under S1 and S2. Note that the methods {vec-Lloyd} and {SKM} \citep{witten2010framework} refer to directly applying Lloyd's algorithm and sparse K-means on vectorized data, i.e., on rows of $\scrM_3(\bcalX^{(S_1)})$  or $\scrM_3(\bcalX^{(S_2)})$, whereas DTC\citep{sun2019dynamic}  and TBM \citep{wang2019multiway} are both tensor-based clustering methods.}
\label{table:EEG}
\end{table}


\subsubsection{Malaria parasite genes networks dataset}\label{sec:num_malaria}
We then consider the \textit{var} genes networks of the human malaria
parasite \textit{Plasmodium falciparum} constructed by
\cite{larremore2013network} via mapping $n=9$ highly variable regions
(HVRs) to a multi-layer network. 
\begin{figure}
	\centering
	\begin{subfigure}[b]{.18\linewidth}
		\includegraphics[width=\linewidth]{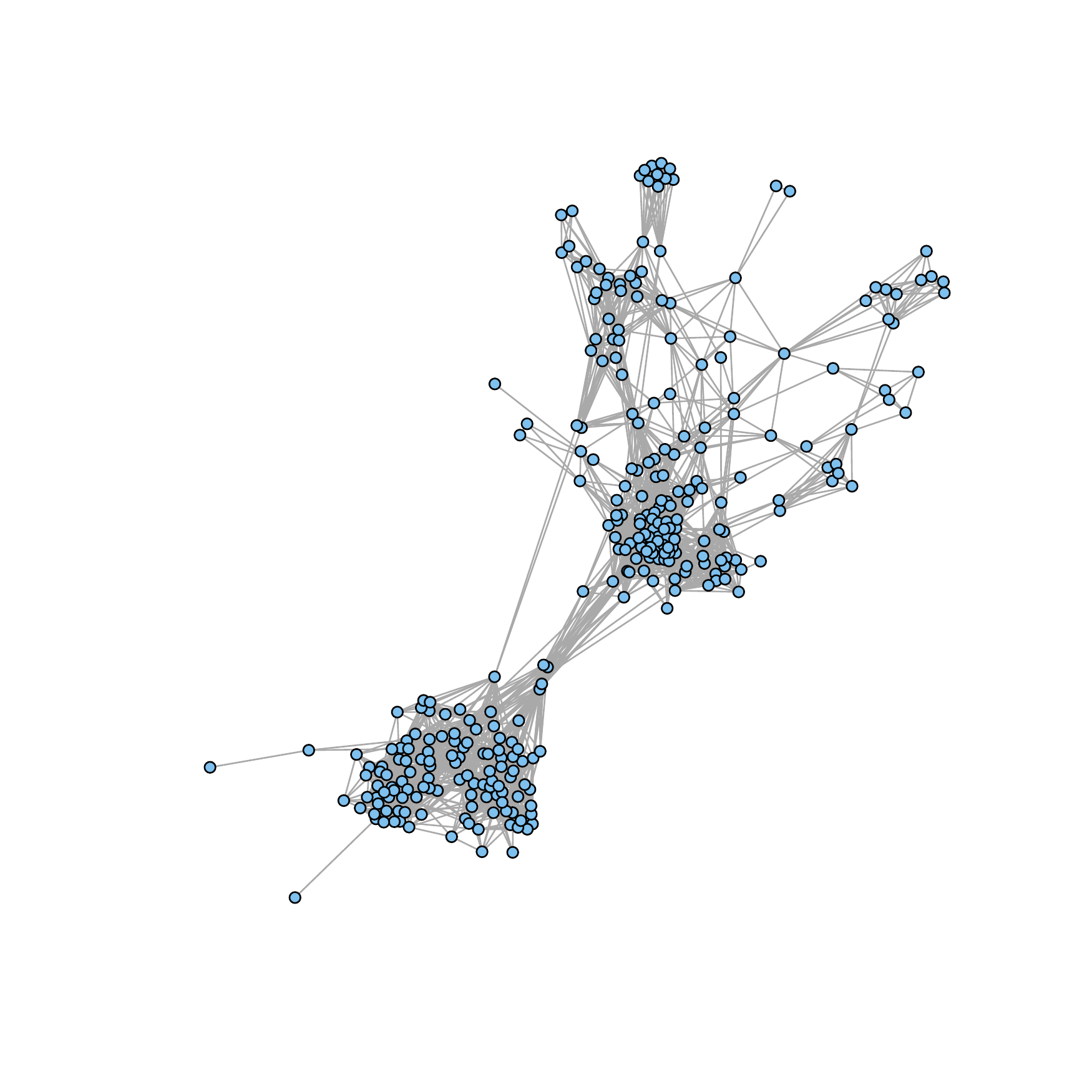}
		\caption{HVR1.}
	\end{subfigure}
	\begin{subfigure}[b]{.18\linewidth}
		\includegraphics[width=\linewidth]{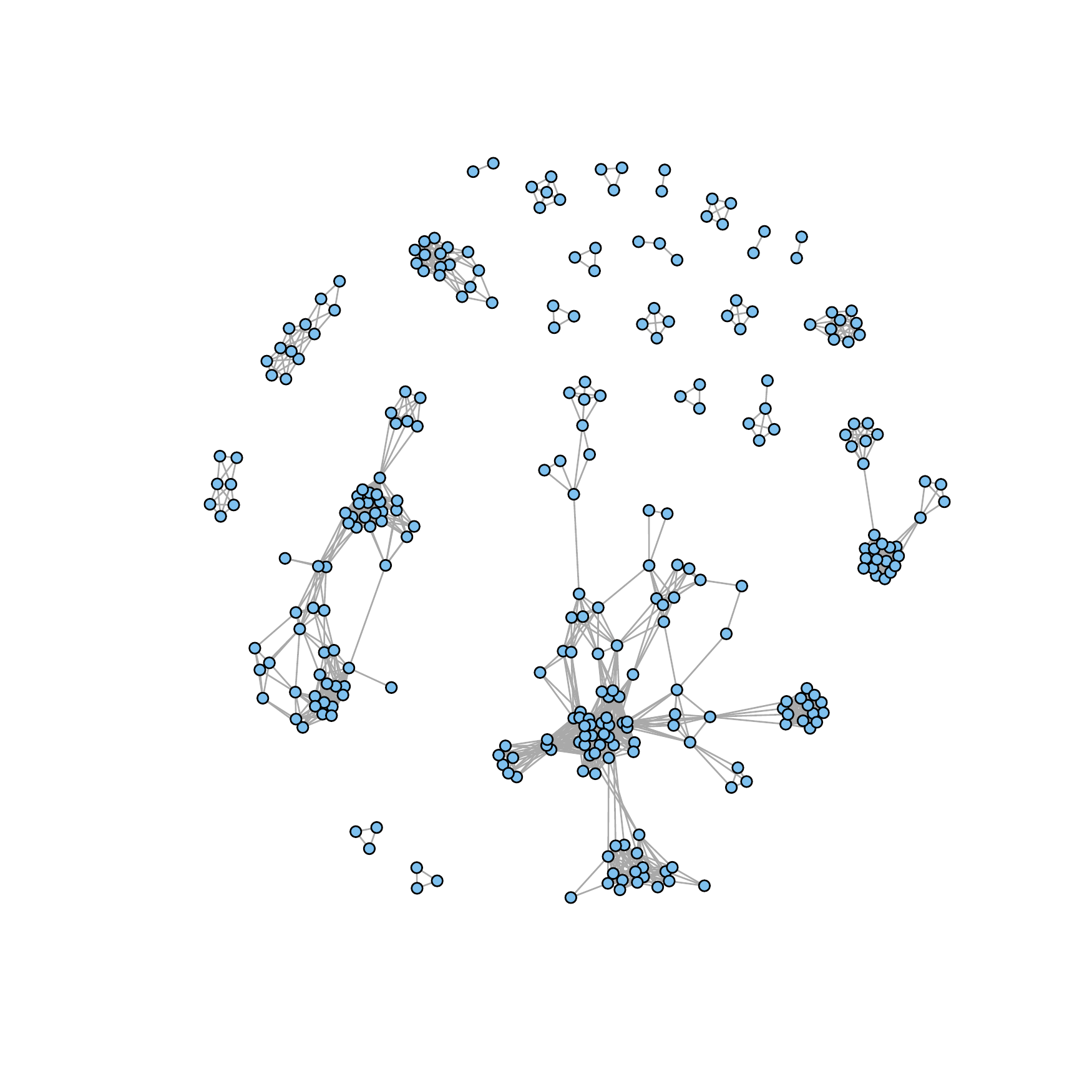}
		\caption{HVR2.}
	\end{subfigure}
	\begin{subfigure}[b]{.18\linewidth}
		\includegraphics[width=\linewidth]{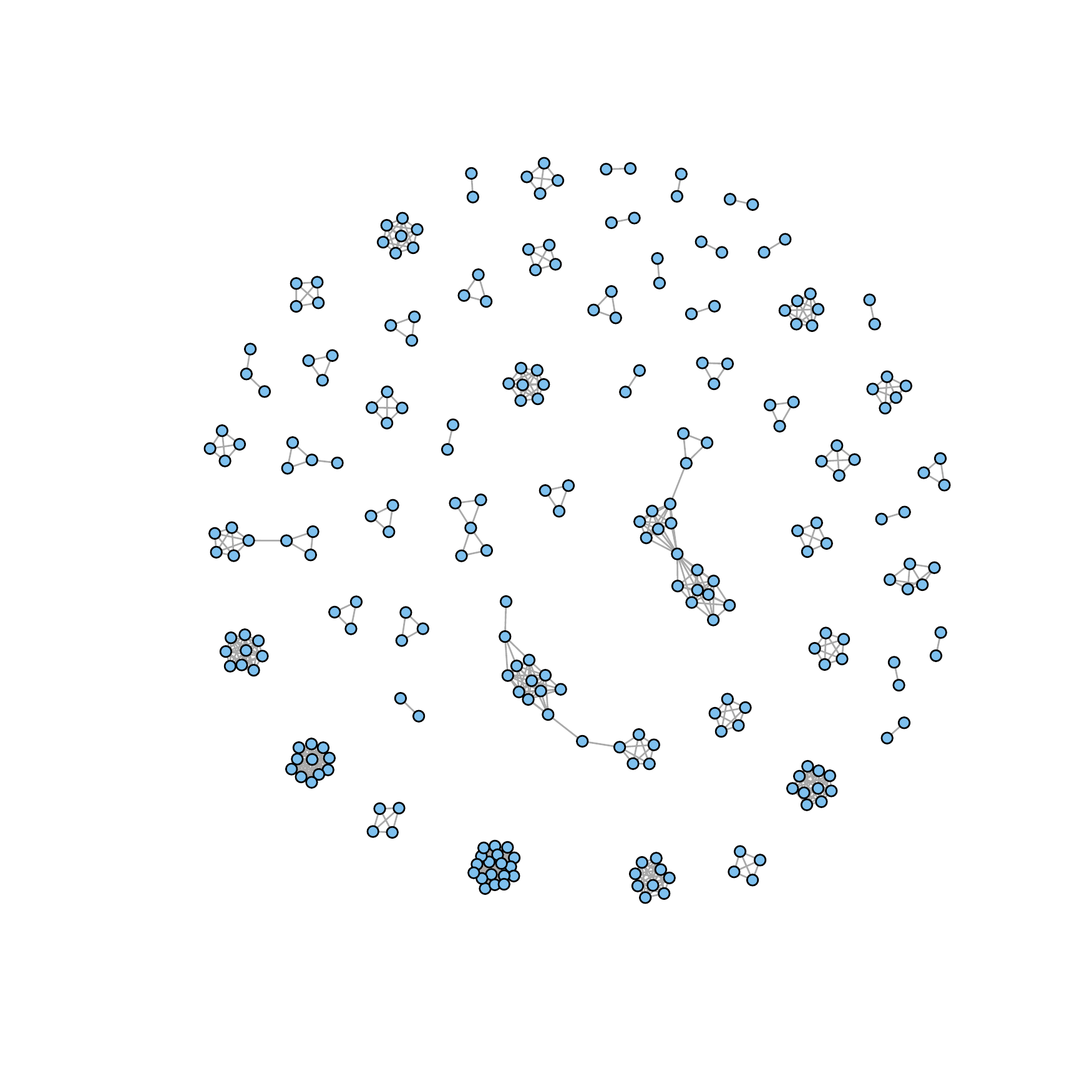}
		\caption{HVR3.}
	\end{subfigure}
	\begin{subfigure}[b]{.18\linewidth}
		\includegraphics[width=\linewidth]{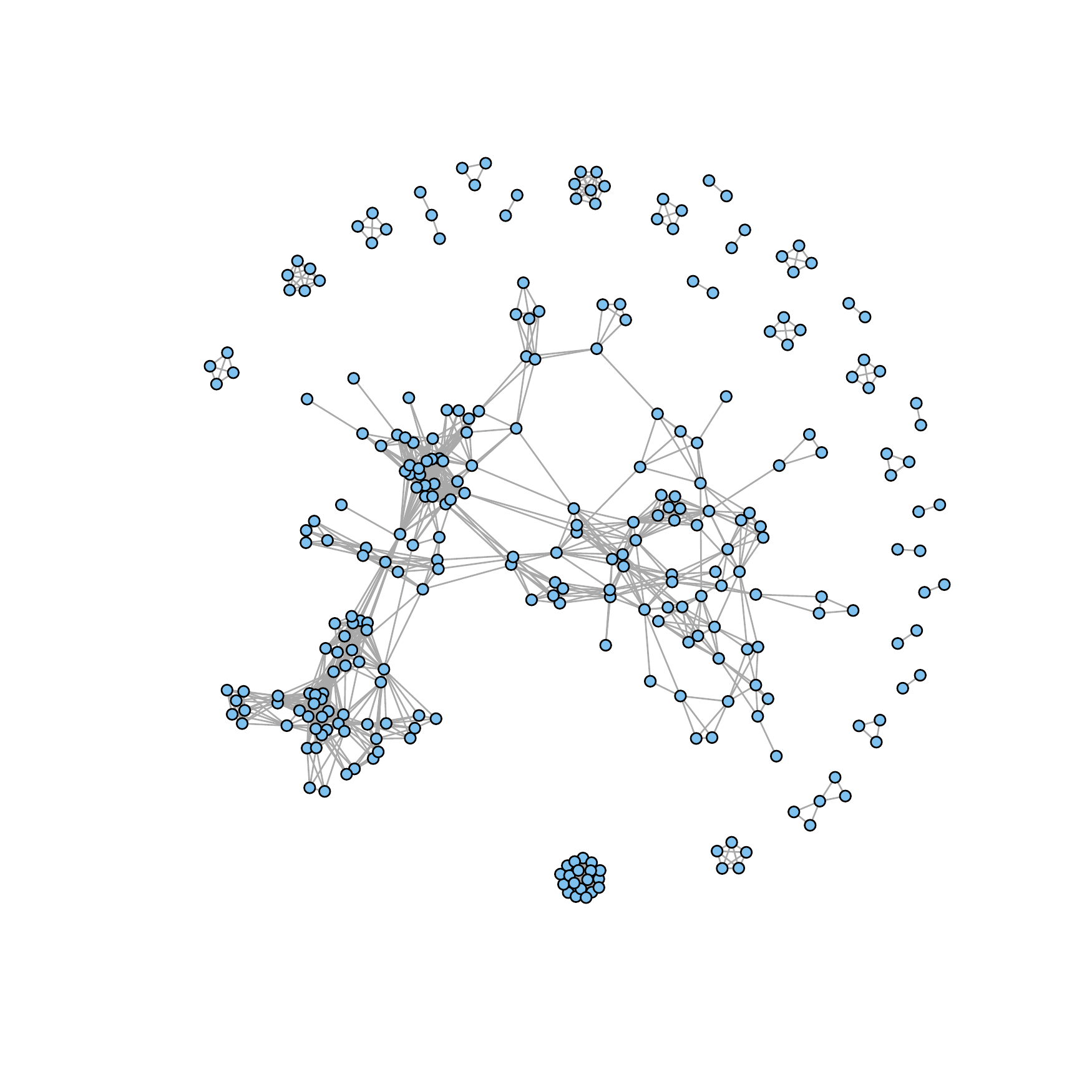}
		\caption{HVR4.}
	\end{subfigure}
	\begin{subfigure}[b]{.18\linewidth}
		\includegraphics[width=\linewidth]{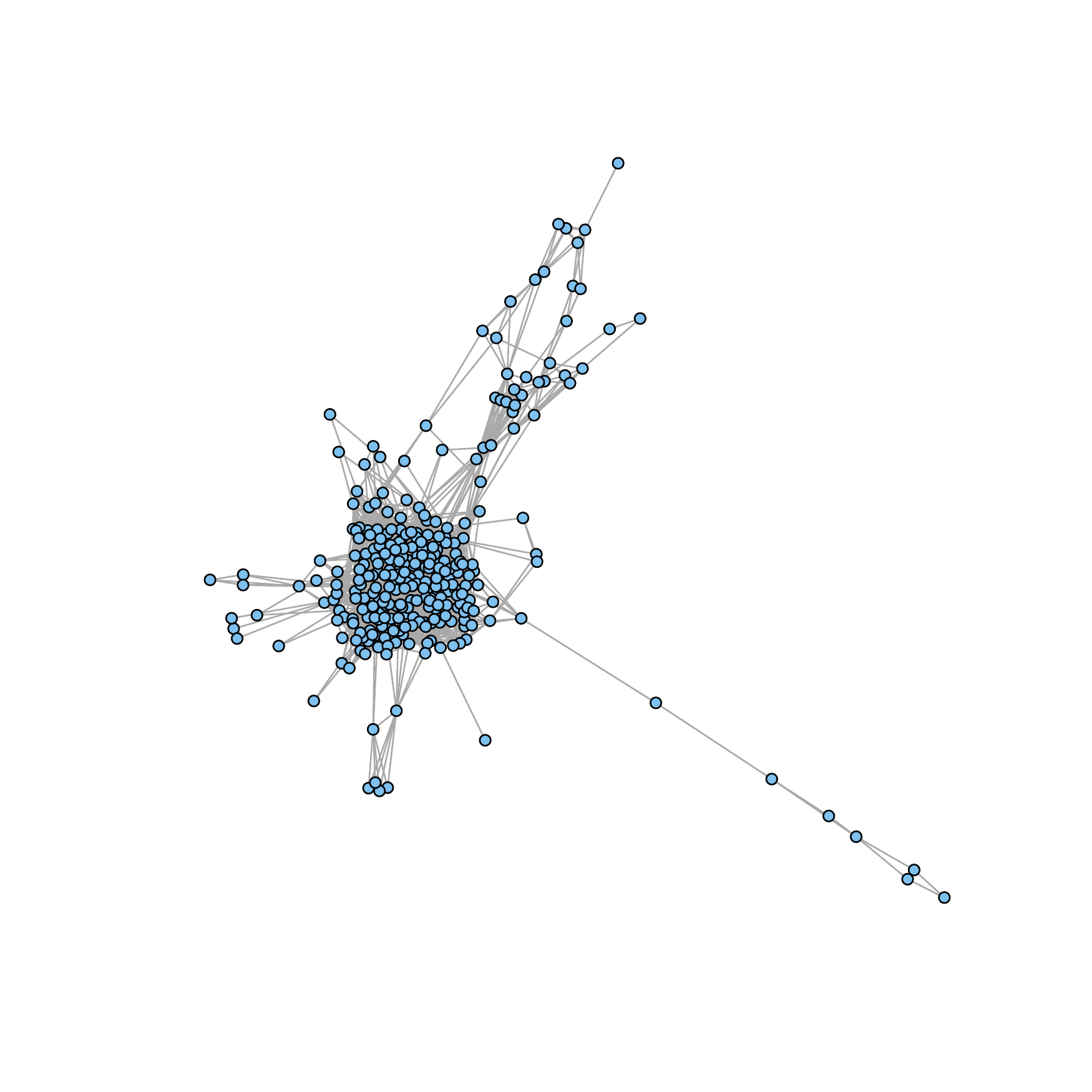}
		\caption{HVR5.}
	\end{subfigure}
	\begin{subfigure}[b]{.18\linewidth}
		\includegraphics[width=\linewidth]{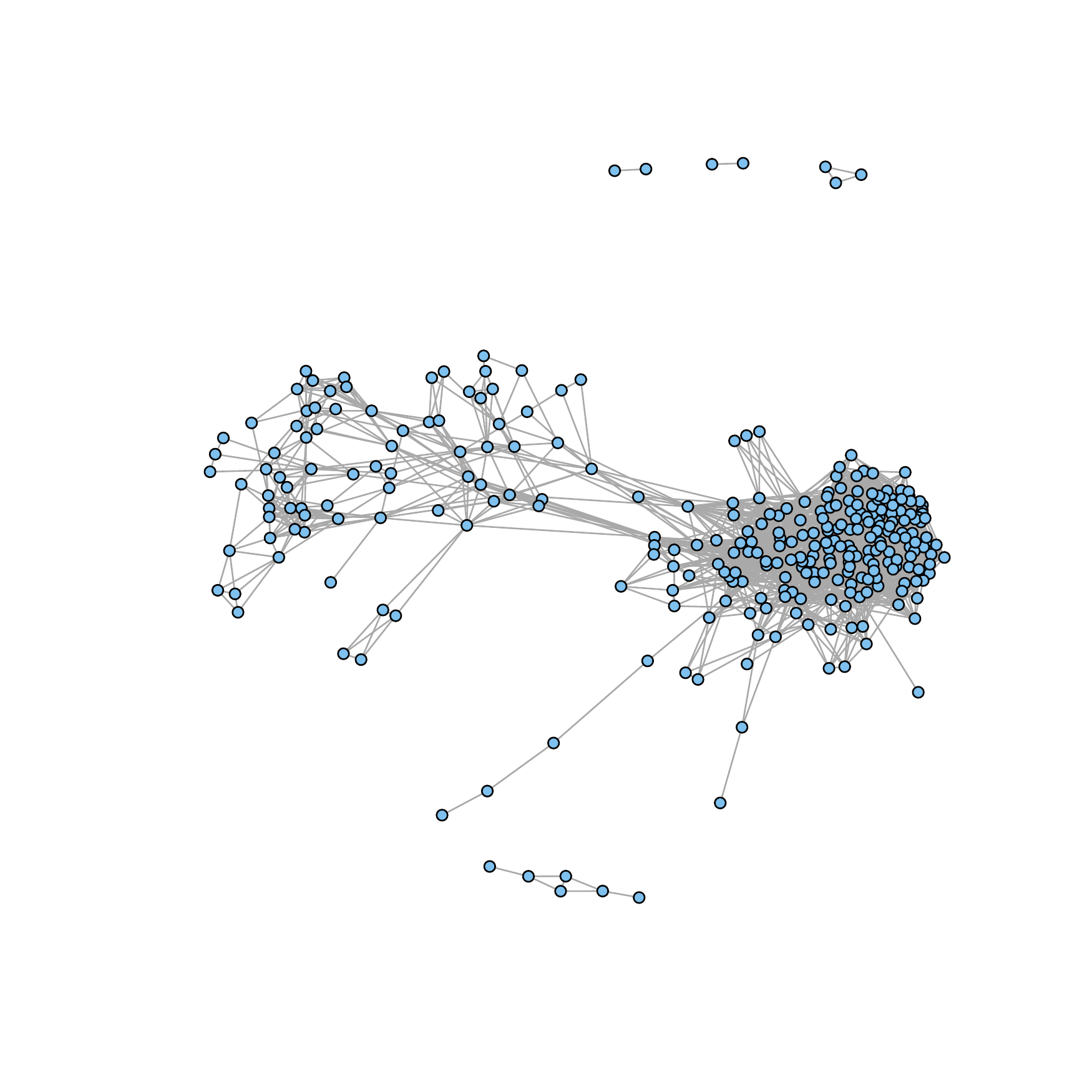}
		\caption{HVR6.}
	\end{subfigure}
	\begin{subfigure}[b]{.18\linewidth}
		\includegraphics[width=\linewidth]{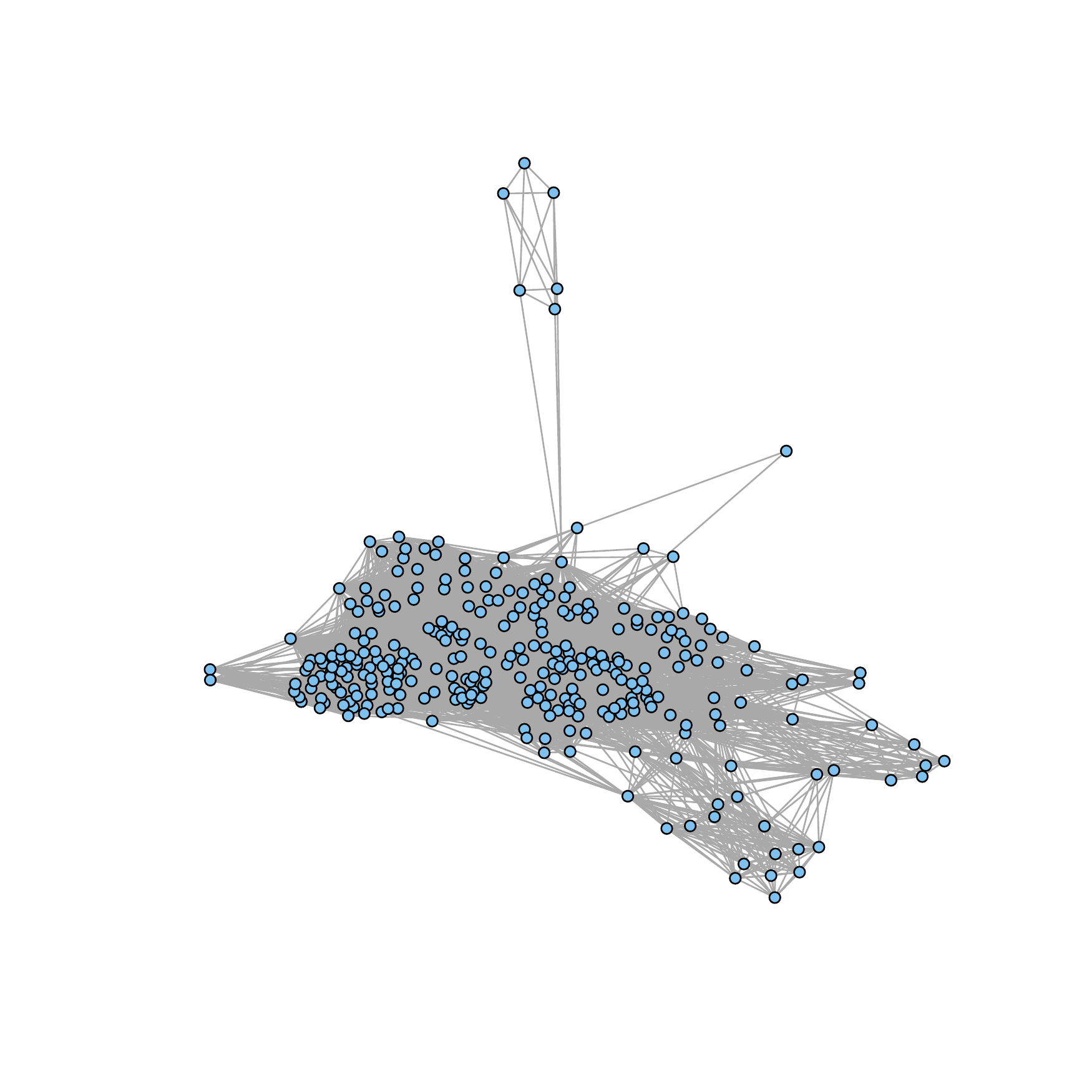}
		\caption{HVR7.}
	\end{subfigure}
	\begin{subfigure}[b]{.18\linewidth}
		\includegraphics[width=\linewidth]{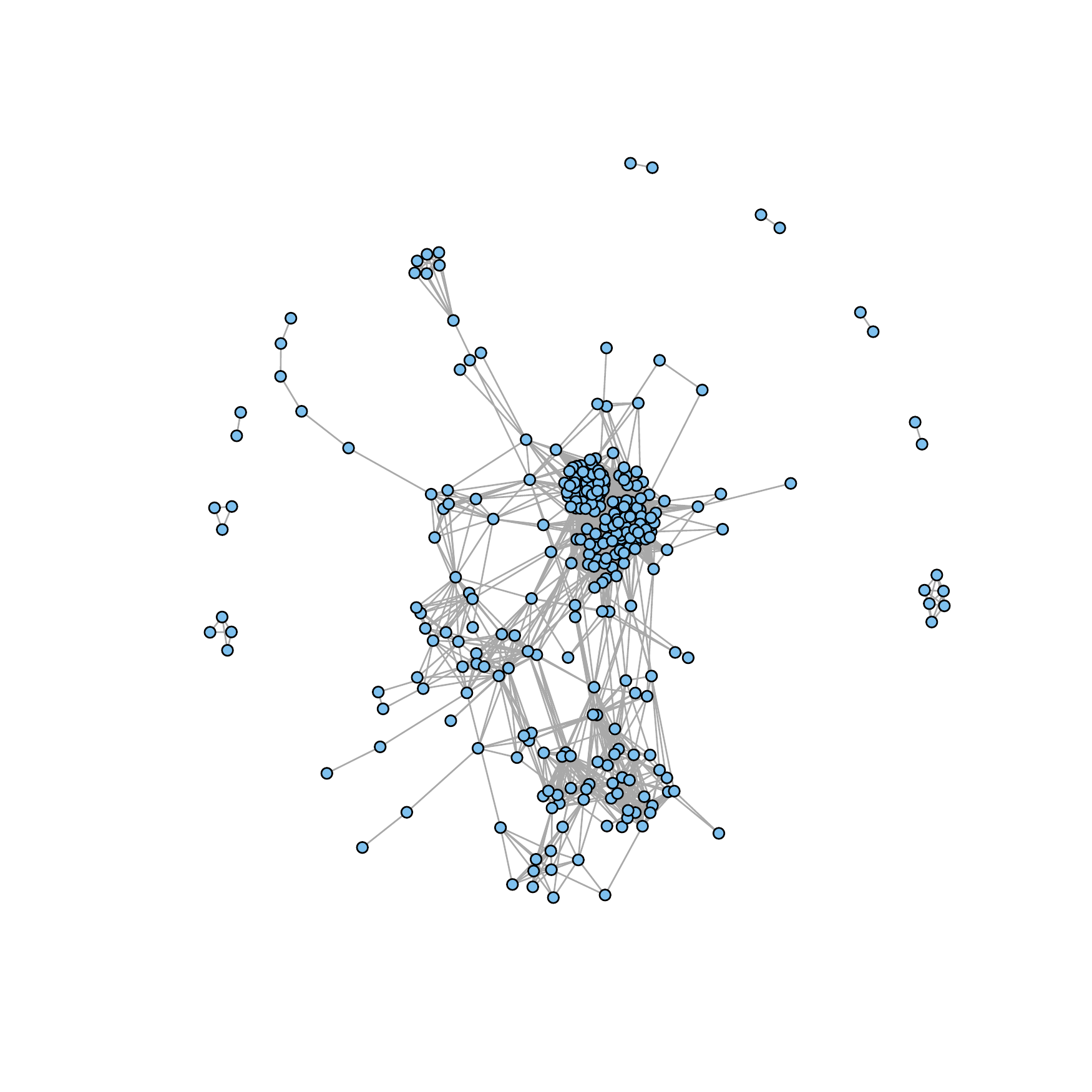}
		\caption{HVR8.}
	\end{subfigure}
	\begin{subfigure}[b]{.18\linewidth}
		\includegraphics[width=\linewidth]{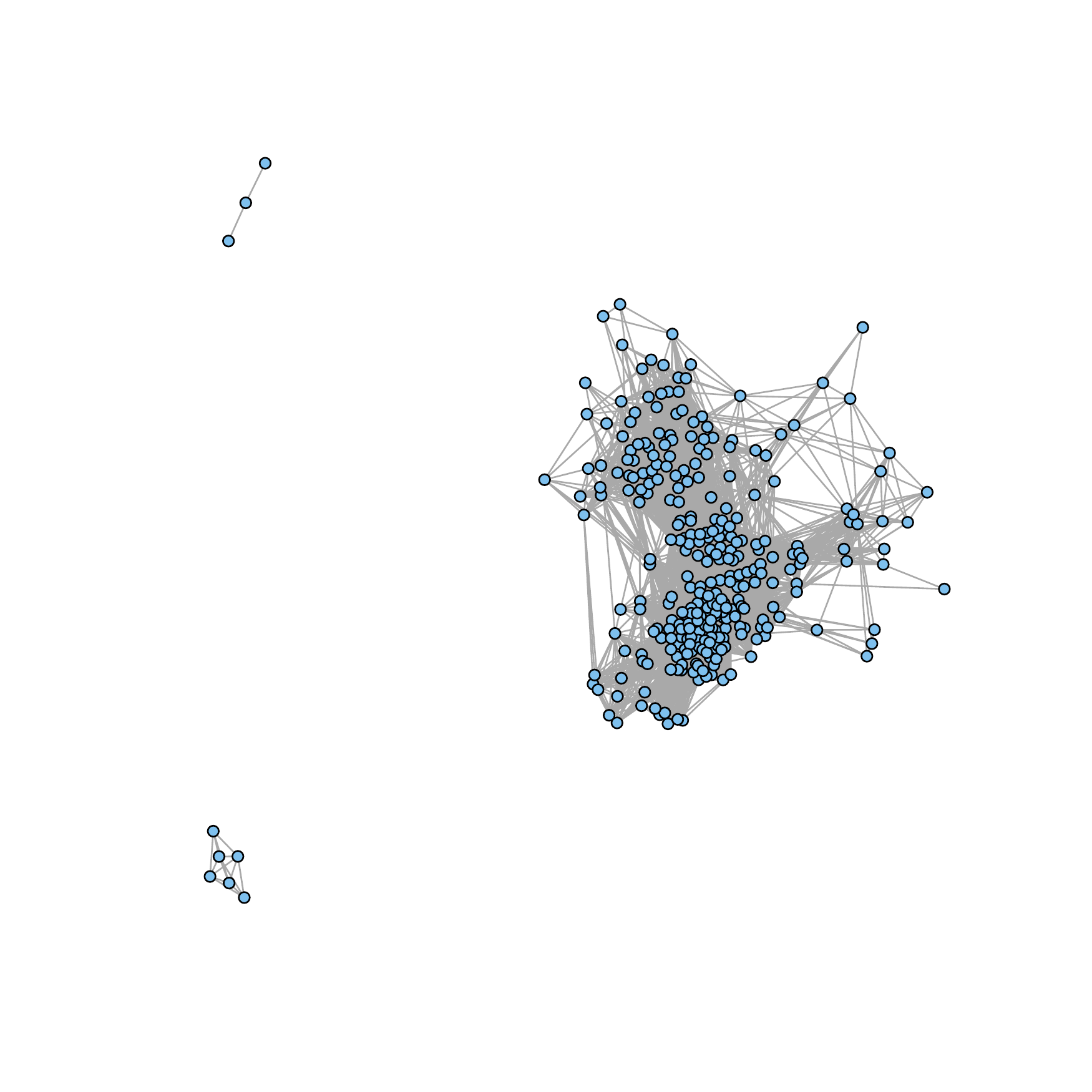}
		\caption{HVR9.}
	\end{subfigure}
	\caption{Malaria parasite genes networks dataset: 9 highly variable regions (HVRs) represented by their adjacency matrices \citep{jing2021community}}
	\label{fig:HVR}
\end{figure}
Following the practice in \cite{jing2021community}, we focus on $d_1=d_2=212$ common nodes appearing on all $9$ layers and obtain a multi-layer network adjacency tensor $\bcalX\in\{0,1\}^{212\times 212\times 9}$ with each layer being the associated adjacency matrix. Unfortunately, the method in \cite{larremore2013network} needs to discard $3$ out of $9$ HVRs due to their extreme sparse structures, referring to region $\{2,3,4\}$ in Figure \ref{fig:HVR}. This later had been remedied by the tensor-decomposition-based method TWIST in  \cite{jing2021community}. In term of clustering all layers, we expect our algorithm would have a comparable performance in contrast with the results in \cite{jing2021community}. Specifically, \cite{jing2021community} obtain a hierarchical structure with $6$ clusters of all layers by repeatedly clustering the embedding vectors. Following their practice, by setting $(r_\bU,r_\bV,K)=(15,15,6)$, we apply Algorithm \ref{alg:initial} on $\bcalX$, and find that the $9$ HVRs fall in to the following clusters:  $\{1\},\{2,3,4,5\},\{6\},\{7\},\{8\},\{9\}$. The result is exactly the same as that in \cite{jing2021community} but our method avoid repeated clustering. We remark that our tensor-based spectral initialization already produces a good initial clustering on this dataset, and thus further low-rank Lloyd's iterations seem unnecessary. In sharp contrast, it would lead to unsatisfactory result if we directly apply K-means with $K=6$ on the embedding matrix obtained by TWIST. This further demonstrates the validity and flexibility of our proposed lr-Lloyd's algorithm.

\subsubsection{UN comtrade trade flow networks dataset}
In the last example, we consider the international commodity trade flow data in $2019$ in terms of countries/regions and different types of commodities, collected by \cite{lyu2021latent} from \textit{UN comtrade Database}\footnote{The dataset is publicly available at \url{https://comtrade.un.org}.}. Following the data processing procedure in \cite{lyu2021latent}, we pick out top $d_1=d_2=48$ countries/regions ranked by  exports and obtain a weighted adjacency tensor $\wt\bcalX\in \RR^{48\times 48 \times 97}$, where $n=97$ layers represent different categories of commodities\footnote{The categories are based on 2-digit HS code in \url{https://www.foreign-trade.com/reference/hscode.htm}.}. The entry $\wt\bcalX(i_1,i_2,i_3)$ indicates the amount of exports from country $i_1$ to country $i_2$ in terms of commodity type $i_3$. To have a comparable magnitude across different entries, our data tensor is obtained after transformation $\bcalX=\log(\wt\bcalX+1)$. We emphasize that  in \cite{lyu2021latent} the edges of $\bcalX$ have to be further converted to binary under their framework, which might cause undesirable information loss. We apply Algorithm \ref{alg:Lloyd} that is initialized by Algorithm \ref{alg:initial} with parameters $(r_\bU,r_\bV,K)=(3,3,2)$ and $(r_1,r_2)=(2,2)$. These choices produce most interpretable result as summarized in Table \ref{table:UNC}. It is intriguing to notice that  cluster 1 mainly consists of products of low durability including animal \& vegetable products and part of foodstuffs, whereas cluster 2 contains most industrial products that might indicate a trend of global trading. These findings are consistent with \cite{lyu2021latent}.

\begin{table}[]
\centering
\resizebox{\textwidth}{!}{%
\begin{tabular}{lll}
\hline
Commodity cluster 1 & Commodity cluster 2\\ \hline
\textbf{01-05  Animal \& Animal Products (100\%)} & 15 Vegetable Products (13.73\%)\\
\textbf{06-14  Vegetable Products (86.27\%)}& 19-22 Foodstuffs (60.82\%)\\
\textbf{16-18, 23-24  Foodstuffs (39.18\%)} & \textbf{25,27 Mineral Products (86.68\%)}\\
26 Mineral Products (13.32\%)& \textbf{28-30,32-35,38 Chemicals \& Allied Industries (96.46\%)}\\
31,36-37 Chemicals \& Allied Industries (3.54\%)& \textbf{39-40 Plastics / Rubbers (100\%)}\\
41,43 Raw Hides, Skins, Leather, \& Furs (23.01\%)& \textbf{42 Raw Hides, Skins, Leather, \& Furs (76.99\%)}\\
45-47 Wood \& Wood Products (15.13\%)&\textbf{44,48-49 Wood \& Wood Products (84.87\%)}\\
50-55,57-58,60 Textiles (23.40\%)& \textbf{56,59,61-63 Textiles (65.97\%)}\\
65-67 Footwear / Headgear (17.45\%) & \textbf{64 Footwear / Headgear (82.55\%)}\\
75,78-81 Metals (6.44\%) & \textbf{68-71 Stone / Glass (100\%)}\\
86,89 Transportation (5.50\%) & \textbf{72-74,76,82-83 Metals (93.56\%)}\\
91-93,97 Miscellaneous (8.19\%) & \textbf{84-85 Machinery / Electrical (100\%)}\\ 
&  \textbf{87-88 Transportation (94.50\%)}\\
&\textbf{90,94-96,99 Miscellaneous (91.81\%)}\\ \hline
\end{tabular}%
}
\caption{Clustering result of UN comtrade network. The number in brackets is the percentage of the amount of exports in the corresponding type of commodity.}
\label{table:UNC}
\end{table}

\bibliographystyle{plainnat}
\bibliography{references} 

\newpage

\appendix

\section{Proofs of Main Theorems}
Throughout the proofs, we use $c,C,C^\prime$ to represent generic absolute constants, whose actual values may vary in different formulas.
\subsection{Proof of Theorem \ref{thm:main}}
\paragraph*{Step 1: Notations and Good Initialization}
We need to introduce some notations to simplify the presentation of our proof. Recall the \textit{individual signal strength} is defined as 
$$ \lambda=\min_{k\in[K]} \sigma_{\submin}(\bM_k)$$
Note in our setting, we simply have $\lambda\gtrsim  \kappa_0^{-1}r^{-1/2}\max_{a\ne b}\fro{\bM_a-\bM_b}\ge  \kappa_0^{-1}r^{-1/2}\Delta$. 
\\Define the frobenius error  with respect to the true label $\bs^\ast$:
$$
\ell(\bs,\bs^\ast):=\sum_{i=1}^n\fro{\bM_{s_i}-\bM_{s_i^\ast}}^2
$$
as well as the corresponding hamming loss:
$$
h(\bs,\bs^\ast):=\sum_{i=1}^n\ind{s_i\ne s_i^\ast}
$$
A simple relation is that $h(\bs,\bs^\ast)\le \Delta^{-2}\cdot \ell(\bs,\bs^\ast)$ due to the fact 
$$
\sum_{i=1}^n \|\bM_{s_i}-\bM_{s_i^{\ast}}\|_{\rm F}^2\geq \sum_{i=1}^n \ind{s_i\neq s_i^{\ast}}\Delta^2.
$$ 
Note that, by definition $\ell_{\textsf{c}}(\hat\bs^{(0)},\bs^\ast)=\sum_{i=1}^n\fro{\bM_{s_i^{(0)}}-\bM_{\pi(s_i^\ast)}}^2$ 
for some permutation $\pi$,  we can always relabel our $\bM_1,\cdots,\bM_K$ to $\bM_{\pi(1)},\cdots,\bM_{\pi(K)}$ after initialization. Therefore, without loss of generality we can assume $\pi=\text{Id}$ and hence $\ell(\hat\bs^{(0)},\bs^\ast)=\ell_{\textsf{c}}(\hat\bs^{(0)},\bs^\ast)$. As a result of condition \eqref{init-cond}, we also have
\begin{align}\label{ell-init-cond}
h(\hat \bs^{(0)}, \bs^{\ast})\leq \frac{\ell(\hat\bs^{(0)},\bs^\ast)}{\Delta^2} =o\left(\frac{\alpha n}{\kappa_0^2 K}\right)
\end{align}
Note that \eqref{init-cond}  can be equivalently expressed as $\ell(\hat\bs^{(0)},\bs^\ast)\le \tau$ for some $\tau=o\left(\kappa_0^{-2}{\alpha n\Delta^2}/{K}\right)$ and hence $\Delta^2\gg \kappa_0^2K\tau/(\alpha n)$.

\paragraph*{Step 2: Iterative Convergence} We then analyze the convergence property of low-rank Lloyd algorithm. Without loss of generality, given the labelling $\hat\bs^{(t-1)}$ at the $(t-1)$-th iteration, we investigate the behavior of $\hat\bs\supt$, i.e., after one iteration of Lloyd algorithm. 

To simplify the presentation, the subsequent analysis is conducted on the following events, where $C>0$ is some absolute constant. 
$$
\calQ_1=\bigcap_{a\in[K]}\left\{\op{\frac{\sum_{i=1}^n{\ind{ s_j^{\ast}=a}\bE_i}}{\sum_{j=1}^n{\ind{ s_j^{\ast}=a}}}}\le C\sqrt{\frac{d}{n^\ast_{a}}}\right\}
$$
$$
\calQ_2=\bigcap_{I\in[n]}\left\{\op{\frac{1}{\sqrt{|I|}}\sum_{i\in I}\bE_i}\le C\left(\sqrt{d}+\sqrt{n}\right)\right\}
$$
$$
\calQ_3=\bigcap_{i\in[n],a\in[K]}\left\{\left\{\op{\frac{\sum_{j\ne i}^n{\ind{ s_j^{\ast}=a}\bE_j}}{\sum_{j=1}^n{\ind{ s_j^{\ast}=a}}}}\le C\sqrt{\frac{d+\log n}{n^\ast_{a}}}\right\}\bigcap\left\{\op{\bE_i}\le C\sqrt{d+\log n}\right\}\right\}
$$

The following lemma dictates that $\calQ_1\cap \calQ_2\cap\calQ_3$ occurs with high probability. 
\begin{lemma}\label{lem:concentration}
There exists some absolute  constants $C_0,c_0>0$ such that if $d\ge C_0\log K$, then 
	$$\Prob\left(Q_1^c\cup Q_2^c\cup Q_3^c\right)\le  \exp(-c_0d)$$ 
\end{lemma}
Our goal is to establish the following relation between two successive iterations:
\begin{equation}\label{contraction}
	\ell(\hat\bs^{(t)},\bs)\le 2n\cdot \exp\left\{-\big(1-o(1)\big)\frac{\Delta^2}{8}\right\}+\frac{1}{2}\ell(\hat\bs^{(t-1)},\bs)
\end{equation}
and prove that it holds with high probability for all positive integer $t$. 

Suppose for iteration $t-1$, $\ell(\hat\bs^{(t-1)},\bs^\ast)$ satisfies \eqref{init-cond} and $h(\hat\bs^{(t-1)},\bs^\ast)$ satisfies \eqref{ell-init-cond}, which will be validated via induction in the last step.  By the definition of $\hat\s^{(t)}$, we have for each $i\in[n]$:
\begin{align*}
	\fro{\bX_i-\hat\bM\supt_{\hat s_i^{(t)}}}^2\le \fro{\bX_i-\hat\bM\supt_{s_i^\ast}}^2
\end{align*}
Rearranging terms above, we obtain
\begin{equation}\label{main-term}
	\inp{\bE_i}{\hat\bM\supt_{s_i^\ast}-\hat\bM\supt_{\hat s_i^{(t)}}}\le -\frac{1}{2}\fro{\bM_{s_i^\ast}-\bM_{\hat s_i^{(t)}}}^2+\calR\left(\hat s_i^{(t)};\hat\bs^{(t-1)}\right)
\end{equation}
where
\begin{equation*}
	\calR\left(a;\hat\bs^{(t-1)}\right):=\frac{1}{2}\left[\fro{\bM_{s_i^\ast}-\hat\bM\supt_{s_i^\ast}}^2-\fro{\bM_{s_i^\ast}-\hat\bM\supt_{a}}^2+\fro{\bM_{s_i^\ast}-\bM_{a}}^2\right]
\end{equation*}
Without loss of generality, suppose $\hat{s}_i^{(t)}=a$ for some $a\in[K]$. Set $\delta=o(1)$ that is to be determined later. The following fact is obvious. 
\begin{align*}
&~~~~\ind{\hat{s}_i^{(t)}= a}\\
	&= \ind{\hat{s}_i^{(t)}= a}\ind{\inp{\bE_i}{\hat\bM\supt_{s_i^\ast}-\hat\bM\supt_{a}}\le -\frac{1}{2}\fro{\bM_{s_i^\ast}-\bM_{a}}^2+\calR(a;\hat\bs^{(t-1)})}\\
	&\le \ind{\inp{\bE_i}{\bM_{a}-\bM_{s_i^\ast}}\ge \frac{1-\delta}{2}\fro{\bM_{s_i^\ast}-\bM_{a}}^2}\\
	&+\ind{\hat{s}_i^{(t)}= a}\ind{\inp{\bE_i}{\bM_{s_i^\ast}-\hat\bM\supt_{s_i^\ast}}+\inp{\bE_i}{\hat\bM\supt_{a}-\bM_{a}}+\calR(a;\hat\bs^{(t-1)})\ge\frac{\delta}{2}\fro{\bM_{s_i^\ast}-\bM_{a}}^2}\\
	&\le \ind{\inp{\bE_i}{\bM_{a}-\bM_{s_i^\ast}}\ge \frac{1-\delta}{2}\fro{\bM_{s_i^\ast}-\bM_{a}}^2}\\
	&+\ind{\hat s_i^{(t)}=a}\ind{\inp{\bE_i}{\bM_{s_i^\ast}-\hat\bM\supt_{s_i^\ast}}+\inp{\bE_i}{\hat\bM\supt_{a}-\bM_{a}}\ge\frac{\delta}{4}\fro{\bM_{s_i^\ast}-\bM_{a}}^2}\\
	&+\ind{\hat{s}_i^{(t)}= a}\ind{\calR(a;\hat\bs^{(t-1)})\ge\frac{\delta}{4}\fro{\bM_{s_i^\ast}-\bM_{a}}^2}
\end{align*}

By the definition of $\ell(\hat\bs^{(t)},\bs^\ast)$, we have
\begin{align*}
	\ell(\hat\bs^{(t)},\bs^\ast)&=\sum_{i=1}^n\sum_{a\in[K]\backslash\{ s_i^\ast\}}\fro{\bM_a-\bM_{s_i^\ast}}^2\ind{\hat{s}_i^{(t)}=a}\\
	&\le \sum_{i=1}^n\sum_{a\in[K]\backslash\{ s_i^\ast\}}\fro{\bM_a-\bM_{s_i^\ast}}^2\ind{\inp{\bE_i}{\bM_{a}-\bM_{s_i^\ast}}\ge \frac{1-\delta}{2}\fro{\bM_{s_i^\ast}-\bM_{a}}^2}\\
	+\sum_{i=1}^n&\sum_{a\in[K]\backslash\{ s_i^\ast\}}\fro{\bM_a-\bM_{s_i^\ast}}^2\ind{\hat s_i^{(t)}=a}\ind{\inp{\bE_i}{\bM_{s_i^\ast}-\hat\bM\supt_{s_i^\ast}}+\inp{\bE_i}{\hat\bM\supt_{a}-\bM_{a}}\ge\frac{\delta}{4}\fro{\bM_{s_i^\ast}-\bM_{a}}^2}\\
	&+\sum_{i=1}^n\sum_{a\in[K]\backslash\{ s_i^\ast\}}\fro{\bM_a-\bM_{s_i^\ast}}^2\ind{\hat{s}_i^{(t)}= a}\ind{\calR(a;\hat\bs^{(t-1)})\ge\frac{\delta}{4}\fro{\bM_{s_i^\ast}-\bM_{a}}^2}\\
	&=: \xi_{\textsf{err}}+\beta_1(\bs^\ast,\hat\bs^{(t)})+\beta_2(\bs^\ast,\hat\bs^{(t)})
\end{align*}
where we define 
\begin{align*}
	\xi_{\textsf{err}}&:=\sum_{i=1}^n\sum_{a\in[K]\backslash\{ s_i^\ast\}}\fro{\bM_a-\bM_{s_i^\ast}}^2\ind{\inp{\bE_i}{\bM_{a}-\bM_{s_i^\ast}}\ge \frac{1-\delta}{2}\fro{\bM_{s_i^\ast}-\bM_{a}}^2}
\end{align*}
and
\begin{align*}
	\beta_1(\bs^\ast,\hat\bs^{(t)}):=\sum_{i=1}^n&\sum_{a\in[K]\backslash\{ s_i^\ast\}}\fro{\bM_a-\bM_{s_i^\ast}}^2\ind{\hat s\supt_i= a}\\
	&\cdot\ind{\inp{\bE_i}{\bM_{s_i^\ast}-\hat\bM\supt_{s_i^\ast}}+\inp{\bE_i}{\hat\bM\supt_{a}-\bM_{a}}\ge\frac{\delta}{4}\fro{\bM_{s_i^\ast}-\bM_{a}}^2}
\end{align*}
and
\begin{align*}
\beta_2(\bs^{\ast}, \hat\bs\supt):=\sum_{i=1}^n\sum_{a\in[K]\backslash\{ s_i^\ast\}}\ind{\hat s\supt_i=a}\fro{\bM_a-\bM_{s_i^\ast}}^2\ind{\calR(a;\hat\bs^{(t-1)})\ge\frac{\delta}{4}\fro{\bM_{s_i^\ast}-\bM_{a}}^2}
\end{align*}
It suffices to bound $\xi_{\textsf{err}}, \beta_1(\bs^{\ast},\hat\bs^{(t)})$ and $\beta_2(\bs^{\ast},\hat\bs^{(t)})$, respectively. 

\paragraph*{Step 2.1: Bounding $\xi_{\textsf{err}}$.} Let us begin with $\EE \xi_{\textsf{err}}$. By definition, 
\begin{align*}
	\E\xi_{\textsf{err}}&=\sum_{i=1}^n\sum_{a\in[K]\backslash\{ s_i^\ast\}}\fro{\bM_a-\bM_{s_i^\ast}}^2\Prob\left({\inp{\bE_i}{\bM_{a}-\bM_{s_i^\ast}}\ge \frac{1-\delta}{2}\fro{\bM_{s_i^\ast}-\bM_{a}}^2}\right)
\end{align*}
Note that $\inp{\bE_i}{\bM_{a}-\bM_{s_i^\ast}}$ is normal distribution with mean zero and variance $\|\bM_a-\bM_{s_i^{\ast}}\|_{\rm F}^2$. The standard concentration inequality of normal random variable yields
\begin{align*}
\Prob\left({\inp{\bE_i}{\bM_{a}-\bM_{s_i^\ast}}\ge \frac{1-\delta}{2}\fro{\bM_{s_i^\ast}-\bM_{a}}^2}\right)\le\exp\left(-\frac{(1-\delta)^2}{8}\fro{\bM_{s_i^\ast}-\bM_{a}}^2\right)
\end{align*}
Therefore, 
\begin{align*}
\E\xi_{\textsf{err}}&\leq \sum_{i=1}^n\sum_{a\in[K]\backslash\{ s_i^\ast\}}\fro{\bM_a-\bM_{s_i^\ast}}^2\exp\left(-\frac{(1-\delta)^2}{8}\fro{\bM_{s_i^\ast}-\bM_{a}}^2\right).
\end{align*}
Assume $n\gg K$, $\Delta^2\gg \log K$ and let $\delta$ converge to $0$ as slow as possible, we can get 
$$
\E\xi_{\textsf{err}}\leq n\cdot \exp\left\{-\big(1-o(1)\big)\frac{\Delta^2}{8} \right\}
$$
By Markov inequality, 
$$
\PP\left(\xi_{\textsf{err}}\geq \exp(\Delta)\EE\xi_{\textsf{err}}\right)\leq \exp(-\Delta)
$$
We conclude that, with probability at least $1-\exp(-\Delta)$, 
$$
\xi_{\textsf{err}}\leq \exp(\Delta)\EE\xi_{\textsf{err}}\leq n\cdot \exp\left\{-\big(1-o(1)\big)\frac{\Delta^2}{8} \right\}
$$

\paragraph*{Step 2.2: Bounding $\beta_1(\bs^{\ast}, \hat\bs\supt)$}
By definition, 
\begin{align*}
\beta_1(\bs^{\ast},\hat\bs^{(t)})=&\sum_{i=1}^n\sum_{a\in[K]\backslash\{ s_i^\ast\}}\fro{\bM_a-\bM_{s_i^\ast}}^2\ind{\hat s\supt_i=a}
	\cdot\ind{\inp{\bE_i}{\bM_{s_i^\ast}-\hat\bM\supt_{s_i^\ast}}\ge\frac{\delta}{8}\fro{\bM_{s_i^\ast}-\bM_{a}}^2}\\
	+&\sum_{i=1}^n\sum_{a\in[K]\backslash\{ s_i^\ast\}}\fro{\bM_a-\bM_{s_i^\ast}}^2\ind{\hat s\supt_i=a}
	\cdot\ind{\inp{\bE_i}{\hat\bM\supt_{a}-\bM_a}\ge\frac{\delta}{8}\fro{\bM_{s_i^\ast}-\bM_{a}}^2}\\
	=:&\beta_{1,1}(\bs^{\ast}, \hat\bs\supt)+\beta_{1,2}(\bs^{\ast}, \hat\bs\supt)
	\end{align*}
Without loss of generality, we only prove the upper bound of the second term $\beta_{1,2}(\bs^{\ast}, \hat\bs\supt)$. Notice that the labels $\hat\bs\supt$ depend on all the noise matrices $\{\bE_i\}_{i=1}^n$, thus $\hat\bM_a\supt$ is dependent on $\bE_i$. Delicate treatment is necessary to establish a sharp upper bound for $\beta_1(\bs^{\ast}, \hat\bs\supt)$. 

Recall the definition that $\hat\bM_a\supt$ is computed by the best rank-$r_a$ approximation of $\bar{\bX}_a(\hat\bs\suptt):=(n_a^{(t-1)})^{-1}\sum_{i=1}^n \ind{\hat s_i\suptt=a}\bX_i$ with $n_a\suptt:=\sum_{i=1}^{n}\ind{\hat s_i^{(t-1)}=a}$. Denote $\hat\bU_a\supt$ and $\hat\bV_a\supt$ the left and right singular vectors of $\hat\bM_a\supt$. Then we have $\hat\bM_a\supt=\hat\bU_a\supt(\hat\bU_a\supt)^{\top}\bar{\bX}_a(\hat\bs\suptt)\hat\bV_a\supt (\hat\bV_a\supt)^{\top}$. {\it For notation simplicity, we  now  drop the superscript $(t)$ in $\hat\bU_a\supt$, $\hat\bV_a\supt$ and write $\hat\bU_a, \hat\bV_a$ instead.}

Now write 
\begin{align*}
	\hat\bM_{a}\supt-\bM_{a}&=\hat\bU_a\hat\bU_a^\top \bar{\bX}_a(\hat\bs^{(t-1)})\hat\bV_a\hat\bV_a^\top-\bM_a\\
	&=\hat\bU_a\hat\bU_a^\top\left(\frac{\sum_{i=1}^n{\ind{\hat s_i^{(t-1)}=a}(\bM_{s_i^\ast}+\bE_i)}}{\sum_{i=1}^n{\ind{\hat s_i^{(t-1)}=a}}}\right)\hat\bV_a\hat\bV_a^\top-\bM_a
\end{align*}
Recall that $n_a^\ast=\sum_{i=1}^{n}\ind{s_i^\ast=a}$. Denote 
$$
\bEsa:=(n_a^\ast)^{-1}\sum_{i=1}^{n}\ind{s_i^\ast=a}\bE_i \quad {\rm and} \quad \bEa\suptt:=(n_a\suptt)^{-1}\sum_{i=1}^{n}\ind{\hat s_i^{(t-1)}=a}\bE_i
$$
Then we can proceed as
\begin{align*}
	\hat\bM_{a}\supt-\bM_{a}&=\hat\bU_a\hat\bU_a^\top\left(\frac{1}{n_a\suptt}\sum_{i=1}^n{\ind{\hat s_i^{(t-1)}=a}\bM_{s_i^\ast}}+\bEa\suptt\right)\hat\bV_a\hat\bV_a^\top-\bM_a\\
	&=\hat\bU_a\hat\bU_a^\top\left[\bM_a+\frac{1}{n_a\suptt}\sum_{i=1}^n{\ind{\hat s_i^{(t-1)}=a}(\bM_{s_i^\ast}}-\bM_a)+\bEsa+(\bEa\suptt-\bEsa)\right]\hat\bV_a\hat\bV_a^\top-\bM_a\\
	&=\hat\bU_a\hat\bU_a^\top\left(\bM_a+\bEsa+\Delta_\bM\suptt+\Delta_\bE\suptt\right)\hat\bV_a\hat\bV_a^\top-\bM_a
\end{align*}
where we've defined
$$
\Delta_\bM\suptt:=\frac{1}{n_a\suptt}\sum_{i=1}^n{\ind{\hat s_i^{(t-1)}=a}(\bM_{s_i^\ast}}-\bM_a)\quad {\rm and}\quad \Delta_\bE\suptt:=\bEa\suptt-\bEsa
$$
For simplicity, we denote $\Delta\suptt:=\bEsa+\Delta_\bM\suptt+\Delta_\bE\suptt$ and write
\begin{equation}\label{eq:hatMa-Ma}
\hat\bM_{a}\supt-\bM_{a}=\hat\bU_a\hat\bU_a^\top\left(\bM_a+\Delta\suptt\right)\hat\bV_a\hat\bV_a^\top-\bM_a
\end{equation}  
Notice that since $h(\hat\bs^{(t-1)},\bs^\ast)$ satisfies \eqref{init-cond}, we have that
\begin{align*}
	n_a\suptt=\sum_{i=1}^{n}\ind{\hat s_i^{(t-1)}=a}&\ge \sum_{i=1}^{n}\ind{s_i^\ast=a}-\sum_{i=1}^{n}\ind{\hat s_i^{(t-1)}\ne s_i^\ast}\\
	&\geq n_a^\ast-h(\hat\bs^{(t-1)},\bs^\ast)\ge \frac{\alpha n}{K}-\frac{\alpha n}{8K}\ge \frac{7\alpha n}{8K}
\end{align*}
The following lemma is useful whose proof is postponed to Section~\ref{sec:technical_lemmas}. 
\begin{lemma}\label{lem:DeltaM_bound}
Suppose that $h(\hat\bs^{(t-1)},\bs^\ast)$ satisfies \eqref{init-cond}. Then, 
$$
\|\Delta_\bM\suptt\|\leq \frac{C_0K}{\alpha n}\min\left\{\kappa_0\lambda h_a(\hat\bs\suptt, \bs^{\ast}), \frac{\ell_a (\hat\bs\suptt, \bs^{\ast})}{\Delta}\right \}
$$
for some absolute constant $C_0>0$, where we define
$$
h_a(\hat\bs^{(t-1)},\bs^\ast):=\sum_{i=1}^n\ind{\hat s_i^{(t-1)}=a, s_i^{\ast}\ne a}+\sum_{i=1}^n\ind{\hat s_i^{(t-1)}\ne a, s_i^{\ast}= a}
$$
and 
$$
\ell_a (\hat\bs\suptt, \bs^{\ast}):=\sum_{i=1}^n\ind{\hat s_i^{(t-1)}=a, s_i^{\ast}\ne a}\fro{\bM_{\hat s_i^{(t-1)}}-\bM_{s_i^{\ast}}}^2+\sum_{i=1}^n\ind{\hat s_i^{(t-1)}\ne a, s_i^{\ast}= a}\fro{\bM_{\hat s_i^{(t-1)}}-\bM_{s_i^{\ast}}}^2
$$
Moreover, under event $\calQ_1\cap \calQ_2$, there exist absolute constants $C_1, C_2>0$ such that
$$
\op{\bEsa}\leq C_1\sqrt{\frac{dK}{\alpha n}}\quad {\rm and}\quad \|\Delta_\bE\suptt\|\leq C_2 \frac{K\sqrt{(d+n)\cdot h_a(\hat\bs^{(t-1)},\bs^\ast)}}{\alpha n}
$$
\end{lemma}
By Lemma \ref{lem:DeltaM_bound}, we obtain that 
\begin{align*}
	\op{\Delta\suptt}\le c\lambda+C\left(\alpha^{-1/2}K^{1/2}\sqrt{\frac{d}{ n}}+\alpha^{-1}K\sqrt{\frac{h_a(\hat\bs^{(t-1)},\bs^\ast)}{n}}\right)
\end{align*}
Recall that $\sigma_{\min}(\bM_a)\ge \lambda\gtrsim \kappa_{0}^{-1}r^{-1/2}\Delta $ and the condition 
{$\Delta\gg\alpha ^{-1/2}\kappa_0 K^{1/2}r^{1/2}\left((d/n)^{1/2}+1\right )$,}
we have that  $\sigma_{\min}(\bM_a)\ge  C \alpha^{-1/2}K^{1/2}\left((d/n)^{1/2}+1\right)$. Combining the  condition that $h_a(\hat\bs^{(t-1)},\bs^\ast)\le h(\hat\bs^{(t-1)},\bs^\ast)=o\left(\kappa_0^{-2}\alpha n /K\right )$ and the bound for $\Delta\suptt$, we obtain 
\begin{align}\label{eq:signal-cond-formula}
	\sigma_{\min}(\bM_a)> 3\op{\Delta\suptt}
\end{align}
Such signal strength condition is essential to obtain a  delicate representation formula for $\hat \bM_a\supt-\bM_a$ in eq. (\ref{eq:hatMa-Ma}), via the following lemma whose proof is deferred to Section~\ref{sec:technical_lemmas}. 
\begin{lemma}\label{lem:pertub-formula}
	For any rank-r matrix $\bM\in\RR^{d_1\times d_2}$ with compact SVD $\bU\bSigma\bV^\top$, where $\bU\in\OO_{d_1,r}$ and $\bV\in\OO_{d_2,r}$ and $\bSigma=\text{diag}(\sigma_1,\cdots,\sigma_r)$ with $\sigma_1\ge \cdots\ge \sigma_r>0$. Let $\Delta$ be an arbitrary $d_1\times d_2$ perturbation matrix and $\bX=\bM+\Delta$. Denote $\hat\bU\in\OO_{d_1,r},\hat\bV\in\OO_{d_2,r}$ the top-r left and right singular vectors of $\bX$. Suppose that $\sigma_r>3\|\Delta\|$, then we have the following relation:
	$$\mat{\hat\bU\hat\bU^\top -\bU\bU^\top&\bf0\\ \bf0&\hat\bV\hat\bV^\top -\bV\bV^\top}=\mat{\sum_{k\ge 1}\calS^\bU_{\bM,k}(\Delta)&\bf0\\ \bf0&\sum_{k\ge 1}\calS^\bV_{\bM,k}(\Delta)}=\sum_{k\ge 1}\calS_{\bM,k}(\Delta)$$
Here the $k$-th order perturbation term $\calS_{\bM,k}(\Delta)$ is defined as 
\begin{align}\label{eq:pertub-formula}
	\calS_{\bM,k}(\Delta):=\sum_{\bm:m_1+\cdots+m_{k+1}=k}(-1)^{1+\tau(\bm)}\cdot \mathfrak{P}^{-m_1}\Delta^*\mathfrak{P}^{-m_2}\Delta^*\cdots \Delta^*\mathfrak{P}^{-m_{k+1}}
\end{align}

where $\bm=(m_1,\cdots,m_{k+1})$ contains non-negative integers, $\tau(\bm)=\sum_{i=1}^{k+1}\II(m_i>0)$ and
$$\Delta^*:=\mat{\bf0&\Delta\\\Delta^\top&\bf0},\quad \mathfrak{P}^{-k}:=\left\{\begin{array}{cc}
\left(\begin{array}{cc}
\bf0 & \bU \bSigma^{-k} \bV^{\top} \\
\bV \bSigma^{-k} \bU^{\top} & \bf0
\end{array}\right) \text { if } k \text { is odd } \\
\left(\begin{array}{cc}
\bU \bSigma^{-k} \bU^{\top} & \bf0 \\
\bf0 & \bV \bSigma^{-k} \bV^{\top}
\end{array}\right) \text { if } k \text { is even. }
\end{array}\right.$$
for all $k\geq 1$. Specifically, $\mathfrak{P}^0=\mathfrak{P}^\perp$ denotes the orthogonal spectral projector defined by 
$$\mathfrak{P}^{\perp}=\left(\begin{array}{cc}
\bU_{\perp} \bU_{\perp}^{\top} & 0 \\
0 & \bV_{\perp} \bV_{\perp}^{\top}
\end{array}\right)$$
\end{lemma}

By Lemma~\ref{lem:pertub-formula} and \eqref{eq:signal-cond-formula}, we have the following decomposition
\begin{align*}
	\hat\bM_a\supt-\bM_a&=\hat\bU_a\hat\bU_a^\top\left(\bM_a+\Delta\suptt\right)\hat\bV_a\hat\bV_a^\top-\bM_a\\
	&=\left(\hat\bU_a\hat\bU_a^\top-\bU_a\bU_a^\top\right)\bM_a+\bM_a\left(\hat\bV_a\hat\bV_a^\top-\bV_a\bV_a^\top\right)\\
	&+\left(\hat\bU_a\hat\bU_a^\top-\bU_a\bU_a^\top\right)\bM_a\left(\hat\bV_a\hat\bV_a^\top-\bV_a\bV_a^\top\right)+\hat\bU_a\hat\bU_a^\top\Delta\suptt\hat\bV_a\hat\bV_a^\top
\end{align*}
so that we can re-write
\begin{align}
\beta_{1,2}(\bs^{\ast}, \hat\bs\supt)&\leq \sum_{i=1}^n\sum_{a\in[K]\backslash\{ s_i^\ast\}}\fro{\bM_a-\bM_{s_i^\ast}}^2
	\cdot\ind{\inp{\bE_i}{\left(\hat\bU_a\hat\bU_a^\top-\bU_a\bU_a^\top\right)\bM_a}\ge\frac{\delta}{32}\fro{\bM_{s_i^\ast}-\bM_{a}}^2}\notag\\
	+& \sum_{i=1}^n\sum_{a\in[K]\backslash\{ s_i^\ast\}}\fro{\bM_a-\bM_{s_i^\ast}}^2
	\cdot\ind{\inp{\bE_i}{\bM_a\left(\hat\bV_a\hat\bV_a^\top-\bV_a\bV_a^\top\right)}\ge\frac{\delta}{32}\fro{\bM_{s_i^\ast}-\bM_{a}}^2}\notag\\
	+\sum_{i=1}^n\sum_{a\in[K]\backslash\{ s_i^\ast\}}&\fro{\bM_a-\bM_{s_i^\ast}}^2
	\cdot\ind{\inp{\bE_i}{\left(\hat\bU_a\hat\bU_a^\top-\bU_a\bU_a^\top\right)\bM_a\left(\hat\bV_a\hat\bV_a^\top-\bV_a\bV_a^\top\right)}\ge\frac{\delta}{32}\fro{\bM_{s_i^\ast}-\bM_{a}}^2}\notag\\
	+ \sum_{i=1}^n&\sum_{a\in[K]\backslash\{ s_i^\ast\}}\fro{\bM_a-\bM_{s_i^\ast}}^2
	\cdot\ind{\inp{\bE_i}{\hat\bU_a\hat\bU_a^\top\Delta\suptt \hat\bV_a\hat\bV_a^\top}\ge\frac{\delta}{32}\fro{\bM_{s_i^\ast}-\bM_{a}}^2}\label{eq:beta_12_bound}
\end{align}
It suffices to bound each term in the RHS of above equation. 

\subparagraph*{Step 2.2.1: Treating the terms of $\big<\bE_i,  \big(\hat\bU_a\hat\bU_a^\top-\bU_a\bU_a^\top\big)\bM_a\big>$}\label{step:main-1}
By Lemma \ref{lem:pertub-formula}, we have
\begin{align}\label{eq:hatPU-U-1}
	\inp{\bE_i}{\left(\hat\bU_a\hat\bU_a^\top-\bU_a\bU_a^\top\right)\bM_a}=\sum_{k\ge 1}\inp{\bE_i}{\calS^{\bU_a}_{\bM,k}(\Delta\suptt)\bM_a}
\end{align} 
The RHS of (\ref{eq:hatPU-U-1}) is the sum of infinite series.  It turns out that delicate treatments are  necessary for general $k\ge 1$. Now we write 
\begin{align}\label{eq:hatPU-U-2}
\sum_{i=1}^n&\sum_{a\in[K]\backslash\{ s_i^\ast\}}\fro{\bM_a-\bM_{s_i^\ast}}^2
	\cdot\ind{\inp{\bE_i}{\left(\hat\bU_a\hat\bU_a^\top-\bU_a\bU_a^\top\right)\bM_a}\ge\frac{\delta}{32}\fro{\bM_{s_i^\ast}-\bM_{a}}^2}\notag\\
	\leq& \sum_{i=1}^n \sum_{a\in[K]\backslash\{ s_i^\ast\}}\fro{\bM_a-\bM_{s_i^\ast}}^2
	\sum_{k\ge 1}\ind{\inp{\bE_i}{\calS^{\bU_a}_{\bM,k}(\bEsa)\bM_a}\ge\frac{\delta}{ 2^{k+6}}\fro{\bM_{s_i^\ast}-\bM_{a}}^2}\notag\\
	+& \sum_{i=1}^n \sum_{a\in[K]\backslash\{ s_i^\ast\}}\fro{\bM_a-\bM_{s_i^\ast}}^2
	\sum_{k\ge 1}\ind{\inp{\bE_i}{\frak S_{\bU_a, k}\suptt\bM_a}\ge\frac{\delta}{2^{k+6}}\fro{\bM_{s_i^\ast}-\bM_{a}}^2}
\end{align}
where $\frak S_{\bU_a, k}\suptt:=\calS^{\bU_a}_{\bM,k}\left(\bEsa+\Delta_\bM\suptt+\Delta_\bE\suptt\right )-\calS^{\bU_a}_{\bM,k}(\bEsa)$. We start with bounding the first term on RHS of \eqref{eq:hatPU-U-2}. According to \eqref{eq:pertub-formula} in Lemma~\ref{lem:pertub-formula}, the $k$-th order perturbation term $\calS^{\bU_a}_{\bM,k}(\bEsa)\bM_a$ can be written as a sum of $2k\choose k$ series. For notational simplicity we define for any $\bB\in \RR^{d_1\times d_2}$,
\begin{align*}
	\cM({\bB}):=&\Big\{\bU_a^\top \bB\bV_{a},\bU_a^\top \bB\bV_{a\perp},\bU_{a\perp}^\top \bB\bV_{a},\bU_{a\perp}^\top \bB\bV_{a\perp},\\
&\bV_a^\top \bB^\top\bU_{a},\bV_{a}^\top \bB^\top\bU_{a\perp},\bV_{a\perp}^\top \bB^\top\bU_a,\bV_{a\perp}^\top \bB^\top\bU_{a\perp}\Big \}
\end{align*} 
By a careful inspection on \eqref{eq:pertub-formula} and the fact that $\bU^\top _{a\perp}\bM_a=0$, only terms of the form 
\begin{equation*}
\bU \bW_1\bW_2\cdots \bW_{2k-1}\bV_a^\top	
\end{equation*}
survive in the $2k\choose k$ series, where $\bU\in\{\pm \bU_{a},\pm \bU_{a\perp}\}$ and  $\bW_j\in\left\{\bSigma^{-1}\right \}\bigcup \cM({\bEsa})$ for $j\in[2k-1]$. Moreover, we have $\left|\{j:\bW_j\in\cM({\bEsa})\}\right|=k$ and $\left|\{j:\bW_j=\bSigma^{-1}\}\right|=k-1$. Without loss of generality for $i\in[n]$, we are going to bound the term
\begin{align}\label{eq:hatPU-U-inner-general}
	\inp{\bE_i}{\bU \bW_1\bW_2\cdots \bW_{2k-1}\bV_a^\top}
\end{align} 
To decouple the dependence of $\bE_i$ and $\bU \bW_1\bW_2\cdots \bW_{2k-1}\bV_a^\top$, we write 
$\bEsa=\bEsai+\bEsani$,
where $\bEsai=(n_a^\ast)^{-1}\bE_i\ind{s_i^\ast=a}$ and $\bEsani = (n_a^\ast)^{-1}\sum_{j\ne i}^n\ind{s_j^\ast=a}\bE_j$. Then for any $\bW_j\in\cM({\bEsa})$, we can decompose $\bW_j$ as 
\begin{align}\label{eq:hatPU-U-X-decomp}
\bW_j=\bW_{j,i}+\bW_{j,-i}	
\end{align}
with $\bW_{j,i}\in \calM(\bEsai)$ and  $\bW_{j,-i}\in \calM(\bEsani)$. Note that on $\calQ_3$, Lemma \ref{lem:concentration} implies that
\begin{align*}
	\op{\bEsai}\lesssim (n_a^\ast)^{-1}\sqrt{d+\log n},\quad  \op{\bEsani}\lesssim\sqrt{\frac{d+\log n}{n_a^*}}
\end{align*}
Since there are $k$ $\bW_j$'s belonging to $\cM({\bEsa})$, we can substitute \eqref{eq:hatPU-U-X-decomp} back into \eqref{eq:hatPU-U-inner-general} and obtain $2^k$ terms. These terms can be categorized into $2$ cases which will be treated separately. 
\begin{enumerate}
	\item Term of form 
	\begin{align*}
	\inp{\bE_i}{\bU \bY_{1}\bY_{2}\cdots \bY_{2k-1}\bV_a^\top}
	\end{align*} 
	where $\bY_{j}\in \left\{\bSigma^{-1}\right \}\bigcup \calM(\bEsani)$,  $\left|\{j:\bY_j=\bSigma^{-1}\}\right|=k-1$, and $\left|\{j:\bY_j\in\cM({\bEsani})\}\right|=k$. In this case, $\bE_i$ is independent of $\bU \bY_{1}\bY_{2}\cdots \bY_{2k-1}\bV_a^\top$. Then we have 
	\begin{align*}
		\fro{\bU \bY_{1}\bY_{2}\cdots \bY_{2k-1}\bV_a^\top}^2\le r_a\prod_{j=1}^{2k-1}\op{\bY_j}^2\le C^k\left(\frac{d+\log n}{n_a^*}\right)^k\frac{r_a}{\lambda^{2k-2}} 
	\end{align*}
	By general Hoeffding's inequality, we thus obtain
	\begin{align*}
	&\PP\left(\inp{\bE_i}{\bU \bY_{1}\bY_{2}\cdots \bY_{2k-1}\bV_a^\top}\ge \frac{\delta}{2^{4k+6}}\fro{\bM_{s_i^\ast}-\bM_{a}}^2\right)\\
	&\le  \EE\left(\exp\left(-\frac{c\delta^2\fro{\bM_{s_i^\ast}-\bM_{a}}^4}{2^{8k}\fro{\bU \bY_{1}\bY_{2}\cdots \bY_{2k-1}\bV_a^\top}^2}\right )\II\left(\fro{\bU \bY_{1}\bY_{2}\cdots \bY_{2k-1}\bV_a^\top}^2\le  C^k\left(\frac{d+\log n}{n_a^*}\right)^k\frac{r_a}{\lambda^{2k-2}}\right )\right)\\
	&\le \exp\left(-\frac{\delta^2\fro{\bM_{s_i^\ast}-\bM_{a}}^4\lambda^{2(k-1)}n_a^{*k}}{C^{k}r_a(d+\log n)^k}\right )\le \exp\left(-\delta^2\fro{\bM_{s_i^\ast}-\bM_{a}}^2\cdot \frac{\Delta^2}{\alpha^{-1} Kr(d+\log n)/n}\cdot \left(C^{\prime}\right )^k\right )
	\end{align*} 
	for some large constant $C^\prime>0$, where the last inequality holds due the condition $\lambda^2\gtrsim \alpha ^{-1}K(d+\log n)/n$. Therefore, we have that   
\begin{align*}
\EE &\sum_{i=1}^n \sum_{a\in[K]\backslash\{ s_i^\ast\}} \fro{\bM_a-\bM_{s_i^\ast}}^2
	\cdot\ind{\inp{\bE_i}{\bU \bY_{1}\bY_{2}\cdots \bY_{2k-1}\bV_a^\top}\ge \frac{\delta}{2^{4k+6}}\fro{\bM_{s_i^\ast}-\bM_{a}}^2} \\
	& \leq n\exp\left(-\Delta^2\cdot \frac{c\delta^2\Delta^2}{\alpha^{-1} Kr(d+\log n)/n}\cdot \left(C^{\prime}\right )^k\right)
\end{align*}
where we've set $\delta=o(1)$ in the way that it converges to $0$ sufficiently slowly compared to $\Delta^2/\left[Kr(d+\log n)(\alpha n)^{-1}\right ]$ .   By Markov inequality,  we get with probability at least $1-\exp\big(-\delta\left(C^{\prime}\right )^{k/2}\left[\Delta^2/\left[Kr(d+\log n)(\alpha n)^{-1}\right ]\right ]^{1/2}\Delta\big)$ that  
\begin{align*}
\sum_{i=1}^n &\sum_{a\in[K]\backslash\{ s_i^\ast\}} \fro{\bM_a-\bM_{s_i^\ast}}^2
	\cdot\ind{\inp{\bE_i}{\bU \bY_{1}\bY_{2}\cdots \bY_{2k-1}\bV_a^\top}\ge \frac{\delta}{2^{4k+6}}\fro{\bM_{s_i^\ast}-\bM_{a}}^2} \\
	&\leq n\exp\left(-\Delta^2\cdot \frac{c\delta^2\Delta^2}{\alpha^{-1} Kr(d+\log n)/n}\cdot \left(C^{\prime}\right )^k\right)
\end{align*}
	\item Terms of form 
	\begin{align*}
	\inp{\bE_i}{\bU \bY_{1}\bY_{2}\cdots \bY_{2k-1}\bV_a^\top}
	\end{align*} 
	where $\bY_{j}\in \left\{\bSigma^{-1}\right \}\bigcup \calM(\bEsai)\bigcup \calM(\bEsani)$, $\left|\{j:\bY_j=\bSigma^{-1}\}\right|=k-1$, $\left|\{j:\bY_j\in\cM({\bEsai})\}\right|=k_1$, $\left|\{j:\bY_j\in\cM({\bEsani})\}\right|=k_2$, $k_1+k_2=k$ and $k_1\ge 1,k_2\ge 0$. Notice that 
	\begin{align*}
	\fro{\bY_{1}\bY_{2}\cdots \bY_{2k-1}}\le \frac{r_a^{1/2}}{\lambda^{k-1}}\left(\frac{d+\log n}{n_a^*}\right)^{k_2/2}\frac{(d+\log n)^{k_1/2}}{(n_a^{\ast})^{k_1}}
	\end{align*} 
	This implies that 
	\begin{align*}
	2^{4k+6}\inp{\bE_i}{\bU \bY_{1}\bY_{2}\cdots \bY_{2k-1}\bV_a^\top}&\le 2^{4k+6}\fro{\bU^\top\bE_i\bV_a}\fro{\bY_{1}\bY_{2}\cdots \bY_{2k-1}}\\
	&\le \frac{C^{k}r_a(d+\log n)^{k/2+1/2}}{\lambda^{k-1}(n_a^{\ast})^{k/2+k_1/2}}\le C\frac{\alpha^{-1}Kr(d+\log n)}{n}
	\end{align*} 
	where the last inequality holds as $\lambda^2\gtrsim\alpha^{-1}K(d+\log n)/n$ and $k_1\ge 1$. Using the condition $\Delta^2\gg \alpha^{-1}Kr(d+\log n)/n$ and $\delta\rightarrow 0$ sufficiently slowly, we get that 
\begin{align*}
\sum_{i=1}^n &\sum_{a\in[K]\backslash\{ s_i^\ast\}} \fro{\bM_a-\bM_{s_i^\ast}}^2
	\cdot\ind{\inp{\bE_i}{\bU \bY_{1}\bY_{2}\cdots \bY_{2k-1}\bV_a^\top}\ge \frac{\delta}{2^{4k+6}}\fro{\bM_{s_i^\ast}-\bM_{a}}^2}=0
\end{align*}
\end{enumerate}
Collecting the above two facts, we conclude that in the $2^k$ terms we obtained by substituting \eqref{eq:hatPU-U-X-decomp} into \eqref{eq:hatPU-U-inner-general}, one term can be bounded exponentially (case 1) and the remaining $2^k-1$ terms vanish (case 2). Thus for \eqref{eq:hatPU-U-inner-general}, we get with probability at least $1-\exp\big(- \delta\left(C^{\prime}\right )^{k/2}\left[\Delta^2/\left[Kr(d+\log n)(\alpha n)^{-1}\right ]\right ]^{1/2}\Delta\big)$ that 
\begin{align*}
\sum_{i=1}^n &\sum_{a\in[K]\backslash\{ s_i^\ast\}} \fro{\bM_a-\bM_{s_i^\ast}}^2
	\cdot\ind{\inp{\bE_i}{\bU \bW_{1}\bW_{2}\cdots \bW_{2k-1}\bV_a^\top}\ge \frac{\delta}{2^{3k+6}}\fro{\bM_{s_i^\ast}-\bM_{a}}^2} \\
	&\leq n\exp\left(-\Delta^2\cdot \frac{c\delta^2\Delta^2}{\alpha^{-1} Kr(d+\log n)/n}\cdot \left(C^{\prime}\right )^k\right)
\end{align*}
Recall that the $k$-th order perturbation $\calS^{\bU_a}_{\bM,k}(\bEsa)\bM_a$ can be written as summation of at most $2k\choose k$ terms of form \eqref{eq:hatPU-U-inner-general}. Applying a union bound and a simple fact that ${2k\choose k}\le 4^k$, we can conclude that with probability at least $1-4^k\exp\big(- \delta\left(C^{\prime}\right )^{k/2}\left[\Delta^2/\left[Kr(d+\log n)(\alpha n)^{-1}\right ]\right ]^{1/2}\Delta\big)$,
\begin{align}\label{eq:hatPU-U-SM-k}
\sum_{i=1}^n &\sum_{a\in[K]\backslash\{ s_i^\ast\}} \fro{\bM_a-\bM_{s_i^\ast}}^2
	\cdot\ind{\inp{\bE_i}{\calS^{\bU_a}_{\bM,k}(\bEsa)\bM_a}\ge \frac{\delta}{2^{k+6}}\fro{\bM_{s_i^\ast}-\bM_{a}}^2} \nonumber\\
	&\leq n\cdot 4^k\exp\left(-\Delta^2\cdot \frac{c\delta^2\Delta^2}{\alpha^{-1} Kr(d+\log n)/n}\cdot \left(C^{\prime}\right )^k\right)
\end{align}
Now a union bound over all $k\ge 1$ gives that with probability  at least $1-\sum_{k\ge 1}4^k\exp\big(- \delta\left(C^{\prime}\right )^{k/2}\left[\Delta^2/\left[Kr(d+\log n)(\alpha n)^{-1}\right ]\right ]^{1/2}\Delta\big)$, \eqref{eq:hatPU-U-SM-k} holds for any $k\ge 1$. Notice that
\begin{align*}
	&\sum_{k\ge 1}4^k\exp\left(- \delta\left(C^{\prime}\right )^{k/2}\left(\frac{\Delta^2}{\alpha^{-1}Kr(d+\log n)/n}\right )^{1/2}\Delta\right )\\
	&\le \sum_{k\ge 1}\exp\left(- 2^k\cdot \delta \left(\frac{\Delta^2}{\alpha^{-1}Kr(d+\log n)/n}\right )^{1/2}\Delta\right )\\
	& \le \exp\left(- \delta \left(\frac{\Delta^2}{\alpha^{-1}Kr(d+\log n)/n}\right )^{1/2}\Delta\right )
\end{align*}
where the first inequality holds as $C^\prime$ is sufficiently large (e.g., $C^\prime >5$) such that $\left(C^\prime\right )^{k/2}\ge k\log 4$. Hence with probability at least $1-\exp\big(- \delta\left[\Delta^2/\left[Kr(d+\log n)(\alpha n)^{-1}\right ]\right ]^{1/2}\Delta\big)$ we have that
\begin{align*}
\sum_{i=1}^n &\sum_{a\in[K]\backslash\{ s_i^\ast\}} \fro{\bM_a-\bM_{s_i^\ast}}^2
	\sum_{k\ge 1}\ind{\inp{\bE_i}{\calS^{\bU_a}_{\bM,k}(\bEsa)\bM_a}\ge \frac{\delta}{2^{k+6}}\fro{\bM_{s_i^\ast}-\bM_{a}}^2} \nonumber\\
	&\leq n\exp\left(-\Delta^2\cdot \frac{c\delta^2\Delta^2}{\alpha^{-1} Kr(d+\log n)/n} \right)
\end{align*}
It remains to bound the second term on RHS of \eqref{eq:hatPU-U-2}. Notice that
\begin{align}
  \sum_{i=1}^n& \sum_{a\in[K]\backslash\{ s_i^\ast\}} \fro{\bM_a-\bM_{s_i^\ast}}^2
	\cdot\ind{\inp{\bE_i}{ \frak S_{\bU_a, k}\suptt\bM_a}\ge\frac{\delta}{2^{k+6}}\fro{\bM_{s_i^\ast}-\bM_{a}}^2}\notag\\
	\leq& \sum_{i=1}^n \sum_{a\in[K]}\sum_{b\in[K]\backslash\{a\}} \ind{s_i^{\ast}=b}\fro{\bM_a-\bM_b}^2
	\cdot\ind{\inp{\bE_i}{ \frak S_{\bU_a, k}\suptt\bM_a}\ge\frac{\delta}{2^{k+6}}\fro{\bM_b-\bM_{a}}^2}\notag\\
	\leq& \sum_{i=1}^n \sum_{a\in[K]}\sum_{b\in[K]\backslash\{a\}}  \ind{s_i^{\ast}=b}\fro{\bM_a-\bM_{b}}^2\cdot \frac{2^{2k+12}\inp{\bE_i}{\frak S_{\bU_a, k}\suptt\bM_a}^2}{\delta^2\fro{\bM_b-\bM_a}^4}\notag\\
	\leq&  \sum_{a\in[K]}\sum_{b\in[K]\backslash\{a\}}\op{\frak S_{\bU_a, k}\suptt\bM_a}^2\cdot \frac{2^{2k+12}\sum_{i=1}^n\ind{s_i^{\ast}=b} \inp{\bE_i} {\frak S_{\bU_a, k}\suptt\bM_a/\op{\frak S_{\bU_a, k}\suptt\bM_a}}^2}{\delta^2\fro{\bM_b-\bM_a}^2}\label{eq:beta_12_first_1_1}
\end{align}
The following lemma is needed whose proof is deferred to Section~\ref{sec:technical_lemmas}.

\begin{lemma}\label{lem:EiXibound}
There exist absolute constants $c_1,C_1>0$ such that,  for any fixed $b\in[K]$ and $d_1, d_2$ and $r$,  the following inequality holds with probability at least $1-\exp(-c_1d)$:
$$
\sup_{\substack{\bXi\in\RR^{d_1\times d_2}, {\rm rank}(\bXi)\leq r\\ \|\bXi\|\leq 1}} \sum_{i=1}^n\ind{s_i^{\ast}=b}\inp{\bE_i}{ \bXi}^2\leq C_1r(dr+n_b^\ast)
$$
\end{lemma}
We denote the event in Lemma \ref{lem:EiXibound} by $\calQ_4$ and proceed  on $\calQ_4$. By Lemma \ref{lem:EiXibound} and \eqref{eq:beta_12_first_1_1}, we obtain that
\begin{align}\label{eq:hatPU-U-SM-k-delta}
	\sum_{i=1}^n& \sum_{a\in[K]\backslash\{ s_i^\ast\}} \fro{\bM_a-\bM_{s_i^\ast}}^2
	\cdot\ind{\inp{\bE_i}{ \frak S_{\bU_a, k}\suptt\bM_a}\ge\frac{\delta}{2^{k+6}}\fro{\bM_{s_i^\ast}-\bM_{a}}^2}\notag\\
	\leq&  \sum_{a\in[K]}\sum_{b\in[K]\backslash\{a\}}\frac{C^kr(dr+n)}{\delta^2\Delta^2}\op{\frak S_{\bU_a, k}\suptt\bM_a}^2
\end{align}
It suffices for us to have an upper bound for $\op{\frak S_{\bU_a, k}\suptt\bM_a}^2$. Recall that by definition $\frak S_{\bU_a, k}\suptt\bM_a=\calS^{\bU_a}_{\bM,k}\left(\bEsa+\Delta_\bM\suptt+\Delta_\bE\suptt\right )\bM_a-\calS^{\bU_a}_{\bM,k}(\bEsa)\bM_a$, consisting of at most $(3^k-1){2k\choose k}$ terms in form of  
\begin{align*}
\bU \bW_1\bW_2\cdots \bW_{2k-1}\bV_a^\top	
\end{align*}
where $\bU\in\{\pm \bU_{a},\pm \bU_{a\perp}\}$ and  $\bW_j\in\left\{\bSigma^{-1}\right \}\bigcup \cM\left({\bEsa}\right )\bigcup \cM\left({\Delta_\bM\suptt}\right )\bigcup \cM\left({\Delta_\bE\suptt}\right )$ for $j\in[2k-1]$ with $\left|\left\{j:\bW_j=\bSigma^{-1}\right \}\right|=k-1$, $\left|\left\{j:\bW_j\in\cM\left({\bEsa}\right )\right \}\right|=k_1$, $\left|\left\{j:\bW_j\in\cM\left({\Delta_\bM\suptt}\right )\right \}\right|=k_2$, $\left|\left\{j:\bW_j\in\cM\left({\Delta_\bE\suptt}\right )\right \}\right|=k_3$ and  $k_1+k_2+k_3=k$, $k_1,k_2,k_3\ge 0$, $k_1\le k-1$. By Lemma \ref{lem:DeltaM_bound}, we have that 
\begin{align*}
	&\op{\bW_j}\le C\sqrt{\frac{dK}{\alpha n}}=:\calR_1, \quad \forall j\in \left\{l:\bW_l\in\cM\left({\bEsa}\right )\right \}\\
	&\op{\bW_j}\le \frac{CK}{\alpha n}\cdot \min\left\{\kappa_0\lambda h_a(\hat\bs\suptt, \bs^{\ast}), \frac{\ell_a (\hat\bs\suptt, \bs^{\ast})}{\Delta}\right \}=:\calR_2, \quad \forall j\in \left\{l:\bW_l\in\cM\left({\Delta_\bM\suptt}\right )\right \}\\
	&\op{\bW_j}\le \frac{CK\sqrt{(d+n)\cdot h_a(\hat\bs^{(t-1)},\bs^\ast)}}{\alpha n}=:\calR_3, \quad \forall j\in \left\{l:\bW_l\in\cM\left({\Delta_\bE\suptt}\right )\right \}
\end{align*}

 Using $k_1\le k-1$, we obtain that 
\begin{align*}
	&\op{\bU \bW_1\bW_2\cdots \bW_{2k-1}\bV_a^\top}^2\le{\lambda^{-2(k-1)}}\max_{k_1\in[k-1]}\calR_1^{2k_1}\left(\calR_2^{2(k-k_1)}+\calR_3^{2(k-k_1)}\right)\\
	&\le \lambda^{-2(k-1)}\left[\calR_2^{2k}+\calR_3^{2k}+\calR_1^{2(k-1)}\left(\calR^2_2+\calR^2_3\right)\right]\\
		&\le \frac{C^{2k}}{\lambda^{2(k-1)}}\bigg[\frac{K^{2k}}{\alpha^{2k} n^{2k}}\frac{\ell^{2k}_a(\hat\bs\suptt, \bs^{\ast})}{\Delta^{2k}} +\frac{K^{2k}[(d+n)\cdot h_a(\hat\bs^{(t-1)},\bs^\ast)]^k}{\alpha^{2k} n^{2k}}\\
	&+\left(\frac{dK}{\alpha n}\right )^{k-1}\left(\frac{K^2}{\alpha^2 n^2}\frac{\ell^2_a(\hat\bs\suptt, \bs^{\ast})}{\Delta^2}+\frac{K^2(d+n)\cdot h_a(\hat\bs^{(t-1)},\bs^\ast)}{\alpha^2 n^2}\right)\bigg]
\end{align*}
{Combining the above fact, \eqref{eq:hatPU-U-SM-k-delta} and the upper bound $2{2k\choose k}\le 2^{2k+1}$, we have that }
\begin{align}\label{eq:hatPU-U-SM-k-delta-1}
\sum_{i=1}^n& \sum_{a\in[K]\backslash\{ s_i^\ast\}} \fro{\bM_a-\bM_{s_i^\ast}}^2
	\cdot\ind{\inp{\bE_i}{ \frak S_{\bU_a, k}\suptt\bM_a}\ge\frac{\delta}{2^{k+6}}\fro{\bM_{s_i^\ast}-\bM_{a}}^2}\notag\\
	 &\le \sum_{a\in[K]}\sum_{b\in[K]\backslash\{a\}}2^{2k+1}\cdot \frac{C^kr(dr+n)}{\delta^2\Delta^2}\lambda^{-2(k-1)}\left[\calR_2^{2k}+\calR_3^{2k}+\calR_1^{2(k-1)}\left(\calR^2_2+\calR^2_3\right)\right]
	 \end{align}
The first term of \eqref{eq:hatPU-U-SM-k-delta-1} can be  bounded as 
	 \begin{align*}
&\sum_{a\in[K]}\sum_{b\in[K]\backslash\{a\}}(4C)^{2k}\cdot \frac{r(dr+n)}{\delta^2\Delta^2}\lambda^{-2(k-1)}\calR_2^{2k}\\
&\overset{(a)}{\leq} \sum_{a\in[K]}\sum_{b\in[K]\backslash\{a\}}C^{\prime 2k}\frac{r(dr+n)}{\delta^2\Delta^2}\frac{K^{2k}}{\alpha^{2k} n^{2k}}\frac{h^{2(k-1)}_a(\hat\bs\suptt, \bs^{\ast})\kappa_0^{2(k-1)}\lambda^{2(k-1)}\ell^2_a(\hat\bs\suptt, \bs^{\ast})}{\lambda^{2(k-1)}\Delta^2}\\
&\overset{(b)}{\leq} \frac{1}{4^{k+2}}\sum_{a\in[K]}\sum_{b\in[K]\backslash\{a\}}\frac{r(dr+n)}{\delta^2\Delta^2}\frac{K^{2}}{\alpha^{2} n^{2}}\frac{\ell^2_a(\hat\bs\suptt, \bs^{\ast})}{\Delta^2}\\
&\overset{(c)}{\leq} \frac{1}{4^{k+2}}\sum_{a\in[K]}\sum_{b\in[K]\backslash\{a\}}\frac{\alpha^{-1}Kr(dr/n+1)}{\delta^2\Delta^2}\ell_a(\hat\bs\suptt, \bs^{\ast})\\
&\overset{(d)}{\leq} \frac{1}{4^{k+2}}\ell(\hat\bs\suptt, \bs^{\ast})
\end{align*}
where we've used in (a) that the definition of $\calR_2$, in (b) that $h(\hat\bs\suptt, \bs^{\ast})\lesssim \kappa_0^{-1}(\alpha n/K)$, in (c) that $\ell(\hat\bs\suptt,\bs^{\ast})\le \Delta^2(\alpha n/{K})$, and in (d) that 
{$\Delta^2\gg\alpha^{-1} K^2r\left({dr}/n+1\right)$}. 
\\The second term of \eqref{eq:hatPU-U-SM-k-delta-1} can be  bounded as 
	 \begin{align*}
&\sum_{a\in[K]}\sum_{b\in[K]\backslash\{a\}}(4C)^{2k}\cdot \frac{r(dr+n)}{\delta^2\Delta^2}\lambda^{-2(k-1)}\calR_3^{2k}\\
&\leq \sum_{a\in[K]}\sum_{b\in[K]\backslash\{a\}}C^{\prime 2k}\frac{r(dr+n)}{\delta^2\Delta^2}\frac{K^{2k}}{\alpha^{2k} n^{2k}}\frac{(d^{k}+n^{k})h^{k}_a(\hat\bs\suptt, \bs^{\ast})}{\lambda^{2(k-1)}}\\
&\overset{(a)}{\leq}  \sum_{a\in[K]}\sum_{b\in[K]\backslash\{a\}}C^{\prime 2k}\frac{r(dr+n)}{\delta^2\Delta^2}\frac{\kappa_0^{2(k-1)}r^{k-1}K^{2k}}{\alpha^{2k} n^{2k}}\frac{(d^{k}+n^{k})h^{k}_a(\hat\bs\suptt, \bs^{\ast})}{\Delta^{2(k-1)}}\\
&\overset{(b)}{\leq} \sum_{a\in[K]}\sum_{b\in[K]\backslash\{a\}}C^{\prime 2k}\frac{r(dr+n)^2}{\delta^2\Delta^4}\frac{\kappa_0^{2(k-1)}K^{k}}{\alpha^{k} n^{k}}h^{k-1}_a(\hat\bs\suptt, \bs^{\ast})\ell_a(\hat\bs\suptt, \bs^{\ast})\\
&\overset{(c)}{\leq} \frac{1}{4^{k+2}}\sum_{a\in[K]}\sum_{b\in[K]\backslash\{a\}}\frac{\alpha^{-1}Kr(dr/n+1)^2}{\delta^2\Delta^4}\ell_a(\hat\bs\suptt, \bs^{\ast})\\
&\overset{(d)}{\leq} \frac{1}{4^{k+2}}\ell(\hat\bs\suptt, \bs^{\ast})
\end{align*}
where we've used in (a) that  $\lambda^2\ge \kappa_0^{-2}r^{-1}\Delta^2$, in (b) that $h_a(\hat\bs\suptt, \bs^{\ast})\Delta^2\le \ell_a(\hat\bs\suptt, \bs^{\ast})$ and 
$\Delta^2\ge C\alpha^{-1}Kr(d/n+1)$, in (c) that $h(\hat\bs\suptt, \bs^{\ast})\lesssim \kappa_0^{-2}(\alpha n/K)$,
in (d) that $\Delta^2\gg\alpha^{-1/2} Kr^{1/2}\left({dr}/n+1\right)$.
\\The last term of \eqref{eq:hatPU-U-SM-k-delta-1} can be  bounded as 
	 \begin{align*}
&\sum_{a\in[K]}\sum_{b\in[K]\backslash\{a\}}(4C)^{2k}\cdot \frac{r(dr+n)}{\delta^2\Delta^2}\lambda^{-2(k-1)}\calR_1^{2(k-1)}\left(\calR^2_2+\calR^2_3\right)\\
&\overset{(a)}{\leq}  \frac{1}{4^{k+2}}\sum_{a\in[K]}\sum_{b\in[K]\backslash\{a\}}\frac{r(dr+n)}{\delta^2\Delta^2}\left(\calR^2_2+\calR^2_3\right)\\
&\overset{(b)}{\leq} \frac{1}{4^{k+2}}\sum_{a\in[K]}\sum_{b\in[K]\backslash\{a\}}\frac{r(dr+n)}{\delta^2\Delta^2}\frac{K^2}{\alpha^2 n^2}\left(\frac{\ell^2_a(\hat\bs\suptt, \bs^{\ast})}{\Delta^2}+\frac{(d+n)\ell_a(\hat\bs\suptt, \bs^{\ast})}{\Delta^2}\right)\\
&\overset{(c)}{\leq} \frac{1}{4^{k+2}}\sum_{a\in[K]}\sum_{b\in[K]\backslash\{a\}}\left[\frac{\alpha^{-1}Kr(dr/n+1)}{\delta^2\Delta^2}\ell_a(\hat\bs\suptt, \bs^{\ast})+\frac{\alpha^{-2}K^2r(dr/n+1)^2}{\delta^2\Delta^4}\ell_a(\hat\bs\suptt, \bs^{\ast})\right]\\
&\overset{(d)}{\leq} \frac{2}{4^{k+2}}\ell(\hat\bs\suptt, \bs^{\ast})
\end{align*}
where we've used in (a) that $\lambda^2\gtrsim \alpha^{-1}Kd/n$, in (b) that $h_a(\hat\bs\suptt, \bs^{\ast})\Delta^2\le \ell_a(\hat\bs\suptt, \bs^{\ast})$, in (c) that $\ell(\hat\bs\suptt,\bs^{\ast})\le \Delta^2(\alpha n/{K})$. 
\\Collecting the above bounds and \eqref{eq:hatPU-U-SM-k-delta-1}, we conclude that the second term on RHS of \eqref{eq:hatPU-U-2} can bounded as
\begin{align*}
\sum_{i=1}^n& \sum_{a\in[K]\backslash\{ s_i^\ast\}} \fro{\bM_a-\bM_{s_i^\ast}}^2
	\sum_{k\ge 1}\ind{\inp{\bE_i}{ \frak S_{\bU_a, k}\suptt\bM_a}\ge\frac{\delta}{2^{k+6}}\fro{\bM_{s_i^\ast}-\bM_{a}}^2}\notag\\
	 &\le \sum_{k\ge 1}\frac{1}{4^{k+2}}\ell(\hat\bs\suptt, \bs^{\ast})\le\frac{1}{32}\ell(\hat\bs\suptt, \bs^{\ast}) 
\end{align*}

\subparagraph*{Step 2.2.2: Treating the terms of $\big<\bE_i,  \bM_a\big(\hat\bV_a\hat\bV_a^\top-\bV_a\bV_a^\top\big)\big>$}\label{step:main-2}
By symmetry, we can bound $\big<\bE_i,  \bM_a\big(\hat\bV_a\hat\bV_a^\top-\bV_a\bV_a^\top\big)\big>$ the same way as $\big<\bE_i,  \left(\hat\bU_a\hat\bU_a^\top-\bU_a\bU_a^\top\right)\bM_a\big>$, and  the proof is omitted.

\subparagraph*{Step 2.2.3: Treating the terms of $\big<\bE_i,  \big(\hat\bU_a\hat\bU_a^\top-\bU_a\bU_a^\top\big) \bM_a\big(\hat\bV_a\hat\bV_a^\top-\bV_a\bV_a^\top\big)\big>$}
By Lemma \ref{lem:pertub-formula}, we obtain that 
\begin{align}
		&\sum_{i=1}^n\sum_{a\in[K]\backslash\{ s_i^\ast\}}\fro{\bM_a-\bM_{s_i^\ast}}^2
	\cdot\ind{\inp{\bE_i}{\left(\hat\bU_a\hat\bU_a^\top-\bU_a\bU_a^\top\right)\bM_a\left(\hat\bV_a\hat\bV_a^\top-\bV_a\bV_a^\top\right)}\ge\frac{\delta}{32}\fro{\bM_{s_i^\ast}-\bM_{a}}^2}\notag\\
	\le & \sum_{i=1}^n\sum_{a\in[K]\backslash\{s_i^\ast\}}\fro{\bM_a-\bM_{s_i^\ast}}^2
	\sum_{k,l\ge 1}\ind{\inp{\bE_i}{\calS^{\bU_a}_{\bM,k}(\bEsa)\bM_a\calS^{\bV_a}_{\bM,l}(\bEsa)}\ge\frac{\delta}{2^{2k+7}}\fro{\bM_{s_i^\ast}-\bM_{a}}^2}\notag\\
	+ & \sum_{i=1}^n\sum_{a\in[K]\backslash\{ s_i^\ast\}}\fro{\bM_a-\bM_{s_i^\ast}}^2
	\sum_{k,l\ge 1}\ind{\inp{\bE_i}{\frak S_{\bU_a, k}\suptt\bM_a\calS^{\bV_a}_{\bM,l}(\bEsa)}\ge\frac{\delta}{2^{2k+7}}\fro{\bM_{s_i^\ast}-\bM_{a}}^2}\notag\\
	+ & \sum_{i=1}^n\sum_{a\in[K]\backslash\{ s_i^\ast\}}\fro{\bM_a-\bM_{s_i^\ast}}^2
	\sum_{k,l\ge 1}\ind{\inp{\bE_i}{\calS^{\bU_a}_{\bM,k}(\Delta\suptt)\bM_a\frak S_{\bV_a, l}\suptt}
\ge\frac{\delta}{2^{2k+7}}\fro{\bM_{s_i^\ast}-\bM_{a}}^2}\label{eq:cross-order-term}
\end{align}
where we define $\frak S_{\bV_a, k}\suptt:=\calS^{\bV_a}_{\bM,k}\left(\bEsa+\Delta_\bM\suptt+\Delta_\bE\suptt\right )-\calS^{\bU_a}_{\bM,k}(\bEsa)$ similar to $\frak S_{\bU_a, k}\suptt$.  
We start by bounding the first term on RHS of \eqref{eq:cross-order-term}. Using Lemma~\ref{lem:pertub-formula}, for any $k,l\ge 1$, $\calS^{\bU_a}_{\bM,k}(\bEsa)\bM_a\calS^{\bV_a}_{\bM,l}(\bEsa)$ can be written as a sum of at most ${2k\choose k}^2$ series, with all non-zero terms taking the form of 
\begin{equation*}
\bU \bW_1\bW_2\cdots \bW_{4k-1}\bV^\top	
\end{equation*}
 where $\bU\in\{\pm \bU_{a},\pm \bU_{a\perp}\}$ and $\bV\in\{ \bV_{a}, \bV_{a\perp}\}$,  and  $\bW_j\in\left\{\bSigma^{-1}\right \}\bigcup \cM({\bEsa})$ for $j\in[4k-1]$. Moreover, we have $\left|\{j:\bW_j\in\cM({\bEsa})\}\right|=2k$ and $\left|\{j:\bW_j=\bSigma^{-1}\}\right|=2k-1$.  Notice that by setting $\tilde k=2k$, this reduces to the case when we treat $\calS^{\bU_a}_{\bM,\tilde k}(\bEsa)\bM_a$. Following the same argument line by line (except for adjusting the constants accordingly), we can arrive at with probability at least $1-\exp\big(- \delta\left[\Delta^2/\left[Kr(d+\log n)(\alpha n)^{-1}\right ]\right ]^{1/2}\Delta\big)$,
\begin{align*}
\sum_{i=1}^n &\sum_{a\in[K]\backslash\{ s_i^\ast\}} \fro{\bM_a-\bM_{s_i^\ast}}^2
	\sum_{k\ge 1}\ind{\inp{\bE_i}{\calS^{\bU_a}_{\bM,k}(\bEsa)\bM_a\calS^{\bV_a}_{\bM,l}(\bEsa)}\ge \frac{\delta}{2^{2k+7}}\fro{\bM_{s_i^\ast}-\bM_{a}}^2} \nonumber\\
	&\leq n\exp\left(-\Delta^2\cdot \frac{c\delta^2\Delta^2}{\alpha^{-1} Kr(d+\log n)/n} \right)
\end{align*}
For the second and third terms on RHS of \eqref{eq:cross-order-term}, using Lemma \ref{lem:EiXibound} we obtain that
\begin{align*}
	&\sum_{i=1}^n\sum_{a\in[K]\backslash\{ s_i^\ast\}}\fro{\bM_a-\bM_{s_i^\ast}}^2
	\cdot \ind{\inp{\bE_i}{\frak S_{\bU_a, k}\suptt\bM_a\calS^{\bV_a}_{\bM,l}(\bEsa)}\ge\frac{\delta}{2^{2k+7}}\fro{\bM_{s_i^\ast}-\bM_{a}}^2}\notag\\
	&+\sum_{i=1}^n\sum_{a\in[K]\backslash\{ s_i^\ast\}}\fro{\bM_a-\bM_{s_i^\ast}}^2
	\cdot \ind{\inp{\bE_i}{\calS^{\bU_a}_{\bM,k}(\Delta\suptt)\bM_a\frak S_{\bV_a, l}\suptt}
\ge\frac{\delta}{2^{2k+7}}\fro{\bM_{s_i^\ast}-\bM_{a}}^2}\\
	&\leq  \sum_{a\in[K]}\sum_{b\in[K]\backslash\{a\}}\frac{C(dr+n)}{\delta^2\Delta^2}\left(\op{\frak S_{\bU_a, k}\suptt\bM_a\calS^{\bV_a}_{\bM,l}(\bEsa)}^2+\op{\calS^{\bU_a}_{\bM,k}(\Delta\suptt)\bM_a\frak S_{\bV_a, l}\suptt}^2\right)
\end{align*}
By definition, $\frak S_{\bU_a, k}\suptt\bM_a\calS^{\bV_a}_{\bM,l}(\bEsa)$ consists of at most $2\cdot 3^k{2k\choose k}$ terms and $\calS^{\bU_a}_{\bM,k}(\Delta\suptt)\bM_a\frak S_{\bV_a, l}\suptt$ consists of at most $2\cdot 3^{2k}{2k\choose k}$ terms, each being in form of  
\begin{align*}
\bU \bW_1\bW_2\cdots \bW_{4k-1}\bV^\top	
\end{align*}
where $\bU\in\{\pm \bU_{a},\pm \bU_{a\perp}\}$, $\bV\in\{ \bV_{a}, \bV_{a\perp}\}$ and  $\bW_j\in\left\{\bSigma^{-1}\right \}\bigcup \cM\left({\bEsa}\right )\bigcup \cM\left({\Delta_\bM\suptt}\right )\bigcup \cM\left({\Delta_\bE\suptt}\right )$ for $j\in[4k-1]$ with $\left|\left\{j:\bW_j=\bSigma^{-1}\right \}\right|=2k-1$, $\left|\left\{j:\bW_j\in\cM\left({\bEsa}\right )\right \}\right|=k_1$, $\left|\left\{j:\bW_j\in\cM\left({\Delta_\bM\suptt}\right )\right \}\right|=k_2$, $\left|\left\{j:\bW_j\in\cM\left({\Delta_\bE\suptt}\right )\right \}\right|=k_3$ and  $k_1+k_2+k_3=2k$, $k_1,k_2,k_3\ge 0$, $k_1\le 2k-1$. Again, this reduces to exact the case of $\frak S_{\bU_a,2k}\suptt\bM_a$. Following the same proof and adjusting constants therein, we can conclude that 
\begin{align*}
&\sum_{i=1}^n \sum_{a\in[K]\backslash\{ s_i^\ast\}} \fro{\bM_a-\bM_{s_i^\ast}}^2
	\sum_{k\ge 1}\ind{\inp{\bE_i}{ \frak S_{\bU_a, k}\suptt\bM_a\calS^{\bV_a}_{\bM,l}(\bEsa)}\ge\frac{\delta}{2^{2k+7}}\fro{\bM_{s_i^\ast}-\bM_{a}}^2}\notag\\
	& +\sum_{i=1}^n\sum_{a\in[K]\backslash\{ s_i^\ast\}} \fro{\bM_a-\bM_{s_i^\ast}}^2
	\sum_{k\ge 1}\ind{\inp{\bE_i}{\calS^{\bU_a}_{\bM,k}(\Delta\suptt)\bM_a\frak S_{\bV_a, l}\suptt}\ge\frac{\delta}{2^{2k+7}}\fro{\bM_{s_i^\ast}-\bM_{a}}^2}\notag\\
	 &\le\frac{1}{32}\ell(\hat\bs\suptt, \bs^{\ast}) 
\end{align*}

\subparagraph*{Step 2.2.4: Treating the terms of $\big<\bE_i,   \hat\bU_a\hat\bU_a^\top\Delta\suptt \hat\bV_a\hat\bV_a^\top\big>$}
The following decomposition is obvious:
\begin{align}
	&\inp{\bE_i}{\hat\bU_a\hat\bU_a^\top\Delta\suptt\hat\bV_a\hat\bV_a^\top}\notag\\
	&=\inp{\bE_i}{\bU_a\bU_a^\top(\bEsa+\Delta_\bM\suptt+\Delta_\bE\suptt)\bV_a\bV_a^\top}\notag\\
	&+\inp{\bE_i}{\left(\hat\bU_a\hat\bU_a^\top-\bU_a\bU_a^\top\right)(\bEsa+\Delta_\bM\suptt+\Delta_\bE\suptt)\bV_a\bV_a^\top}\notag\\
	&+\inp{\bE_i}{\bU_a\bU_a^\top(\bEsa+\Delta_\bM\suptt+\Delta_\bE\suptt)\left(\hat\bV_a\hat\bV_a^\top-\bV_a\bV_a^\top\right)}\notag\\
	&+\inp{\bE_i}{\left(\hat\bU_a\hat\bU_a^\top-\bU_a\bU_a^\top\right)(\bEsa+\Delta_\bM\suptt+\Delta_\bE\suptt)\left(\hat\bV_a\hat\bV_a^\top-\bV_a\bV_a^\top\right)}\label{eq:delta-order-term}
\end{align}
The first term above, i.e., $\inp{\bE_i}{\bU_a\bU_a^\top(\bEsa+\Delta_\bM\suptt+\Delta_\bE\suptt)\bV_a\bV_a^\top}$, is essentially the same as $\inp{\bE_i}{\calS_{\bM,1}^{\bU_a}(\Delta\suptt)\bM_a}$. For the second term of \eqref{eq:delta-order-term}, we further have  
\begin{align*}
	&\inp{\bE_i}{\left(\hat\bU_a\hat\bU_a^\top-\bU_a\bU_a^\top\right)(\bEsa+\Delta_\bM\suptt+\Delta_\bE\suptt)\bV_a\bV_a^\top}\\
	&=\sum_{k\ge 1}\inp{\bE_i}{\calS^{\bU_a}_{\bM,k}(\Delta\suptt)(\bEsa+\Delta_\bM\suptt+\Delta_\bE\suptt)\bV_a\bV_a^\top}
\end{align*}
Note that 
\begin{align*}
	&\inp{\bE_i}{\calS^{\bU_a}_{\bM,k}(\Delta\suptt)(\bEsa+\Delta_\bM\suptt+\Delta_\bE\suptt)\bV_a\bV_a^\top}\\
	&=\inp{\bE_i}{\calS^{\bU_a}_{\bM,k}(\bEsa)\bEsa\bV_a\bV_a^\top}+\inp{\bE_i}{\calS^{\bU_a}_{\bM,k}(\bEsa)(\Delta_\bM\suptt+\Delta_\bE\suptt)\bV_a\bV_a^\top}\\
	&+\inp{\bE_i}{\frak S_{\bU_a,k}\suptt(\bEsa+\Delta_\bM\suptt+\Delta_\bE\suptt)\bV_a\bV_a^\top}
\end{align*}
Here, $\calS^{\bU_a}_{\bM,k}(\bEsa)\bEsa\bV_a\bV_a^\top$ is of the same structure as $\calS^{\bU_a}_{\bM,k+1}(\bEsa)\bM_a$, $\calS^{\bU_a}_{\bM,k}(\bEsa)(\Delta_\bM\suptt+\Delta_\bE\suptt)\bV_a\bV_a^\top$ and $\frak S_{\bU_a,k}\suptt(\bEsa+\Delta_\bM\suptt+\Delta_\bE\suptt)\bV_a\bV_a^\top$ are 
of the same structure as $\frak S_{\bU_a, k+1}\suptt\bM_a$. By symmetry, the third term of \eqref{eq:delta-order-term} can be handled similarly. For the last term of \eqref{eq:delta-order-term}, it can be decomposed as
\begin{align*}
	&\inp{\bE_i}{\left(\hat\bU_a\hat\bU_a^\top-\bU_a\bU_a^\top\right)(\bEsa+\Delta_\bM\suptt+\Delta_\bE\suptt)\left(\hat\bV_a\hat\bV_a^\top-\bV_a\bV_a^\top\right)}\\
	&=\sum_{k,l\ge 1}\inp{\bE_i}{\calS^{\bU_a}_{\bM,k}(\Delta\suptt)(\bEsa+\Delta_\bM\suptt+\Delta_\bE\suptt)\calS^{\bV_a}_{\bM,l}(\Delta\suptt)}\\
	&=\sum_{k,l\ge 1}\inp{\bE_i}{\calS^{\bU_a}_{\bM,k}(\bEsa)\bEsa\calS^{\bV_a}_{\bM,l}(\bEsa)}+\sum_{k,l\ge 1}\inp{\bE_i}{\frak S_{\bU_a,k}\suptt\bEsa\calS^{\bV_a}_{\bM,l}(\bEsa)}\\
	&+\sum_{k,l\ge 1}\inp{\bE_i}{\calS^{\bU_a}_{\bM,k}(\Delta\suptt)\bEsa\frak S_{\bV_a,l}\suptt}+\sum_{k,l\ge 1}\inp{\bE_i}{\calS^{\bU_a}_{\bM,k}(\Delta\suptt)(\Delta_\bM\suptt+\Delta_\bE\suptt)\calS^{\bV_a}_{\bM,l}(\Delta\suptt)}
\end{align*}
Notice that $\calS^{\bU_a}_{\bM,k}(\bEsa)\bEsa\calS^{\bV_a}_{\bM,l}(\bEsa)$ is of the same structure as $\calS^{\bU_a}_{\bM,k+1}(\bEsa)\bM_a\calS^{\bV_a}_{\bM,l}(\bEsa)$, $\frak S_{\bU_a,k}\suptt\bEsa\calS^{\bV_a}_{\bM,l}(\bEsa)$ is  of the same structure as $\frak S_{\bU_a, k+1}\suptt\bM_a\calS^{\bV_a}_{\bM,l}(\bEsa)$, $\calS^{\bU_a}_{\bM,k}(\Delta\suptt)\bEsa\frak S_{\bV_a,l}\suptt$ is  of the same structure as $\calS^{\bU_a}_{\bM,k+1}(\Delta\suptt)\bM_a\frak S_{\bV_a, l}\suptt$. It suffices to note that the last term consists of at most $2\cdot 3^{2k}{2k\choose k}$ terms, each being in form of  
\begin{align*}
\bU \bW_1\bW_2\cdots \bW_{4k-1}\bV^\top	
\end{align*}
where $\bU\in\{\pm \bU_{a},\pm \bU_{a\perp}\}$, $\bV\in\{ \bV_{a}, \bV_{a\perp}\}$ and  $\bW_j\in\left\{\bSigma^{-1}\right \}\bigcup \cM\left({\bEsa}\right )\bigcup \cM\left({\Delta_\bM\suptt}\right )\bigcup \cM\left({\Delta_\bE\suptt}\right )$ for $j\in[4k-1]$ with $\left|\left\{j:\bW_j=\bSigma^{-1}\right \}\right|=2k-1$, $\left|\left\{j:\bW_j\in\cM\left({\bEsa}\right )\right \}\right|=k_1$, $\left|\left\{j:\bW_j\in\cM\left({\Delta_\bM\suptt}\right )\right \}\right|=k_2$, $\left|\left\{j:\bW_j\in\cM\left({\Delta_\bE\suptt}\right )\right \}\right|=k_3$ and  $k_1+k_2+k_3=2k$, $k_1,k_2,k_3\ge 0$, $k_1\le 2k-1$. This again reduces to the case of $\frak S_{\bU_a,2k}\suptt\bM_a$. 
\\ So far we finish the analysis of $\beta_{1,2}(\bs^{\ast}, \hat\bs\supt)$ and by symmetry the term $\beta_{1,1}(\bs^{\ast}, \hat\bs\supt)$ can be handled in a similar way. 
\paragraph*{Step 2.3: Bounding $\beta_2(\bs^{\ast}, \hat\bs\supt)$}
Recall the definition of $\calR(a;\hat\bs^{(t-1)})$, we have that
\begin{align}
\beta_2(\bs^{\ast}, \hat\bs\supt)=&\sum_{i=1}^n\sum_{a\in[K]\backslash\{ s_i^\ast\}}\ind{\hat s\supt_i=a}\fro{\bM_a-\bM_{s_i^\ast}}^2\ind{\calR(a;\hat\bs^{(t-1)})\ge\frac{\delta}{4}\fro{\bM_{s_i^\ast}-\bM_{a}}^2}\notag\\
\leq &\sum_{i=1}^n\sum_{a\in[K]\backslash\{ s_i^\ast\}}\ind{\hat s\supt_i=a}\fro{\bM_a-\bM_{s_i^\ast}}^2\ind{\frac{1}{2}\fro{\bM_{s_i^\ast}-\hat\bM_{s_i^\ast}\supt}^2\ge\frac{\delta}{12}\fro{\bM_{s_i^\ast}-\bM_{a}}^2}\notag\\
+ &\sum_{i=1}^n\sum_{a\in[K]\backslash\{ s_i^\ast\}}\ind{\hat s\supt_i=a}\fro{\bM_a-\bM_{s_i^\ast}}^2\ind{\frac{1}{2}\fro{\bM_{a}-\hat\bM_a\supt}^2\ge\frac{\delta}{12}\fro{\bM_{s_i^\ast}-\bM_{a}}^2}\notag\\
+ &\sum_{i=1}^n\sum_{a\in[K]\backslash\{ s_i^\ast\}}\ind{\hat s\supt_i=a}\fro{\bM_a-\bM_{s_i^\ast}}^2\ind{\fro{\bM_{s_i^\ast}-\bM_a}\fro{\bM_{a}-\hat\bM_a\supt}\ge\frac{\delta}{12}\fro{\bM_{s_i^\ast}-\bM_{a}}^2}\label{eq:beta_2}
\end{align}
We need to bound three terms on RHS of eq. \eqref{eq:beta_2} separately. It follows from Lemma \ref{lem:DeltaM_bound} that 
\begin{align*}
	&\fro{\bM_{s_i^\ast}-\hat\bM_{s_i^\ast}\supt}^2\le C\left(\frac{K^2}{\alpha^2n^2} \frac{\ell^2_{s_i^\ast}(\hat\bs^{(t-1)},\bs^\ast)}{\Delta^2}+\frac{K^2{(d+n)h_{s_i^\ast}(\hat\bs^{(t-1)},\bs^\ast)}}{\alpha^2n^2}+{\frac{dK}{\alpha n}}\right)
\end{align*}
Then for the first term on RHS of eq. \eqref{eq:beta_2}, we have 
\begin{align*}
	&\sum_{i=1}^n\sum_{a\in[K]\backslash\{ s_i^\ast\}}\ind{\hat{s}_i^{(t)}= a}\fro{\bM_a-\bM_{s_i^\ast}}^2\ind{\frac{1}{2}\fro{\bM_{s_i^\ast}-\hat\bM_{s_i^\ast}\supt}^2\ge\frac{\delta}{12}\fro{\bM_{s_i^\ast}-\bM_{a}}^2}\\
	&\le C^\prime \sum_{i=1}^n\ind{\hat{s}_i^{(t)}\ne s_i^\ast}\fro{\bM_{\hat s_i\supt}-\bM_{s_i^\ast}}^2\max_{a\in[K]\backslash\{ s_i^\ast\}}\frac{ \frac{K^4\ell^4_{s_i^\ast}(\hat\bs^{(t-1)},\bs^\ast)}{\alpha^4n^4\Delta^4}+\frac{K^4(d^2+n^2)h^2_{s_i^\ast}(\hat\bs^{(t-1)},\bs^\ast)}{\alpha^4n^4}+{\frac{d^2K^2}{\alpha^2n^2}}}{\delta^2\fro{\bM_{s_i^\ast}-\bM_{a}}^4}\\
	&\le \ell (\hat\bs^{(t)},\bs^\ast)\frac{C^\prime\max_{b\in[K]}{\left({\frac{K^4\ell^4_{b}(\hat\bs^{(t-1)},\bs^\ast)}{\alpha^4n^4\Delta^4}+\frac{K^4{(d^2+n^2)\ell^2_{b}(\hat\bs^{(t-1)},\bs^\ast)}}{\alpha^4n^4\Delta^4}}+{\frac{d^2K^2}{\alpha^2n^2}}\right)}}{\delta^2\Delta^4}\\
	&\le \frac{1 }{6}\ell(\hat\bs^{(t)},\bs^\ast) 
\end{align*}
where in the last inequality we've used $\Delta^2\gg \tau {K}/({\alpha n})$, $\Delta^2\gg \alpha^{-1}K\left({d}/n+1\right)$ and $\ell(\hat\bs^{(t-1)},\bs^\ast) \le \tau$.  Similarly, we can bound the second term on RHS of eq. \eqref{eq:beta_2} as
\begin{align*}
	&\sum_{i=1}^n\sum_{a\in[K]\backslash\{ s_i^\ast\}}\ind{\hat{s}_i^{(t)}= a}\fro{\bM_a-\bM_{s_i^\ast}}^2\ind{\frac{1}{2}\fro{\bM_{a}-\hat\bM_{a}\supt}^2\ge\frac{\delta}{12}\fro{\bM_{s_i^\ast}-\bM_{a}}^2}\le \frac{1}{6}\ell(\hat\bs^{(t)},\bs^\ast)
\end{align*}
It remains to consider the last term on RHS of eq. \eqref{eq:beta_2}, which has the following bound:
\begin{align*}
	&\fro{\bM_{s_i^\ast}-\bM_a}\fro{\bM_{a}-\hat\bM_a\supt}\\
	\le & C\fro{\bM_{s_i^\ast}-\bM_a}\left(\frac{K\ell_{s_i^\ast}(\hat\bs^{(t-1)},\bs^\ast)}{\alpha n\Delta }+\frac{K\sqrt{(d+n)h_{s_i^\ast}(\hat\bs^{(t-1)},\bs^\ast)}}{\alpha n}+\sqrt{\frac{dK}{\alpha n}}\right)
\end{align*}
Hence we can obtain that
\begin{align*}
	&\sum_{i=1}^n\sum_{a\in[K]\backslash\{ s_i^\ast\}}\ind{\hat{s}_i^{(t)}= a}\fro{\bM_a-\bM_{s_i^\ast}}^2\ind{\fro{\bM_{s_i^\ast}-\bM_a}\fro{\bM_{a}-\hat\bM_a\supt}\ge\frac{\delta}{12}\fro{\bM_{s_i^\ast}-\bM_{a}}^2}\\
	&\le  C^\prime \ell(\hat\bs^{(t)},\bs^\ast)\cdot \frac{\max_{b\in[K]}\left({\frac{K^2\ell^2_{b}(\hat\bs^{(t-1)},\bs^\ast)}{\alpha^2n^2\Delta^2} +\frac{K^2{(d+n)\ell_{b}(\hat\bs^{(t-1)},\bs^\ast)}}{\alpha^2n^2\Delta^2}}+{\frac{dK}{\alpha n}}\right)}{\delta^2\Delta^2}\\
	&\le \frac{1}{6}\ell(\hat\bs^{(t)},\bs^\ast)
\end{align*}
provided that $\Delta^2\gg \tau {K}/({\alpha n})$, $\Delta^2\gg \alpha^{-1}K\left({d}/{n}+1\right)$ and $\ell(\hat\bs^{(t-1)},\bs^\ast) \le \tau$. Collecting the above facts, we conclude that 
\begin{align*}
\beta_2(\bs^{\ast}, \hat\bs\supt)\le \frac{1}{2}\ell(\hat\bs^{(t)},\bs^\ast)
\end{align*}
\paragraph*{Step 3: Obtaining contraction property} Collecting all pieces in the previous steps, we arrive at with probability at least $1-\exp(-\Delta)$: 
\begin{align*}
	\ell(\hat\bs^{(t)},\bs^\ast)&\le n\exp\left(-(1-o(1))\frac{\Delta^2}{8}\right)+\frac{1}{4}\ell(\bs^\ast,\hat\bs^{(t-1)})+\frac{1}{2}\ell(\bs^\ast,\hat\bs^{(t)})
\end{align*}
as $\Delta^2/\left[Kr(d+\log n)(\alpha n)^{-1}\right ]\to\infty$.
As a consequence, we obtain the contraction property \eqref{contraction}. \\
To finish the proof for any $t\ge 1$, we use a mathematical induction step. At iteration $t=1$, the conclusion holds via above argument together with the initialization conditions  \eqref{init-cond} and \eqref{ell-init-cond}. Now suppose at iteration $t-1$ for $t\ge 2$, $\ell(\hat\bs^{(t-1)},\bs^\ast)$ satisfies \eqref{init-cond} and $h(\hat\bs^{(t-1)},\bs^\ast)$ satisfies \eqref{ell-init-cond}, via above argument  we can obtain  $\ell(\bs^\ast,\hat\bs^{(t)})\le 2n\exp\left(-(1-o(1))\frac{\Delta^2}{8}\right)+\ell(\bs^\ast,\hat\bs^{(t-1)})/2\le \tau $ as long as $\Delta^2\gg |\log(\tau/n)|$, which is automatically met by the condition for $\ell(\hat\bs^{(t-1)},\bs^\ast)$.
Moreover, we also have 
$$h(\bs^\ast,\hat\bs^{(t)})\le \Delta^{-2}\ell(\bs^\ast,\hat\bs^{(t)})\le \frac{\tau}{\Delta^{2}}=o\left(\frac{\alpha n}{\kappa_0^2K}\right) $$
This implies the conditions $\ell(\bs^\ast,\hat\bs^{(t)})\le \tau$ and $h(\bs^\ast,\hat\bs^{(t)})\le \kappa_0^{-2}\alpha n/8K$ hold for all $t\ge 0$ and hence \eqref{contraction} holds for all $t\ge 1$. Using the relation $h(\bs^\ast,\hat\bs^{(t)})\le \Delta^{-2}\ell(\bs^\ast,\hat\bs^{(t)})$ and the condition $\Delta^2\gg \kappa_0^2K\tau/(\alpha n)$, with probability greater than $1-\exp(-\Delta)$, for each $t\ge 0$ we have that
\begin{equation*}
	n^{-1}\cdot h(\hat\bs^{(t)},\bs)\le \exp\left(-(1-o(1))\frac{\Delta^2}{8}\right)+2^{-t}
\end{equation*}
The proof is completed by applying a union bound accounting for the events $\calQ_1,\calQ_2,\calQ_3,\calQ_4$.

\subsection{Proof of Theorem \ref{thm:spec-initialization}}
We first characterize the error of $\hat\bU$ and $\hat\bV$ and without loss of generality, we only consider $\hat\bU$. Following the same argument in the proof of Theorem 1 in \cite{zhang2018tensor}, one can obtain that there exists some absolute constant $c_0,C_0>0$ such that if $\sigma_{\min}(\scrM_1(\bcalM))\ge C_0(dr_\bU)^{1/2}n^{1/4}$, then with probability at least $1-\exp(-c_0(n\wedge d))$:
\begin{equation*}
	\fro{\sin\Theta(\hat \bU,\bU^\ast)}\le \frac{C(dr_\bU)^{1/2}\left[\sigma_{\min}(\scrM_1(\bcalM))+(dn)^{1/2}\right]}{\sigma^2_{\min}(\scrM_1(\bcalM))}\le \frac{1}{4\sqrt{2}}
\end{equation*}

Combined with the bound for $\hat\bV$, we conclude that if  $\max\{\sigma_{\min}(\scrM_1(\bcalM)),\sigma_{\min}(\scrM_2(\bcalM))\}\ge C_0(dr_\bU)^{1/2}n^{1/4}$, then with probability at least $1-\exp(-c_0(n\wedge d))$:
\begin{equation}\label{eq:init-UV-small}
	\max\left\{\fro{\sin\Theta(\hat \bU,\bU^\ast)},\fro{\sin\Theta(\hat \bV,\bV^\ast)}\right\}\le  \frac{1}{4\sqrt{2}}
\end{equation}
Denote the above event by $\calQ_{0,1}$ and we proceed  on $\calQ_{0,1}$.
\\We then analyze the performance of spectral clustering based on $\hat\bcalG=\bcalX\times_1\hat \bU\hat \bU^\top \times_2\hat\bV\hat \bV^\top $. Our proof is based on the proof for Lemma 4.2 in \cite{loffler2021optimality} with slight modification. Let $\bcalG:=\bcalM\times_1\hat \bU\hat \bU^\top \times_2\hat\bV\hat \bV^\top $ denote the signal part of $\hat\bcalG$ (also $\bG:=\scrM_3(\bcalG)$) and $\bfrakM=[vec(\hat \bM_{\hat s_1^{(0)}}),\cdots,vec(\hat \bM_{\hat s_n^{(0)}})]^\top \in\RR^{n\times d_1d_2}$ denote the corresponding k-means solution.  We claim the following lemma, whose proof is deferred to Section \ref{sec:technical_lemmas}.
\begin{lemma}\label{lem:init-two-facts} Suppose $\calQ_{0,1}$ holds. Then we have the following facts:
	\begin{enumerate}
	\item[(I)] $\bfrakM$, the k-means solution,  is close to $\bG$, i.e.,  there exists some absolute constants $c_0,C_0>0$ such that with probability at least $1-\exp(-c_0d)$:
	$$\fro{\bfrakM-\bG}\le C_0\sqrt{K}\left(\sqrt{dKr+n}\right)$$
	\item[(II)] The rows of $\bG$ belonging to different clusters is well-separated, i.e.
	\begin{align*}
\fro{\bcalG\times_3 (\be_i^\top-\be_j^\top)}\ge \frac{\Delta}{2}
\end{align*}
for any $i,j\in[n],s_i^\ast\ne s_j^\ast$.
\end{enumerate}
\end{lemma}
We proceed  on the event $\calQ_{0,2}:=\{\text{(I) holds}\}$.
Define the following set 
$$S=\left\{i\in[n]:\op{[\bfrakM]_{i\cdot}-[\bG]_{i\cdot}}\ge 
\frac{\Delta}{4}\right\}$$
Then by construction we have
$$|S|\le \frac{\fro{\bfrakM-\bG}^2}{(\Delta/4)^2}\le \frac{\alpha n}{2K}$$
where the last inequality is due to the condition $\Delta^2\ge {32C_0^2\alpha^{-1} K^2}\left({dKr}/{n}+1\right)$.\\
We claim that all indices in $S^c$ are correctly clustered. To see this, let
$$N_k=\{i\in[n]:s_i^\ast=k,i\in S^c\}$$
The following two facts hold:
\begin{itemize}
	\item For each $k\in[K]$, $|N_k|\ge n_k^\ast-|S|\ge \alpha n/(2K)>0$
	\item For each pair $a,b\in[K], a\ne b$, there cannot exist some $i\in N_a$ and $j\in N_b$ such that $\hat s_i^{(0)}=\hat s_j^{(0)}$. Otherwise we have $\hat \bM_{\hat s_i^{(0)}}=\hat \bM_{\hat s_j^{(0)}}$ and it follows that 
\begin{align*}
\op{[\bG]_{i\cdot}-[\bG]_{j\cdot}}&\le \op{[\bG]_{i\cdot}-[\bfrakM]_{i\cdot}}+\op{[\bfrakM]_{i\cdot}-[\bfrakM]_{j\cdot}} +\op{[\bfrakM]_{j\cdot}-[\bG]_{j\cdot}} \\
	&< \frac{\Delta}{2}
\end{align*}
which contradicts (II).
\end{itemize}
The above two facts imply that sets $\{\hat s_i^{(0)}:i\in N_k\}$ are disjoint for all $k\in[K]$. Therefore, there exists a permutation $\pi$ such that $\sum_{i\in S^c}\II\left(\hat s_i^{(0)}\ne \pi(s_i^\ast)\right)=0$, i.e., indices in $S^c$ are correctly clustered. Therefore, we have that 

$$n^{-1}\cdot h_{\textsf{c}}(\hat\bs^{(0)},\bs^\ast)\le n^{-1}\cdot |S|\le \frac{C{K}}{\Delta^2}\left(\frac{dKr}{n}+1\right)$$
Moreover, we have
\begin{align*}
	n^{-1}\cdot \ell_{\textsf{c}}(\hat\bs^{(0)},\bs^\ast)&\le  \frac{1}{n}\sum_{i=1}^n\fro{\bM_{\hat s_i^{(0)}}-\bM_{\pi(s_i^\ast)}}^2\II\left(\hat s_i^{(0)}\ne \pi(s_i^\ast)\right)\\
	&\le  \frac{1}{n}|S|\gamma^2\Delta^2\le C{\gamma^2K}\left(\frac{dKr}{n}+1\right)
\end{align*}
The proof is completed by taking union bound over $\calQ_{0}^c:=\calQ_{0,1}^c\bigcup \calQ_{0,2}^c$.

\subsection{Proof of Theorem \ref{thm:minimax-lower-bound}}
We essentially follow a similar argument of \cite{gao2018community}. Without loss of generality we assume $\fro{\bM_1-\bM_2}=\Delta$. Consider the $\bs^\ast\in[K]^n$ such that $n_1^\ast\le n_2^\ast\le \cdots  \le n_K^\ast$ and $n_1^\ast=n_2^\ast=\lfloor \alpha n/K\rfloor$. For every $k\in[K]$, we can choose a subset $\mathfrak{
N}_k\subset\{i\in[n]:s_i^*=k\}$ with cardinality $\lceil n_k^\ast-\frac{\alpha n}{4K^2}\rceil$. And let $\mathfrak{N}=\bigcup_{k=1}^K \mathfrak{N}_k$ denote the collection of samples in $\mathfrak{N}_k$'s. Define the following parameter space for $\bs$: 
$$\bS^\ast=\{\bs\in[K]^n:s_i=s_i^\ast\text{~for~} i \in \mathfrak{N}\}$$
For any two $\bs,\bs^\prime\in \bS^\ast$ such that $\bs\ne \bs^\prime$, we have
$$\frac{1}{n}\sum_{i=1}^n\II(s_i\ne s_i^\prime)\le \frac{K}{n}\frac{\alpha n}{4K^2}=\frac{\alpha}{4K} $$
Meanwhile, for any permutation $\pi\ne\text{Id}$ from $[K]$ to $[K]$, we have
$$\frac{1}{n}\sum_{i=1}^n\II(\pi(s_i)\ne s_i^\prime)\ge \frac{K}{n}
\left(\frac{\alpha n}{K}-\frac{\alpha n}{4K^2}\right)\ge \frac{3\alpha}{4K}$$
Therefore, we conclude that $h_{\textsf{c}}(\bs,\bs^\prime)=h(\bs,\bs^\prime)=\sum_{i=1}^n\II(s_i\ne s_i^\prime)$ for any $\bs,\bs^\prime\in \bS^\ast$. Define the parameter space
\begin{align*}
\Omega(d_1,d_2,n,K,\alpha)=\Big\{(\{\bM_k\}_{k=1}^K,\bs):~&\bM_k\in\mathbb{R}^{d_1\times d_2},\text{rank}(\bM_k)=r_k, \forall k\in[K],\bs\in[K]^n,\\
&\min_{k\in[K]} |\{i\in[n]:s_i=k\}|\ge \alpha n/K,\min_{a\ne b}\fro{\bM_a-\bM_b}\ge \Delta\Big\}	
\end{align*}
and 
\begin{align*}
\Omega_0(d_1,d_2,n,K,\alpha)=\Big\{(\{\bM_k\}_{k=1}^K,\bs):~&\bM_k\in\mathbb{R}^{d_1\times d_2},\text{rank}(\bM_k)=r_k, \forall k\in[K],\bs\in\bS^\ast,\\
&\min_{k\in[K]} |\{i\in[n]:s_i=k\}|\ge \alpha n/K,\min_{a\ne b}\fro{\bM_a-\bM_b}\ge \Delta\Big\}
\end{align*}
Since $\Omega_0\subset\Omega$, we have
\begin{equation}\label{minimax_ineq}
\inf_{\hat \bs}\sup_{\Omega}\E h_{\textsf{c}}(\hat \bs, \bs)\ge\inf_{\hat \bs}\sup_{\Omega_0}\E h_{\textsf{c}}(\hat \bs, \bs)\ge \inf_{\hat \bs}\frac{1}{|\bS^\ast|}\sum_{\bs\in \bS^\ast}\E h_{\textsf{c}}(\hat \bs, \bs)\ge \sum_{i\in \mathfrak{N}^c}\inf_{\hat s_i}\frac{1}{|\bS^\ast|}\sum_{\bs\in \bS^\ast}\Prob(\hat s_i\ne s_i)
\end{equation}
where we consider a uniform prior on $\bS^\ast$ and hence the second inequality holds as minimax risk is lower bounded by Bayes risk, and the last inequality holds since the infimum can be taken over all $\hat \bs$ such that $\hat s_i=s^\ast_i$ for $i\in\mathfrak{N}$. Then it suffices to consider $\inf_{\hat s_i}\frac{1}{|\bS^\ast|}\sum_{\bs\in \bS^\ast}\Prob(\hat s_i\ne s_i)$ for $i\in\mathfrak{N}^c$. Without loss  generality, we assume $1\in \mathfrak{N}^c$ and for any $k\in[K]$ we denote $\bS_k^\ast=\{\bs\in \bS^\ast: s_1=k\}$. It's obvious that $\bS^\ast=\bigcup_{k=1}^K \bS_k^\ast$ and $\bS^\ast_a\bigcap \bS^\ast_b=\phi$ for $a\ne b$. In addition, by the definition of such partition, for any $a\ne b\in[K]$ and $\bs\in \bS^\ast_a$, there exists a unique $\bs^\prime\in \bS^\ast_b$ such that $s_i=s^\prime_i$ for all $i\ne 1$, which implies that $|\bS^\ast_a|=|\bS^\ast_b|$ for all $a,b\in[K]$. Then we have
\begin{align}\label{minimax_ineq2}
\inf_{\hat s_1}\frac{1}{|\bS^\ast|}\sum_{\bs\in \bS^\ast}\Prob(\hat s_1\ne s_1)&= \inf_{\hat s_1}\frac{1}{|\bS^\ast|}\frac{1}{K-1}\sum_{a<b}\left(\sum_{\bs\in \bS^\ast_a}\Prob(\hat s_1\ne a)+\sum_{\bs\in \bS^\ast_b}\Prob(\hat s_1\ne b)\right)\notag\\
&\geq\frac{1}{K(K-1)}\sum_{a<b}\inf_{\hat s_1}\left(\frac{1}{|\bS^\ast_a|}\sum_{\bs\in \bS^\ast_a}\Prob(\hat s_1\ne a)+\frac{1}{|\bS^\ast_b|}\sum_{\bs\in \bS^\ast_b}\Prob(\hat s_1\ne b)\right)\notag\\
&\geq\frac{1}{K(K-1)}\inf_{\hat s_1}\left(\frac{1}{|S_1^*|}\sum_{\bs\in S_1^*}\Prob(\hat s_1\ne 1)+\frac{1}{|S_2	^*|}\sum_{\bs\in S_2^*}\Prob(\hat s_1\ne 2)\right)\notag\\
&\geq\frac{1}{K(K-1)}\frac{1}{|\bS^\ast_{-1}|}\sum_{\bs_{-1}\in \bS^\ast_{-1}}\inf_{\hat s_1}\left(\Prob_{\bs=(1,\bs_{-1})}(\hat s_1\ne 1)+\Prob_{\bs=(2,\bs_{-1})}(\hat s_1\ne 2)\right)\notag\\
&\ge \frac{1}{K(K-1)}\inf_{\hat s_1}\Big(\Prob_{H_0^{(1)}}(\hat s_1= 2)+\Prob_{H_1^{(1)}}(\hat s_1=1)\Big)
\end{align}
where $\bS^\ast_{-1}$ is the collection of the subvectors in $\bS^\ast$ excluding the first coordinate, and we define a simple hypothesis testing for each $i\in[n]$:
\begin{equation*}
H_0^{(i)}: s_i=1 \text{~~~vs.~~~} H_1^{(i)}: s_i=2
\end{equation*}
Hence in \eqref{minimax_ineq2}, we have the form of Type-I error $+$ Type-II error of the above test. 
Notice that $|\{i\in[n]:s_i^*=k\}\backslash \mathfrak{N}_k |\ge \lfloor\alpha n/(4K^2)\rfloor$ and hence $|\mathfrak{N}^c|\geq c_0\alpha n/K $ for some  constant $c_0>0$. Combining this with \eqref{minimax_ineq}, \eqref{minimax_ineq2}, we proceed that
$$\inf_{\hat \bs}\sup_{\Omega}\E h_{\textsf{c}}(\hat \bs, \bs)\ge c_0\frac{\alpha n}{K^3}\frac{1}{|\mathfrak{N}^c|}\sum_{i\in \mathfrak{N}^c}\inf_{\hat s_i}\Big(\Prob_{H_0^{(l)}}(\hat s_i= 2)+\Prob_{H_1^{(l)}}(\hat s_i=1)\Big)$$
According to the Neyman-Pearson lemma, for each $i\in[n]$, the optimal test of $H_0^{(l)}\text{~vs.~} H_1^{(l)}$ is given by the likelihood ratio test with threshold $1$. Let $p_0(\bX_i)$ and $p_1(\bX_i)$ denote the likelihood of $\bX_i$ under $H_0$ and $H_1$, respectively. Then $\frac{p_1(\bX_i)}{p_0(\bX_i)}=\frac{\exp(\fro{\bX_i-\bM_1}^2/2)}{\exp(\fro{\bX_i-\bM_2}^2/2)}$ and hence the infimum is achieved by  $\hat s_i=\argmin_{k\in\{1,2\}}\fro{\bX_{i}-\bM_k}^2$. Therefore,
\begin{align*}
&\inf_{\hat s_i}\left(\frac{1}{2}\Prob_{H_0^{(l)}}(\hat s_i= 2)+\frac{1}{2}\Prob_{H_1^{(l)}}(\hat s_i=1)\right)\\
&=\frac{1}{2}\left(\Prob\left(\fro{\bM_1+\bE_i-\bM_2}^2\le \fro{\bE_i}^2\right)+\Prob\left(\fro{\bM_2+\bE_i-\bM_1}^2\le \fro{\bE_i}^2\right)\right)\\
&=\frac{1}{2}\left(\Prob\left(\frac{1}{2}\fro{\bM_1-\bM_2}^2\le \langle\bM_2-\bM_1,\bE_i\rangle\right)+\Prob\left(\frac{1}{2}\fro{\bM_1-\bM_2}^2\le \langle\bM_1-\bM_2,\bE_i\rangle\right)\right)
\end{align*}
Notice that $\langle\bM_2-\bM_1,\bE_i\rangle\overset{\text{d}}{=}\langle\bM_1-\bM_2,\bE_i\rangle\overset{\text{d}}{=}\calN(0,\sigma^2\fro{\bM_1-\bM_2}^2)$, we can proceed as
\begin{align*}
\inf_{\hat s_i}\left(\frac{1}{2}\Prob_{H_0^{(l)}}(\hat s_i= 2)+\frac{1}{2}\Prob_{H_1^{(l)}}(\hat s_i=1)\right)\geq \frac{\sigma}{\sqrt{2\pi}\fro{\bM_1-\bM_2}}\exp\left(-\frac{\fro{\bM_1-\bM_2}^2}{8\sigma^2}\right)
\end{align*}
where the inequality holds as $\fro{\bM_1-\bM_2}/\sigma\ge 1$. Hence we conclude that 
$$\inf_{\hat \bs}\sup_{\Omega}\E n^{-1}\cdot h_{\textsf{c}}(\hat \bs, \bs)\ge \exp\left(-\frac{\Delta^2}{8\sigma^2}-C\log\frac{\Delta K}{\alpha\sigma}\right)=\exp\left(-(1+o(1))\frac{\Delta^2}{8\sigma^2}\right)$$
provided that $\frac{\Delta^2}{\sigma^2\log(K/\alpha)}\rightarrow \infty$.

\
\subsection{Proof of Theorem \ref{thm:reduction}}
Suppose we are given the data $\{\bX_i\}_{i=1}^n$ generated by {eq:rank-one-model} with $((1-\epsilon)\bM,\bs^\ast) \in\wt\Omega_{\Lambda_{\submin}^{(n)}}$ for any $\epsilon\in(0,1]$. We utilize the sample splitting trick, similar to that in Theorem 2.4 in \cite{loffler2020computationally}, to generate two independent copies $\{\bX_i^{(1)}\}_{i=1}^n$ and $\{\bX_i^{(2)}\}_{i=1}^n$ by
$$\bX_i^{(1)}=\frac{\bX_i+\epsilon^{-1}\wt\bE_i}{\sqrt{1+\epsilon^{-2}}},\quad \bX_i^{(2)}=\frac{\bX_i-\epsilon\wt\bE_i}{\sqrt{1+\epsilon^2}}$$
for $i=1,\cdots,n$ where $\{\wt\bE_i\}_{i=1}^n$ are Gaussian noise matrices independent of $\{\bE_i\}_{i=1}^n$. As a consequence, we have $\bX_i^{(1)}=\frac{s_i^{\ast}\bM}{\sqrt{1+\epsilon^{-2}}}+\bE_i^{(1)}$ and $\bX_i^{(2)}=\frac{s_i^{\ast}\bM}{\sqrt{1+\epsilon^2}}+\bE_i^{(2)}$ with $\bE_i^{(1)}=\frac{\bE_i+\epsilon^{-1}\wt\bE_i}{\sqrt{1+\epsilon^{-2}}}$ and $\bE_i^{(2)}=\frac{\bE_i-\epsilon\wt\bE_i}{\sqrt{1+\epsilon^{2}}}$. Due to the property of Gaussian, $\{\bE_i^{(1)}\}_{i=1}^n$ and $\{\bE_i^{(2)}\}_{i=1}^n$ are independent. We define the following test statistic:
$$T_n=\op{\sum_{i=1}^n\frac{\hat s_i\bX_i^{(1)}}{n}}$$
where $(\hat s_1,\cdots,\hat s_n)=\hat\bs_{\textsf{comp}}(\bcalX^{(2)})$ with $\bcalX^{(2)}$ being the data tensor by stacking $\left\{\bX_i^{(2)}\right\}_{i=1}^n$. By construction, $\{\hat s_i\}_{i=1}^n$ is independent of $\left\{\bE_i^{(1)}\right\}_{i=1}^n$ and hence $\sum_{i=1}^n\frac{\hat s_i\bE_i^{(1)}}{n}\overset{d}{=}\sum_{i=1}^n\frac{\bE_i^{(1)}}{n}$. Under $H_0$, with probability at least $1-\exp(-d)$:
$$T_n=\op{\sum_{i=1}^n\frac{\hat s_i\bX_i^{(1)}}{n}}\le \frac{C_0}{2}\sqrt{\frac{d}{n}}$$
for some absolute constant $C_0>0$. Under $H_1$, we have  
$((1+\epsilon^2)^{-1/2}\bM,\bs^\ast)\in\wt\Omega_{\Lambda_{\submin}^{(n)}}$ since $(1-\epsilon)\le (1+\epsilon^2)^{-1/2}$. By \eqref{eq:reduction-assump} we have that with probability greater than $1-\zeta_n$:
\begin{equation}\label{eq:reduction-assumption-res}
	n^{-1}\cdot h_{\textsf{c}}(\hat \bs_{\textsf{comp}},\bs^\ast)\le \delta_n
\end{equation}
Without loss of generality we assume $h_{\textsf{c}}(\hat \bs_{\textsf{comp}},\bs^\ast)=h(\hat \bs_{\textsf{comp}},\bs^\ast)$. Hence we can obtain with probability at least $1-\zeta_n-\exp(-d)$:
\begin{align*}
	T_n&\ge \op{\sum_{i=1}^n\frac{\hat s_is_i^\ast}{n\sqrt{1+\epsilon^{-2}}}\bM}-\op{\sum_{i=1}^n\frac{\hat s_i\bE_i^{(1)}}{n}}\\
	&\ge \frac{\Lambda_{\submin}^{(n)}(1-2n^{-1} h(\hat\bs_{\textsf{comp}},\bs^\ast))}{\sqrt{n}\cdot \sqrt{1+\epsilon^{-2}}}-\frac{C_0}{2}\sqrt{\frac{d}{n}}\\
	&> \frac{C_0}{2}\sqrt{\frac{d}{n}}
\end{align*}
where we've used \eqref{eq:reduction-assumption-res} and $\Lambda_{\submin}^{(n)}> C_0(1-2\delta_n)^{-1}\sqrt{1+\epsilon^{-2}}d^{1/2}$ in the last inequality. Then the test $\phi_n$ can be defined as 
\begin{align*}
	\phi_n(\bcalX)=\begin{cases}
		1 & \text{if~} T_n> C_0\sqrt{\frac{d}{n}},\\
		0 & \text{otherwise.}
	\end{cases}
\end{align*}
It turns out that 
$$\EE_{Q_n}[\phi_n(\bcalX)]+\sup_{((1-\epsilon)\bM,\bs^\ast)\in \wt\Omega_{\Lambda_{\submin}^{(n)}}}\EE_{(\bM,\bs^\ast)}[1-\phi_n(\bcalX)]\le \zeta_n+\exp(-d)$$
Notice that computing $T_n$ requires only $poly(d,n)$ and the proof is completed  by setting $n,d\rightarrow \infty$.

\subsection{Proof of Theorem \ref{thm:weak-SNR}}
Theorem \ref{thm:weak-SNR-init} can be obtained by modifying the proofs of Theorem \ref{thm:main}, and hence we only sketch the necessary modifications here. Similar to the proof of Theorem \ref{thm:main}, we  have
\begin{align}\label{ell-init-cond}
h(\hat \bs^{(0)}, \bs^{\ast})\leq \frac{\ell(\hat\bs^{(0)},\bs^\ast)}{\Delta^2} =o\left(\frac{\alpha n}{\kappa_0^2 }\right)
\end{align}
as a consequence of condition \eqref{relaxed-init-cond}. 
\\We consider the iterative convergence of Algorithm \ref{alg:weak-SNR-Lloyd}. Following the same argument of Step 2 in the proof of Theorem \ref{thm:main} and adopting the same notation therein, we have the following inequality:
\begin{align*}
	\ell(\hat\bs^{(t)},\bs^\ast)\le  \xi_{\textsf{err}}+\beta_1(\bs^\ast,\hat\bs^{(t)})+\beta_2(\bs^\ast,\hat\bs^{(t)})
\end{align*}
We can bound $\xi_{\textsf{err}}$ the same as \textbf{Step 2.1} in the proof of Theorem \ref{thm:main}. To bound $\beta_1(\bs^{\ast},\hat\bs^{(t)})$, it turns out that, by symmetry, we only need to bound 
\begin{align}\label{eq:relaxed-beta12}
	\beta_{1,2}(\bs^{\ast}, \hat\bs\supt)&:=\sum_{i=1}^n\fro{\bM_1-\bM_{2}}^2\ind{\hat s\supt_i\neq 1}
	\cdot\ind{\inp{\bE_i}{\hat\bM\supt_{1}-\bM_1}\ge\frac{\delta}{8}\fro{\bM_{2}-\bM_{1}}^2}\notag\\
	&+\sum_{i=1}^n\fro{\bM_2-\bM_{1}}^2\ind{\hat s\supt_i\neq 2}
	\cdot\ind{\inp{\bE_i}{\hat\bM\supt_{2}-\bM_2}\ge\frac{\delta}{8}\fro{\bM_{1}-\bM_{2}}^2}
\end{align}
The argument in \textbf{Step 2.2} in the proof of Theorem \ref{thm:main} can be directly applied to the analysis of $\hat\bM_1-\bM_1$, i.e., the first term on RHS of eq. \eqref{eq:relaxed-beta12}, whereas it fails for $\hat\bM_2-\bM_2$ since $\sigma_{\min}(\bM_2)$ {\it can} be arbitrarily close to $0$ and Lemma \ref{lem:pertub-formula} no longer holds. Observe that
\begin{align*}
	\hat\bM_{2}\supt&=\hat\bU_2\hat\bU_2^\top\left(\frac{1}{n_2\suptt}\sum_{i=1}^n{\ind{\hat s_i^{(t-1)}=2}\bM_{s_i^\ast}}+\bEs\suptt\right)\hat\bV_2\hat\bV_2^\top\\
	&=\hat\bU_2\hat\bU_2^\top\left[\bM_2+\frac{1}{n_a\suptt}\sum_{i=1}^n{\ind{\hat s_i^{(t-1)}=2}(\bM_{s_i^\ast}}-\bM_2)+\bEss+(\bEs\suptt-\bEss)\right]\hat\bV_2\hat\bV_2^\top\\
	&=\hat\bU_2\hat\bU_2^\top\left(\bM_2+\bEss+\Delta_\bM\suptt+\Delta_\bE\suptt\right)\hat\bV_2\hat\bV_2^\top
\end{align*}
where 
$$
\Delta_\bM\suptt=\frac{1}{n_2\suptt}\sum_{i=1}^n{\ind{\hat s_i^{(t-1)}=2}(\bM_{s_i^\ast}}-\bM_2)\quad {\rm and}\quad \Delta_\bE\suptt=\bEs\suptt-\bEss
$$
Notice that since $h(\hat\bs^{(t-1)},\bs^\ast)$ satisfies \eqref{eq:weakSNR-h-init-cond}, we have $n_2\suptt\ge {7\alpha n}/{16}$. 
Lemma \ref{lem:DeltaM_bound} implies that under event $\calQ_1\cap \calQ_2$, we have
\begin{align*}
	\op{\hat\bM_{2}\supt}&\le (1+c)\op{\bM_2}+c \left(\alpha^{-1/2}\sqrt{\frac{d}{ n}}+\alpha^{-1}\sqrt{\frac{h(\hat\bs^{(t-1)},\bs^\ast)}{n}}\right)\\
	&\le c^\prime \left(\alpha^{-1/2}\sqrt{\frac{d}{n}}+\alpha^{-1/2}\kappa_0^{-1}\right) 
\end{align*}
for some small universal constant $c^\prime>0$, where the second inequality is due to Assumption \ref{assump:relaxed-lloyd-SNR}. On the other hand, under event $\calQ_1\cap \calQ_2$ and Assumption \ref{assump:relaxed-lloyd-SNR} we also have
\begin{align*}
	\op{\hat\bM_{1}\supt}\ge (1-c)\op{\bM_1}-c^\prime\left(\alpha^{-1/2}\sqrt{\frac{d}{ n}}+\alpha^{-1/2}\kappa_0^{-1}\right)>\op{\hat\bM_{2}\supt}
\end{align*}
By taking a union bound over $\calQ_1\cap \calQ_2$, we conclude that with probability at least $1-\exp(-cd)$ we have $\op{\hat\bM_{2}\supt}<\op{\hat\bM_{1}\supt}$ and hence we set ${\hat\bM_{2}\supt}=0$ afterwards. Then for the second term on RHS of eq. \eqref{eq:relaxed-beta12}, we have
\begin{align*}
	\PP&\left(\inp{\bE_i}{\hat\bM\supt_{2}-\bM_2}\ge\frac{\delta}{8}\fro{\bM_{1}-\bM_{2}}^2\right)=	\PP\left(\inp{\bE_i}{-\bM_2}\ge\frac{\delta}{8}\fro{\bM_{1}-\bM_{2}}^2\right)\\
	&\le \exp\left( -\frac{\delta^2\fro{\bM_{1}-\bM_{2}}^4}{128\fro{\bM_2}^2}\right)\le \exp\left( -c\frac{\lambda_1^2r_1}{\op{\bM_2}^2r_2}\delta^2\fro{\bM_{1}-\bM_{2}}^2\right)
\end{align*}
where the last inequality is due to Assumption \ref{assump:weak-SNR}. Hence the expecatation can be bounded as 
\begin{align*}
\EE\bigg[ &\sum_{i=1}^n\fro{\bM_2-\bM_{1}}^2\ind{\hat s\supt_i\neq 2}
	\cdot\ind{\inp{\bE_i}{\hat\bM\supt_{2}-\bM_2}\ge\frac{\delta}{8}\fro{\bM_{1}-\bM_{2}}^2}\bigg]\\
	&\quad \leq n\Delta^2\exp\left[-c\delta^2r_1r_2^{-1}\left(\lambda_1/\op{\bM_2}\right)^2\Delta^2\right]
\end{align*}
By Markov inequality, with probability at least $1-\exp\left[-\delta\left(\sqrt{r_1/r_2}\lambda_1/\op{\bM_2}\right)\Delta\right]$ we get  
\begin{align*}
\sum_{i=1}^n&\fro{\bM_2-\bM_{1}}^2\ind{\hat s\supt_i\neq 2}
	\cdot\ind{\inp{\bE_i}{\hat\bM\supt_{2}-\bM_2}\ge\frac{\delta}{8}\fro{\bM_{1}-\bM_{2}}^2}
	&\leq n\cdot \exp\left(-\delta (\alpha n/K)^{1/2}\Delta^2\right)\\
	&\leq n\cdot \exp\left[-\delta^2r_1r_2^{-1}\left(\lambda_1/\op{\bM_2}\right)^2\Delta^2\right]
\end{align*}
which holds as long as $\delta\to 0$ sufficiently slowly compared with $\lambda_1^2r_1r_2^{-1}/\op{\bM_2}^2\to\infty$. 
\\It remains to consider $\beta_2(\bs^{\ast},\hat\bs^{(t)})$. Observe that 
\begin{align}
\beta_2(\bs^{\ast}, \hat\bs\supt)\leq &\sum_{i=1}^n\sum_{a\in[2]\backslash\{ s_i^\ast\}}\ind{\hat s\supt_i\neq a}\fro{\bM_a-\bM_{s_i^\ast}}^2\ind{\frac{1}{2}\fro{\bM_{s_i^\ast}-\hat\bM_{s_i^\ast}\supt}^2\ge\frac{\delta}{12}\fro{\bM_{s_i^\ast}-\bM_{a}}^2}\notag\\
+ &\sum_{i=1}^n\sum_{a\in[2]\backslash\{ s_i^\ast\}}\ind{\hat s\supt_i\neq a}\fro{\bM_a-\bM_{s_i^\ast}}^2\ind{\frac{1}{2}\fro{\bM_{s_i^\ast}-\hat\bM_a\supt}^2\ge\frac{\delta}{12}\fro{\bM_{s_i^\ast}-\bM_{a}}^2}\notag\\
+ &\sum_{i=1}^n\sum_{a\in[2]\backslash\{ s_i^\ast\}}\ind{\hat s\supt_i\neq a}\fro{\bM_a-\bM_{s_i^\ast}}^2\ind{\fro{\bM_{s_i^\ast}-\bM_a}\fro{\bM_{a}-\hat\bM_a\supt}\ge\frac{\delta}{12}\fro{\bM_{s_i^\ast}-\bM_{a}}^2}\label{eq:relaxed-beta_2}
\end{align}
The first term on RHS of eq. \eqref{eq:relaxed-beta_2} can be written as 
\begin{align}
	&\sum_{i=1}^n\sum_{a\in[2]\backslash\{ s_i^\ast\}}\ind{\hat s\supt_i\neq a}\fro{\bM_a-\bM_{s_i^\ast}}^2\ind{\frac{1}{2}\fro{\bM_{s_i^\ast}-\hat\bM_{s_i^\ast}\supt}^2\ge\frac{\delta}{12}\fro{\bM_{s_i^\ast}-\bM_{a}}^2}\notag\\
	&=\sum_{i=1}^n\ind{\hat s\supt_i\neq 1}\fro{\bM_1-\bM_{2}}^2\ind{\frac{1}{2}\fro{\bM_{2}-\hat\bM_{2}\supt}^2\ge\frac{\delta}{12}\fro{\bM_{2}-\bM_{1}}^2}\notag\\
	&+\sum_{i=1}^n\ind{\hat s\supt_i\neq 2}\fro{\bM_2-\bM_{1}}^2\ind{\frac{1}{2}\fro{\bM_{1}-\hat\bM_{1}\supt}^2\ge\frac{\delta}{12}\fro{\bM_{1}-\bM_{2}}^2}\label{eq:relaxed-beta_21}
\end{align}
The second term of \eqref{eq:relaxed-beta_21} can be bounded the same way as that in \textbf{Step 2.3} of the proof of Theorem \ref{thm:main}. Note that  $\fro{\bM_{2}-\hat\bM_{2}\supt}^2=\fro{\bM_{2}}^2\le r_2\op{\bM_2}^2=o(r_1\lambda_1^2)=o\left(\fro{\bM_{1}-\bM_{2}}^2\right)$ and hence  the first term of \eqref{eq:relaxed-beta_21} vanishes by setting $\delta$ slowly converging to $0$. It suffices to consider the last term on RHS of eq. \eqref{eq:relaxed-beta_2}. Observe that 
\begin{align}
&\sum_{i=1}^n\sum_{a\in[2]\backslash\{ s_i^\ast\}}\ind{\hat s\supt_i\neq a}\fro{\bM_a-\bM_{s_i^\ast}}^2\ind{\fro{\bM_{s_i^\ast}-\bM_a}\fro{\bM_{a}-\hat\bM_a\supt}\ge\frac{\delta}{12}\fro{\bM_{s_i^\ast}-\bM_{a}}^2}\notag\\
&=\sum_{i=1}^n\ind{\hat s\supt_i\neq 1}\fro{\bM_1-\bM_{2}}^2\ind{\fro{\bM_{2}-\bM_1}\fro{\bM_{1}-\hat\bM_1\supt}\ge\frac{\delta}{12}\fro{\bM_{2}-\bM_{1}}^2}\notag\\
&+\sum_{i=1}^n\ind{\hat s\supt_i\neq 2}\fro{\bM_2-\bM_{1}}^2\ind{\fro{\bM_{1}-\bM_2}\fro{\bM_{2}-\hat\bM_2\supt}\ge\frac{\delta}{12}\fro{\bM_{1}-\bM_{2}}^2}\label{eq:relaxed-beta_23}
\end{align}
The first term of \eqref{eq:relaxed-beta_23} can be bounded the same way as that in \textbf{Step 2.3} of the proof of Theorem \ref{thm:main}, and the second term vanishes as $\fro{\bM_{2}-\hat\bM_{2}\supt}=\fro{\bM_{2}}=o\left(\fro{\bM_{1}-\bM_{2}}\right)$. \\
By mimicing the remaining proofs of Theorem \ref{thm:main}, we can finish the proof of Theorem \ref{thm:weak-SNR}.

\subsection{Proof of Theorem \ref{thm:weak-SNR-init}}
For notational simplicity, we denote the smallest non-trivial singular value of $\bM_1$ as $\lambda_1$. 
Denote the following decomposition of tensor $\bcalM=\bcalM_1+\bcalM_2$, where for $k\in[2]$, the $i$-th slice of $\bcalM_k$ is defined as $[\bcalM_k]_{\cdot\cdot i}=\II(\bs_i^\ast=k)\bM_k$. It turns out that $\bU_1$ is the leading-$r_1$ left singular vectors of $\scrM_1(\bcalM_1)$ and $\bV_1$ is the leading-$r_1$ left singular vectors of $\scrM_2(\bcalM_1)$. We first show that $\hat\bU_1$ and $\hat\bV_1$ are close to $\bU_1$ and $\bV_1$, respectively.  Without loss of generality, we only consider $\hat\bU_1$. A key observation is that $\hat\bU_1$ is also the leading-$r_1$ left eigenvectors $\scrM_1(\bcalX)\scrM_1^\top (\bcalX)$. Then write
\begin{align}\label{eq:sec-mom-decomp}
\scrM_1(\bcalX)\scrM_1^\top(\bcalX)&=\scrM_1(\bcalM)\scrM_1^\top(\bcalM)+\scrM_1(\bcalM)\scrM_1^\top(\bcalE)+\scrM_1(\bcalE)\scrM_1^\top(\bcalM)+\scrM_1(\bcalE)\scrM_1^\top(\bcalE)\notag\\
&=\scrM_1(\bcalM_1)\scrM_1^\top(\bcalM_1)+\scrM_1(\bcalM_1)\scrM_1^\top(\bcalM_2)+\scrM_1(\bcalM_2)\scrM_1^\top(\bcalM_1)\notag\\
&+\scrM_1(\bcalM_2)\scrM_1^\top(\bcalM_2)+\left[\scrM_1(\bcalM_1)+\scrM_1(\bcalM_2)\right]\scrM_1^\top(\bcalE)\notag\\
&+\scrM_1(\bcalE)\left[\scrM_1(\bcalM_1)+\scrM_1(\bcalM_2)\right]^\top+\scrM_1(\bcalE)\scrM_1^\top(\bcalE)
\end{align}
We are going to bound each term on RHS of eq. \eqref{eq:sec-mom-decomp}. The first term $\scrM_1(\bcalM_1)\scrM_1^\top(\bcalM_1)$ is the signal part and we have
$$\sigma_{\min}(\scrM_1(\bcalM_1)\scrM_1^\top(\bcalM_1))=\sigma_{r_1}(\scrM_1(\bcalM_1)\scrM_1^\top(\bcalM_1))\ge n_1^\ast\lambda_1^2$$
For the $2$nd, $3$rd and $4$th term of \eqref{eq:sec-mom-decomp}, we can have 
$$\op{\scrM_1(\bcalM_1)\scrM_1^\top(\bcalM_2)+\scrM_1(\bcalM_2)\scrM_1^\top(\bcalM_1)}\le 2\kappa_0\sqrt{n_1^\ast n_2^\ast} \lambda_1\op{\bM_2}$$
and
$$\op{\scrM_1(\bcalM_2)\scrM_1^\top(\bcalM_2)}\le n_2^\ast\op{\bM_2}^2$$
The $5$th and $6$th term of eq. \eqref{eq:sec-mom-decomp} can be together bounded as 
\begin{align*}
	&\op{\left[\scrM_1(\bcalM_1)+\scrM_1(\bcalM_2)\right]\scrM_1^\top(\bcalE)+\scrM_1(\bcalE)\left[\scrM_1(\bcalM_1)+\scrM_1(\bcalM_2)\right]^\top}\\
	&\le C\left(\kappa_0\sqrt{n_1^\ast}\lambda_1+\sqrt{n_2^\ast}\op{\bM_2}\right)\sqrt{d}
\end{align*}
with probability at least $1-\exp(-cd)$, for some absolute constant $c,C>0$.
Lastly, we notice that $\EE\left(\scrM_1(\bcalE)\scrM_1^\top(\bcalE)\right)=nd_2\bI_{d_1}$, then by \cite{koltchinskii2017concentration}, with probability at least $1-\exp(-d)$ we have 
$$\op{\scrM_1(\bcalE)\scrM_1^\top(\bcalE)-nd_2\bI_{d_1}}\le C\sqrt{n}d$$
Note that $n_2^\ast/n_1^\ast\le 2(1-\alpha/2)/\alpha\le 2\alpha^{-1}$. Collecting all pieces  above, if 
\begin{align}\label{eq:cond-relaxed-init}
	\lambda_1\ge C\left(\kappa_0\alpha^{-1/2}\sqrt{\frac{d}{n}}+\alpha^{-1/2}\frac{d^{1/2}}{n^{1/4}}\right),\quad \lambda_1\ge \kappa_0\alpha^{-1/2}\op{\bM_2} 
\end{align} 
for some large constant $C>0$. Note that $\kappa_0\alpha^{-1/2}\sqrt{\frac{d}{n}}$ in the first condition in \eqref{eq:cond-relaxed-init}  is trivial as we assume $n/\kappa_0^4\ge C$ for some large constant $C>0$, and the second term is implied by the condition on $\sigma_{r_1}(\bM_1)$ together with the assumption
$$\op{\bM_2}\le C\kappa_0^{-1}\frac{d^{1/2}}{n^{1/4}}$$ 
Then with probability greater than $1-\exp(-cd)$ we can have
$\op{\hat\bU_1\hat\bU_1^\top-\bU_1\bU_1^\top}\le 1/4$. Using same analysis on $\hat\bV_1$, we can conclude with probability at least $1-\exp(-cd)$:
\begin{equation}\label{eq:relaxed-uv-small}
	\max\left\{\op{\hat\bU_1\hat\bU_1^\top-\bU_1\bU_1^\top},\op{\hat\bV_1\hat\bV_1^\top-\bV_1\bV_1^\top}\right\}\le \frac{1}{6}
\end{equation}
Define $\hat\bcalG=\bcalX\times_1\hat \bU_1\hat \bU_1^\top \times_2\hat\bV_1\hat \bV_1^\top $, $\bcalG:=\bcalM\times_1\hat \bU_1\hat \bU_1^\top \times_2\hat\bV_1\hat \bV_1^\top $ (also $\bG:=\scrM_3(\bcalG)$) and $\bfrakM:=[vec(\hat \bM_{\hat s_1^{(0)}}),\cdots,vec(\hat \bM_{\hat s_n^{(0)}})]^\top \in\RR^{n\times d_1d_2}$.  We can have the following lemma, which is an analogue to Lemma \ref{lem:init-two-facts}.
\begin{lemma}\label{lem:relaxed-two-facts} Suppose \eqref{eq:relaxed-uv-small} holds. Then we have the following facts:
	\begin{enumerate}
	\item[(I)] $\bfrakM$, the k-means solution,  is close $\bG$, i.e.,  there exists some absolute constants $c_0,C_0>0$ such that with probability at least $1-\exp(-c_0d)$:
	$$\fro{\bfrakM-\bG}\le C_0\left(\sqrt{dr_1+n}\right)$$
	\item[(II)] The rows of $\bG$ belonging to different clusters is well-separated, i.e.
	\begin{align*}
\fro{\bcalG\times_3 (\be_i^\top-\be_j^\top)}\ge \frac{\Delta}{2}
\end{align*}
for any $i,j\in[n],s_i^\ast\ne s_j^\ast$.
\end{enumerate}
\end{lemma}
Following the almost identical argument in the proof of Theorem \ref{thm:spec-initialization} but replacing $\hat\bU$ with $\hat\bU_1$ and $\hat\bV$ with $\hat\bV_1$, with probability at least $1-\exp(-cd)$ we have
\begin{equation}\label{eq:weakSNR-h-init-cond}
	n^{-1}\cdot h_{\textsf{c}}(\hat\bs^{(0)},\bs^\ast)\le \frac{C}{\Delta^2}\left(\frac{dr_1}{n}+1\right)=o\left(\frac{\alpha }{\kappa_0^2}\right)
\end{equation}
where the last equality holds provided that $\Delta^2\gg \kappa_0^2\alpha^{-1}\left({dr_1}/{n}+1\right)$. Since the condition in Theorem \ref{thm:weak-SNR} already implies that
\begin{align*}
	\Delta^2\gtrsim r_1\lambda_1^2\ge C\alpha^{-1}\frac{{dr_1}}{\sqrt{n}}
\end{align*}
Then if $ n/\kappa_0^4\rightarrow \infty$ and $\alpha\Delta^2/\kappa_0^2\rightarrow \infty$, the condition $\Delta^2\gg \kappa_0^2\alpha^{-1}\left({dr_1}/{n}+1\right)$ automatically holds.   
\\As a result, we also have the following holds with probability at least $1-\exp(-c d)$:
\begin{align*}
\ell_{\textsf{c}}(\hat \bs^{(0)}, \bs^{\ast})\leq \Delta^2 h_{\textsf{c}}(\hat\bs^{(0)},\bs^\ast)=o\left(\frac{\alpha n\Delta^2}{\kappa_0^2 }\right)
\end{align*}
which is an analogue to \eqref{ell-init-cond}, where we've used $\gamma=1$ in the two component case. 

\section{Proof of Technical Lemmas}\label{sec:technical_lemmas}

\subsection{Proof of Lemma~\ref{lem:sv-lower-bound}}
Without loss of generality we only proof $j=1$. It follows that
\begin{align*}
\sigma^2_{\min}(\scrM_1(\bcalM))&\ge \kappa_1^{-2}\op{\scrM_1(\bcalM)}^2\ge \kappa_1^{-2}r_\bU^{-1}\sum_{k=1}^Kn_k\fro{\bM_k}^2\\
 &\ge \kappa_1^{-2}r_\bU^{-1}n\lambda^2\ge \kappa_1^{-2}(Kr)^{-1}n\lambda^2
\end{align*}
where the last inequality is due to $r_\bU\le \sum_{k=1}^Kr_k\le Kr$.

\subsection{Proof of Lemma \ref{lem:kappa12}} By definition we have that 
$$\bU^\top\bU=\mat{
\bI_{r_1}&\bU_1^\top \bU_2&\cdots&\bU_1^\top \bU_K\\
\bU_2^\top \bU_1&\bI_{r_2} &\cdots&\bU_2^\top\bU_K\\
\vdots&\vdots&\ddots&\vdots\\
\bU_K^\top\bU_1&\bU_K^\top\bU_2&\cdots &\bI_{r_K}
}$$ 
and $\bW^\top\bW=\text{diag}(n_1^\ast,\cdots,n_K^\ast)$. Hence we have 
$$\bW^\top \bW\otimes\bV^\top \bV=\mat{
n_1^\ast\bU^\top\bU&\mathbf{0}&\cdots&\mathbf{0}\\
\mathbf{0}&n_2^\ast\bU^\top\bU&\cdots&\mathbf{0}\\
\vdots&\vdots&\ddots&\vdots\\
\mathbf{0}&\mathbf{0}&\cdots &n_K^\ast\bU^\top\bU
}$$
Simple calculations give that
\begin{align*}
	\scrM_1(\bcalM)\scrM^\top_1(\bcalM)&=\bU\scrM_1(\bcalS)(\bW^\top \bW\otimes\bV^\top \bV)\scrM_1^\top(\bcalS)\bU^\top\\
	&=\bU\cdot \text{diag}(n_1^\ast\bSigma_1^2,\cdots,n_K^\ast\bSigma_K^2) \cdot \bU^\top 
\end{align*}
As a result, we obtain 
$$\sigma_1(\scrM_1(\bcalM)\scrM^\top_1(\bcalM))\le \sigma_1^2(\bU)\cdot \max_{1\le k\le K}{n_k^\ast\sigma^2_{\submax}(\bSigma_k)}$$
$$\sigma_{r_\bU}(\scrM_1(\bcalM)\scrM^\top_1(\bcalM))\geq \sigma_{r_\bU}^2(\bU)\cdot \min_{1\le k\le K}{n_k^\ast\sigma^2_{\submin}(\bSigma_k)}$$
Hence we conclude that 
$$\kappa_1=\sqrt{\frac{\sigma_1(\scrM_1(\bcalM)\scrM^\top_1(\bcalM))}{\sigma_{r_\bU}(\scrM_1(\bcalM)\scrM^\top_1(\bcalM))}}\le \kappa_0\kappa(\bU)\cdot \sqrt{\frac{n_\submax^\ast}{n_\submin^\ast}}$$
Similarly we can prove that $\scrM_2(\bcalM)\scrM^\top_2(\bcalM)=\bV\cdot \text{diag}(n_1^\ast\bSigma_1^2,\cdots,n_K^\ast\bSigma_K^2) \cdot \bV^\top $  and $\kappa_1\le \kappa_0\kappa(\bU)\cdot ({n_\submax^\ast}/{n_\submin^\ast})^{1/2}$.
\\If $r_\bU=r_\bV=r_1$, by min-max principle for singular values we have 
$$\sigma_\submin(\bU)=\sigma_{r_1}(\bU)=\max_{S\subset \RR^{n},\text{dim}(S)=r_1}\min_{x\in S, \op{x}=1}\op{\mat{\bU_1^\top x\\ \vdots\\\bU_K^\top x}}\ge \max_{S\subset \RR^{n},\text{dim}(S)=r_1}\min_{x\in S, \op{x}=1}\op{\bU_1^\top x}=\sigma_\submin(\bU_1)=1$$
and 
$$\sigma_\submax(\bU)=\max_{x\in \RR^n, \op{x}=1}\op{\mat{\bU_1^\top x\\ \vdots\\\bU_K^\top x}}\leq\sqrt{\sum_{k=1}^K \max_{x\in \RR^n, \op{x}=1}\op{\bU_k^\top x}}=\sqrt{K}$$
Therefore, we have $\kappa(\bU)\le K^{1/2}$ and similarly $\kappa(\bV)\le K^{1/2}$, from which we can  conclude that $\max\{\kappa_1,\kappa_2\}\le \kappa_0 (K^2/\alpha)^{1/2}$.
\\If $r_\bU=r_\bV=\mathring{r}$ and $\bU_k$'s are mutually orthogonal, then $\bU,\bV$ has orthonormal columns and $\kappa(\bU)=\kappa(\bV)=1$. Hence we have  $\max\{\kappa_1,\kappa_2\}\le \kappa_0 (K/\alpha)^{1/2}$.
\subsection{Proof of Lemma \ref{lem:concentration}} Note that for fixed $k\in[K]$, we have $\frac{\sum_{i=1}^n{\ind{ s_i^{\ast}=k}\bE_i}}{\sum_{i=1}^n{\ind{ s_i^{\ast}=k}}}$ has i.i.d. sub-gaussian entries with mean zero and variance $(n_k^\ast)^{-1}$. By random matrix theory there exists some absolute constants $c,C>0$ such that 
\begin{align*}
	\Prob\left(\op{\frac{\sum_{i=1}^n{\ind{ s_i^{\ast}=k}\bE_i}}{\sum_{i=1}^n{\ind{ s_i^{\ast}=k}}}}\ge C\sqrt{\frac{d}{n^\ast_{k}}}\right)\le \exp(-cd)
\end{align*}
Applying a union bound over $[K]$ gives 
\begin{align*}
	\Prob(Q_1^c)=\Prob\left(\bigcup_{k=1}^K\left\{\op{\frac{\sum_{i=1}^n{\ind{ s_i^{\ast}=k}\bE_i}}{\sum_{i=1}^n{\ind{ s_i^{\ast}=k}}}}\ge C\sqrt{\frac{d}{n^\ast_{k}}}\right\}\right)\le K\exp(-cd)\le \exp(-c_0d)
\end{align*}
for some absolute constant $c_0>0$, provided that $d\gtrsim \log K$. To prove the tail bound for $Q_2$, consider fixed set $I\subseteq [n]$, we have for any $t>0$:
\begin{align*}
	\Prob\left(\op{\frac{1}{\sqrt{|I|}}\sum_{i\in I}\bE_i}\le C\left(\sqrt{d}+t\right)\right)\le 2\exp(-t^2)
\end{align*} 
Applying a union bound over all subsets of $[n]$ gives 
\begin{align*}
	\Prob(Q_2^c)=\Prob\left(\bigcup_{I\subseteq [n]}\left\{\op{\frac{1}{\sqrt{|I|}}\sum_{i\in I}\bE_i}\le C\left(\sqrt{d}+t\right)\right\}\right)\le 2\exp(-t^2+n)
\end{align*} 
By choosing $t=C_1\left(\sqrt{n}+\sqrt{d}\right )$ for some absolute constant $C_1>0$, we obtain the desired result. It suffices to prove the bound for $\calQ_3$. Fix $i\in[n]$, then for any $t>0$:
\begin{align*}
	\Prob\left(\op{\frac{\sum_{j\ne i}^n{\ind{ s_j^{\ast}=a}\bE_j}}{\sum_{j=1}^n{\ind{ s_j^{\ast}=a}}}}\ge C\sqrt{\frac{d+t^2}{n^\ast_{a}}}\right)\le 2\exp\left (-t^2\right )
\end{align*}
and
\begin{align*}
	\Prob\left(\op{\bE_i}\ge C\sqrt{d+t^2}\right)\le 2\exp\left (-t^2\right )
\end{align*}
Applying a union bound over $[n]$ and $[K]$ gives 
\begin{align*}
	\Prob\left(\bigcup_{a=1}^K\bigcup_{i=1}^n\left\{\op{\frac{\sum_{j\ne i}^n{\ind{ s_j^{\ast}=a}\bE_j}}{\sum_{j=1}^n{\ind{ s_j^{\ast}=a}}}\right\}}\ge C\sqrt{\frac{d+t^2}{n^\ast_{a}}}\right)\le 2nK\exp\left (-t^2\right )
\end{align*}
and
\begin{align*}
	\Prob\left(\bigcup_{i=1}^n\left\{\op{\bE_i}\ge C\sqrt{d+t^2}\right\}
	\right)\le 2n\exp\left (-t^2\right )
\end{align*}
We can take $t=C_2\sqrt{d}+\log n$ for some absolute constant $C_2>0$ (using $d\gtrsim \log K$) and the proof is completed.  
\subsection{Proof of Lemma~\ref{lem:DeltaM_bound}}
By definition, $\op{\Delta_\bM\suptt}$ can be bounded by
\begin{align*}
	\op{\Delta_\bM\suptt}&=\op{\frac{1}{n_a\suptt}\sum_{i=1}^n{\ind{\hat s_i^{(t-1)}=a}(\bM_{s_i^\ast}}-\bM_a)}\\
	&=\op{\frac{1}{n_a\suptt}\sum_{i=1}^n{\ind{\hat s_i^{(t-1)}=a,s_i^\ast\ne a }(\bM_{s_i^\ast}}-\bM_a)}\\
	&\le \frac{8K}{7\alpha n}\sum_{i=1}^n\ind{\hat s_i^{(t-1)}=a,s_i^\ast\ne a}\op{\bM_{s_i^\ast}-\bM_a}\\
	&\le \frac{8K}{7\alpha n}\cdot \frac{\ell_a(\hat\bs^{(t-1)},\bs^\ast)}{\Delta}
\end{align*}
An alternative bound for $\op{\Delta_\bM\suptt}$:
\begin{align*}
	\op{\Delta_\bM\suptt}&=\op{\frac{1}{n_a\suptt}\sum_{i=1}^n{\ind{\hat s_i^{(t-1)}=a}(\bM_{s_i^\ast}}-\bM_a)}\\
	&=\op{\frac{1}{n_a\suptt}\sum_{i=1}^n{\ind{\hat s_i^{(t-1)}=a,s_i^\ast\ne a }(\bM_{s_i^\ast}}-\bM_a)}\\
	&\le \frac{8K}{7\alpha n}\sum_{i=1}^n\ind{\hat s_i^{(t-1)}=a,s_i^\ast\ne a}\op{\bM_{s_i^\ast}-\bM_a}\\
	&\le \frac{16\kappa_0K}{7\alpha n}\cdot \lambda\cdot h_a(\hat\bs^{(t-1)},\bs^\ast)
\end{align*}
where we've used $h_a(\hat\bs^{(t-1)},\bs^\ast)\le \sum_{a\in[K]}h_a(\hat\bs^{(t-1)},\bs^\ast)=h(\hat\bs^{(t-1)},\bs^\ast)$
and the condition \eqref{init-cond}. 
In other words, we have the following bound for $\Delta_\bM^{(t-1)}$ that will be utilized repeatedly later:
\begin{align}\label{deltaM-bound}
	\op{\Delta_\bM\suptt}&\leq \frac{16K}{7\alpha n}\cdot\min\left\{\kappa_0 \lambda h_a(\hat\bs^{(t-1)},\bs^\ast),\frac{\ell _a(\hat\bs^{(t-1)},\bs^\ast)}{\Delta}\right\}
\end{align} 
Moreover, under $\calQ_1$ we have
\begin{align*}
	\op{\bEsa}\lesssim\sqrt{\frac{d}{n_a^\ast}}\lesssim \sqrt{\frac{dK}{\alpha n}}
\end{align*}
and it remains to bound $\op{\Delta_\bE\suptt}$. Note that
\begin{align*}
	\op{\Delta_\bE\suptt}&=\op{\frac{1}{n_a\suptt}\sum_{i=1}^n\ind{\hat s_i^{(t-1)}=a}\bE_i-\frac{1}{n_a^\ast}\sum_{i=1}^n\ind{ s_i^{\ast}=a}\bE_i}\\
	&\le \op{\frac{1}{n_a\suptt}\sum_{i=1}^n\left[\ind{\hat s_i^{(t-1)}=a}-\ind{ s_i^{\ast}=a}\right]\bE_i}+\op{\frac{n_a^\ast-n_a\suptt}{n_a\suptt n_a^\ast}\sum_{i=1}^n\ind{ s_i^{\ast}=a}\bE_i}\\
	&\le \op{\frac{1}{n_a\suptt}\sum_{i=1}^n\ind{\hat s_i^{(t-1)}=a, s_i^{\ast}\ne a}\bE_i}+\op{\frac{1}{n_a\suptt}\sum_{i=1}^n\ind{ \hat s_i^{(t-1)}\ne a, s_i^{\ast}=a}\bE_i}\\
	&+ \frac{1}{n_a\suptt}\cdot \left|\sum_{i=1}^n\ind{ s_i^{\ast}=a, \hat s_i^{(t-1)}\ne a}\right|\op{\frac{1}{n_a^\ast}\sum_{i=1}^n\ind{ s_i^{\ast}=a}\bE_i}\\
	&+ \frac{1}{n_a\suptt}\cdot \left|\sum_{i=1}^n\ind{ s_i^{\ast}\ne a, \hat s_i^{(t-1)}=a}\right|\op{\frac{1}{n_a^\ast}\sum_{i=1}^n\ind{ s_i^{\ast}=a}\bE_i}\\
	&\overset{(a)}{\lesssim} \frac{K\sqrt{(d+n)h_a(\hat\bs^{(t-1)},\bs^\ast)}}{\alpha n}+\frac{K}{n}h_a(\hat\bs^{(t-1)},\bs^\ast)\sqrt{\frac{dK}{\alpha n}}\\
	&\overset{(b)}{\lesssim} \frac{K\sqrt{(d+n)h_a(\hat\bs^{(t-1)},\bs^\ast)}}{\alpha n}
\end{align*}
where in (a) we've used the fact that $\calQ_2$ holds and (b) is due to that fact that $h_a(\hat\bs^{(t-1)},\bs^\ast)\lesssim \alpha n/K$.
\subsection{Proof of Lemma \ref{lem:pertub-formula}}
The conclusion directly follows from dilation, i.e., define 
\begin{align*}
	\bX^*:=\mat{\bf0&\bX\\\bX^\top&\bf0},\quad \bM^*:=\mat{\bf0&\bM\\\bM^\top&\bf0},\quad\Delta^*:=\mat{\bf0&\Delta\\\Delta^\top&\bf0}
\end{align*}
and applying Theorem 1 in \cite{xia2021normal}.
\subsection{Proof of Lemma \ref{lem:EiXibound}}
To decouple the potential dependency of $\bE_i$ and $\bXi$, we employ the technical tool in \cite{mendelson2016upper}, for which we need to introduce additional  notations. Let $\calF\subset L_2$ be a class of function defined on some measure $\mu$. Denote $\EE \op{G}_\calF:=\EE\sup_{f\in\calF} G_f$ where $\{G_f:f\in\calF\}$ is the centered canonical gaussian process indexed by $\calF$. A class $\calF$ is $L$-subgaussian if for every $f,h\in\calF\cup\{0\}$, $\op{f-h}_{\psi_2}\le L\op{f-h}_{L_2}$. Here $\op{\cdot}_{\psi_2}$ is the standard $\psi_2$ norm (sub-Gaussian norm). The following lemma is adapted  from \cite{mendelson2016upper}.
\begin{lemma}[Theorem 1.13 in \cite{mendelson2016upper}]\label{lem:product-process}
	Let $\calF$ be a $L$-subgaussian class. There exists an absolute constant $c_0$ and for every $q>4$ there exists a constant $c_1(q)$ that depends only on $q$ for which the following holds. Let $\calF$ be a class of functions on $(\Omega,\mu)$, set $u\ge \max\{8,\sqrt{q}\}$ and consider an integer $s_0\ge 0$. Then, with probability at least $1-2\exp(-c_02^{s_0}u^2)$, for every $f\in \calF$,
	\begin{align*}
		\left|\sum_{i=1}^n(f^2(X_i)-\EE f^2)\right|\le c_1(q)\left(u^2\tilde\Lambda^2_{s_0,u}(\calF)+u\sqrt{n}\left(d_q(\calF)\tilde\Lambda_{s_0,u}(\calF)\right)\right)
	\end{align*} 
	where  $\tilde\Lambda_{s_0,u}(\calF)$ (see a formal definition in \cite{mendelson2016upper}) can be further bounded by
	\begin{align*}
		\tilde\Lambda_{s_0,u}(\calF)\le c_2L\left(\EE\op{G}_{\calF}+2^{s_0/2}d_q(\calF)\right)
	\end{align*}
	and $d_p(\calF):=\sup_{f\in \calF}\op{f}_{L_p}$ for any $p>0$.
\end{lemma}

In our case, denote $\mu$ as the distribution of each $\bE_i$. Define $$\calX_r =\left\{\bX\in \RR^{d_1\times d_2}, \text{rank}(\bX)\le r,\op{\bX}\le 1\right\}$$
and 
 $\calF_r:=\left\{f:f(\cdot)=\inp{\cdot}{\bX},\bX\in \calX_r\right \}$ on $\mu $. Observe that for any $f,g\in \calF_r$ and any $\bE\in\RR^{d_1\times d_2}\sim \mu $ having the same distribution as $\bE_i$, 
\begin{align*}
	\op{f(\bE)-g(\bE)}_{\psi_2}=\op{\inp{\bE}{\bX_1-\bX_2}}_{\psi_2}\lesssim\fro{\bX_1-\bX_2}=\op{\inp{\bE}{\bX_1-\bX_2}}_{L_2}
\end{align*}
This indicates that $\calF_r$ is $L$-subgaussian class with $L\le C$ for some absolute constant $C>0$. Also notice that for any fixed $q\ge 2$ 
\begin{align*}
	d_q(\calF_r)=\sup_{f\in \calF_r}\op{f}_{L_q}=\sup_{\bX\in \calX_r}\op{\inp{\bE}{\bX}}_{L_q}\lesssim\sup_{\bX\in \calX_r}\fro{\bX}\le \sqrt{r}
\end{align*}
where we've used the  the moment characterization of the $\psi_{2}$ norm, and that
\begin{align*}
	\EE\op{G}_{\calF_r}=\EE \sup_{f\in \calF_r}G_f=\EE\sup_{\bX\in \calX_r}|\inp{\bZ}{\bX}|\le r\EE \op{\bZ}\lesssim r\sqrt{d}
\end{align*}
where $\bZ\in\RR^{d_1\times d_2}$ has i.i.d. standard normal entries.
As a result, by choosing $s_0$ such that $2^{s_0}\asymp d$, $q=5$, $u=8$, we can  apply Lemma \ref{lem:product-process} and obtain that with probability at least $1-\exp(-cd)$,
\begin{align*}
	\sup_{\bX\in \calX_r}\sum_{i:s_i^\ast=b}&\left(\inp{\bE_i}{\bX}^2-\fro{\bX}^2\right)=\sup_{f\in \calF_r}\sum_{i:s_i^\ast=b}\left(f^2(\bE_i)-\EE f^2(\bE_i)\right)\\
	&\lesssim r\left(dr+\sqrt{dr\cdot n_b^\ast}\right)
\end{align*} 
Hence with the same probability,
\begin{align*}
	\sup_{\substack{\bXi\in\RR^{d_1\times d_2}, {\rm rank}(\bXi)\leq r\\ \|\bXi\|\leq 1}} &\sum_{i=1}^n\ind{s_i^{\ast}=b}\inp{\bE_i}{ \bXi}^2=\sup_{\bX\in \calX_r} \sum_{i:s_i^\ast=b}\left(\inp{\bE_i}{\bX}^2-\fro{\bX}^2\right) +n_b^\ast\sup_{\bX\in \calX_r}\fro{\bX}^2\\
	&\lesssim r\left(dr+\sqrt{dr\cdot n_b^\ast}+n_b^\ast\right)\lesssim r\left(dr+n_b^\ast\right)
\end{align*} 
by noticing that $\fro{\bX}^2\le r$ for any $\bX\in\calX_r$.

\subsection{Proof of Lemma~\ref{lem:init-two-facts}}
We first prove (I). By definition of k-means
\begin{align*}
	\fro{\bfrakM-\bG}&\le \fro{\bfrakM-\hat \bG}+\fro{\hat \bG-\bG}\le2 \fro{\hat \bG-\bG}\le 2 \sqrt{2K}\op{\hat \bG-\bG}
\end{align*}
It suffices to notice that 
\begin{align*}
	\op{\hat \bG-\bG}&=\op{\scrM_3(\bcalX\times_1\hat \bU\hat \bU^\top \times_2\hat\bV\hat \bV^\top -\bcalM\times_1\hat \bU\hat \bU^\top \times_2\hat\bV\hat \bV^\top)}\\
	&=\op{\scrM_3(\bcalE)(\hat\bV\hat \bV^\top\otimes \hat \bU\hat \bU^\top)}= \op{\scrM_3(\bcalE)(\hat\bV\otimes \hat \bU )}\\
	&\le C\left(\sqrt{d(r_\bU+r_\bV)}+\sqrt{n}\right)
\end{align*}
where the last inequality holds with probability at least $1-\exp(-cd)$ by Lemma 5 in \cite{zhang2018tensor}. Hence there exists some $C_0>0$, and with  probability at least $1-\exp(-cd)$ we have
\begin{align*}
	\fro{\bfrakM-\bG}&\le C_0\sqrt{K}\left(\sqrt{dKr+n}\right)
\end{align*}
for some absolute constant $C_0>0$. 
\\It remains to prove (II). By definition of $\bcalG$, we obtain
\begin{align*}
&\fro{\bcalG\times_3 (\be_i^\top-\be_j^\top)}\\
&=\fro{\left[\bcalM\times_3 (\be_i^\top-\be_j^\top)\right]\times_1\hat \bU\hat \bU^\top \times_2\hat\bV\hat \bV^\top}\\
&\ge \fro{\left[\bcalM\times_3 (\be_i^\top-\be_j^\top)\right]\times_1 \bU \bU^\top \times_2\bV \bV^\top}-\fro{\left[\bcalM\times_3 (\be_i^\top-\be_j^\top)\right]\times_1(\hat \bU\hat \bU^\top-\bU \bU^\top )\times_2\bV \bV^\top}\\
&-\fro{\left[\bcalM\times_3 (\be_i^\top-\be_j^\top)\right]\times_1\hat \bU\hat \bU^\top \times_2(\hat\bV\hat \bV^\top-\bV \bV^\top)}\\
&\ge \Delta-\frac{\Delta}{4}-\frac{\Delta}{4}\ge \frac{\Delta}{2}
\end{align*}
where we've used the fact that $\calQ_0$ holds and the equivalence between $\sqrt{2}\fro{\sin\Theta(\bU_1,\bU_2)}$ and projection distance 	$\fro{\bU_1\bU_1^\top-\bU_2\bU_2^\top}$.

\subsection{Proof of Lemma~\ref{lem:relaxed-two-facts}}
The proof of (I) is identical to that in the proof of Lemma \ref{lem:init-two-facts} and hence we only show (II).  By definition of $\bcalG$, we obtain
\begin{align*}
&\fro{\bcalG\times_3 (\be_i^\top-\be_j^\top)}\\
&=\fro{\left[\bcalM\times_3 (\be_i^\top-\be_j^\top)\right]\times_1\hat \bU_1\hat \bU_1^\top \times_2\hat\bV_1\hat \bV_1^\top}\\
&\ge \fro{\left[\bcalM_1\times_3 (\be_i^\top-\be_j^\top)\right]\times_1 \bU_1 \bU_1^\top \times_2\bV_1 \bV_1^\top}\\
&-\fro{\left[\bcalM_1\times_3 (\be_i^\top-\be_j^\top)\right]\times_1(\hat \bU_1\hat \bU_1^\top-\bU_1 \bU_1^\top )\times_2\bV_1 \bV_1^\top}\\
&-\fro{\left[\bcalM_1\times_3 (\be_i^\top-\be_j^\top)\right]\times_1\hat \bU_1\hat \bU_1^\top \times_2(\hat\bV_1\hat \bV_1^\top-\bV_1 \bV_1^\top)}\\
&-\fro{\left[\bcalM_2\times_3 (\be_i^\top-\be_j^\top)\right]\times_1 \hat\bU_1 \hat\bU_1^\top \times_2\hat\bV_1 \hat\bV_1^\top}\\
&\overset{(a)}{\ge} \fro{\bM_1}-\frac{\fro{\bM_1}}{6}-\frac{\fro{\bM_1}}{6}-\fro{\bM_2}\\
&\overset{(b)}{\ge} \frac{5}{9}\fro{\bM_1}\ge \frac{\Delta}{2}
\end{align*}
where in (a) we've used \eqref{eq:relaxed-uv-small}, (b) and (c) are due to the facts that $\fro{\bM_1}\ge 9\fro{\bM_2}$ and $\Delta=\fro{\bM_1-\bM_2}\le 10/9\cdot \fro{\bM_1}$, by properly choosing the absolute constant $C$ in Assumption \ref{assump:weak-SNR} and the proof is therefore completed.


%
%
%

\end{document}